\documentclass[12pt]{amsart}
\usepackage{amsmath,amssymb,amscd,array, mathrsfs }

\usepackage{amsmath,amscd,amssymb,amsthm,array}
\usepackage{color}

\usepackage{amsmath,amssymb,mathrsfs,amsthm, tikz-cd,mathrsfs}
\usepackage{comment}
\def\url#1{\expandafter\s

\tring\csname #1\endcsname}

\def\mmat #1,#2,#3,#4,{\text{\small\arraycolsep=3pt $
\begin{pmatrix}#1&#2\\#3&#4\end{pmatrix}$}}

\usepackage{hyperref}
\usepackage{capt-of}

\usepackage{multirow}
\usepackage[all]{xy}

\usepackage{comments}
\newComments\SBe{Said}{blue}
\newComments\SBo{Sofiane}{blue}
\newComments\AM{Nacer}{blue}
\newComments\DL{DL}{red}
\newComments\QEh{QEh}{blue}

\def\mmat #1,#2,#3,#4,{\text{\small\arraycolsep=3pt $
\begin{pmatrix}#1&#2\\#3&#4\end{pmatrix}$}}

\usepackage{lscape}
\usepackage{tikz-cd}
\usepackage{enumerate}

\usepackage{DLdef1}
\usepackage{multicol}

\def\mmat #1,#2,#3,#4,{\text{\small\arraycolsep=3pt $
\begin{pmatrix}#1&#2\\#3&#4\end{pmatrix}$}}

\renewcommand {\ssbegin}[2][*]
 {\refstepcounter{subsection}%
\if#1*
\addcontentsline{toc}{subsection}{\thesubsection.~#2}%
\else
\addcontentsline{toc}{subsection}{\thesubsection.\hskip 1pc #2. #1}%
\fi
 \def \secno {\gdef \secno {}{\ssecfont
\thesubsection.\hskip 2ex}%
 }%
 \begin{#2}}

\renewcommand {\sssbegin}[2][*]
 {\refstepcounter{subsubsection}
\if#1*
\addcontentsline{toc}{subsubsection}{\thesubsubsection.\hskip 1pc #2}%
\else
\addcontentsline{toc}{subsubsection}{\thesubsubsection.\hskip 1pc #2. #1}
\fi
 \def \secno {\gdef \secno {}{\ssecfont \thesubsubsection.\hskip 2ex}%
 }%
 \begin{#2}}

\renewcommand {\parbegin}[2][*]
 {\refstepcounter{paragraph}
\if#1*
\addcontentsline{toc}{paragraph}{\theparagraph.\hskip 1pc #2}%
\else
\addcontentsline{toc}{paragraph}{\theparagraph.\hskip 1pc #2. #1}
\fi
 \def \secno {\gdef \secno {}{\ssecfont \theparagraph.\hskip 2ex}%
 }%
 \begin{#2}}

\rmnameii{etr}{etr}
\rmnameii{evv}{ev}
\DeclareMathOperator{\h}{\mathcal{H}}

\DeclareMathOperator{\K}{\mathbb{K}}
\newcommand{\Z}{\mathbb{Z}}

\newcommand {\A}{{\cal{A}}}
\newcommand {\w}{\omega}

\newcommand{\Ll}{{\mathrm{L}}}
\newcommand{\Rr}{{\mathrm{R}}}

    \newcommand{\Om}{\Omega}

\newcommand{\al}{\alpha}

\newcommand{\la}{\lambda}
\newcommand{\De}{\Delta}

\newcommand{\om}{\omega}
\newcommand{\prs}{\langle\;,\;\rangle}
\def\br{[\;,\;]}

\newcommand{\g}{\mathfrak{g}}
\begin{document}

\title[Flat quasi-Frobenius Lie superalgebras]{Flat quasi-Frobenius Lie superalgebras}


\author{Sofiane Bouarroudj}

\address{Division of Science and Mathematics, New York University Abu Dhabi, P.O. Box 129188, Abu Dhabi, United Arab Emirates.}
\email{sofiane.bouarroudj@nyu.edu}

\author{Hamza El Ouali}
\address{Universit\'e Cadi-Ayyad,
	Facult\'e des sciences et techniques,
	B.P. 549, Marrakech, Maroc.}
\email{eloualihamza11@gmail.com}




\begin{abstract} A non-associative  superalgebra is called {\it pre-symplectic} if it is equipped with a non-degenerate, anti-symmetric bilinear form. It is called {\it quasi-Frobenius} if, in addition, is a Lie superalgebra and the form is closed. We introduce the Levi-Civita product associated with pre-symplectic superalgebras and establish its existence and uniqueness. We then introduce the symplectic product associated with quasi-Frobenius Lie superalgebras. We prove that while such a product always exists, it is not unique. We therefore define a \emph{natural symplectic product} that depends only on the Lie structure and the bilinear form. When the curvature of this product vanishes, the superalgebra is called a \emph{flat quasi-Frobenius Lie superalgebra}.  

In this paper, we study flat quasi-Frobenius Lie superalgebras and introduce the notion of a \emph{flat double extension}. We prove that the double extension process characterizes such  superalgebras. More precisely, every flat orthosymplectic (resp. periplectic) quasi-Frobenius Lie superalgebra can be obtained by a sequence of flat double extensions starting from an abelian one (resp. the trivial one). Moreover, we show that every flat quasi-Frobenius Lie superalgebra is nilpotent with a degenerate center. We apply our results to obtain a complete classification of such  superalgebras of total dimension at most five.

\end{abstract}


\maketitle

\thispagestyle{empty}
\setcounter{tocdepth}{2}

\section{Introduction}
\subsection{The Riemannian setting}
A \emph{pseudo-Riemannian Lie group} is a Lie group $G$ endowed with a left-invariant pseudo-Riemannian metric $\mu$.  
Its Lie algebra $\mathfrak{g}=T_eG$, equipped with the scalar product $\prs=\mu_e$, is called a \emph{pseudo-Riemannian Lie algebra}.  
Whenever there exists a torsion-free connection compatible with the metric, such a connection is necessarily unique and is called the \emph{Levi-Civita connection} (see for instance \cite{NOM, ON}). The Levi-Civita connection induces a bilinear product $(u,v)\mapsto u\star v$ on $\mathfrak{g}$, called the \emph{Levi-Civita product}, defined by Koszul’s formula:
\[
2\langle u\star v, w\rangle
= \langle [u,v], w\rangle
+ \langle [w,u], v\rangle
+ \langle [w,v], u\rangle,
\qquad \text{for all }u,v,w\in\mathfrak{g}.
\]
This product satisfies the following properties: 
\[
u\star v - v\star u = [u,v],\qquad
\langle u\star v, w\rangle = -\,\langle v, u\star w\rangle,
\qquad \text{for all }u,v,w\in\mathfrak{g}.
\]
In the setting of pseudo-Euclidean non-associative superalgebras (in particular, Leibniz superalgebras)  this construction has been generalized in \cite{BBE2}.

The curvature of $(\mathfrak{g},[\; , \;], \prs)$ at the identity is
\[
K(u,v):=\Ll_{[u,v]}-[\Ll_u,\Ll_v],\qquad \text{for all }u,v\in\mathfrak{g},
\]
where $\Ll_u$ is the left multiplication defined by
$
\Ll_u(v)=u\star v,$ for all $u,v\in\mathfrak{g}.$

A pseudo-Riemannian Lie group $(G,\mu)$ is called \emph{flat} when its curvature vanishes identically; in this case, $(\mathfrak{g}, [\; , \;], \prs)$ will be referred to as a \emph{flat pseudo-Riemannian Lie algebra}. In fact, flatness is equivalent to 
$(\g, \star)$ being a left-symmetric algebra.

The study of flat Lie algebras originated with the Euclidean case. A foundational result was proved by Milnor in 1976 \cite{M}, who characterized such algebras: a Lie algebra $\g$ is flat Euclidean if and only if it admits an orthogonal decomposition
$
\g
=
\mathfrak{b}
\oplus 
\mathfrak{u}
$
where 
$\mathfrak u$ is an abelian ideal, 
$\mathfrak b$ is an abelian subalgebra, and for every 
$b\in \mathfrak b$, the adjoint map 
$\ad_{\mathfrak b}$ is skew-symmetric with respect to the inner product. 

In recent years, a number of papers have investigated flat pseudo-Riemannian Lie algebras with various signatures, using in particular the method of double extension introduced in \cite{AM}.  
Further developments include the study of flat Lorentzian Lie algebras \cite{AM, ABL, BL} as well as the study of flat pseudo-Euclidean nilpotent Lie algebras with the signature $(2,n-2)$, see \cite{BL1}. Furthermore, in \cite{BBE2}, the authors studied flat pseudo-Euclidean left-Leibniz superalgebras.
\subsection{The symplectic setting.}
A \emph{symplectic Lie group} is a Lie group $G$ endowed with a left-invariant symplectic form $\Omega$.  
The associated Lie algebra $\mathfrak{g}=T_eG$, equipped with $\omega=\Omega_e$, is a \emph{symplectic Lie algebra} (called quasi-Frobenuis in this paper).  A connection $\nabla$ on $G$ is called a \emph{symplectic connection} if it is torsion-free and is compatible with $\Omega$, i.e., $\nabla \Om=0$.   
Such a symplectic connection always exists, but it is  not unique (see \cite{BCGRS}).  Whenever a symplectic connection is given, one may define a bilinear product 
$
(u,v) \mapsto u \bullet v
$
on $\mathfrak{g}$ that satisfies the following properties: 
\[
u\bullet v - v\bullet u = [u,v],\qquad
\omega(u\bullet v, w) = -\,\omega(v, u\bullet w),
\qquad \text{for all }u,v,w\in\mathfrak{g}.
\]
In \cite{BB}, the authors introduced a {\it natural} symplectic connection on symplectic Lie groups, defined by
\[
\omega(u\star v, w)
= \tfrac{1}{3}\,\omega([u,v],w)
+ \tfrac{1}{3}\,\omega([u,w],v),
\qquad \text{for all }u,v,w\in\mathfrak{g}.
\]
As in the pseudo-Riemannian case, one defines the curvature of this connection. A symplectic Lie group $(G,\Omega)$ is called \emph{flat} when its curvature vanishes identically.  
In this case, the induced Lie algebra $(\mathfrak{g}, [\; , \;], \omega)$ is a \emph{flat quasi-Frobenius Lie algebra}. In fact, flatness is equivalent to the fact that the symplectic product defines a left-symmetric structure on $\g$. 

Despite its importance, the flat symplectic case has received little attention from experts. General results on flat quasi-Frobenius Lie algebras appear in \cite{V}. Meanwhile, the symplectic connection above was studied in \cite{BELL}; further contributions can be found in \cite{A, BC, MR1, NB}.  

The study of such structures is motivated by a well-known open problem of Milnor \cite{M}, which seeks a classification and description of flat affine Lie groups, i.e. Lie groups carrying a left-invariant flat affine structure.

Recall that if $(G,\Omega)$ is a symplectic Lie group, then the left-invariant connection $\nabla$ defined by
\begin{equation}\label{eq:flat-conn}
\omega(\nabla_u v, w) = -\omega(v,[u,w]), \qquad \text{for all }u,v,w\in\mathfrak{g},
\end{equation}
is indeed flat, i.e. its curvature vanishes.  Note, in addition, that $\nabla$ is symplectic (that is, $\nabla\Omega = 0$) if and only if $G$ is abelian, equivalently its Lie algebra $\mathfrak{g}$ is abelian (see \cite{Ch}).
\subsection{Outline of the paper}
Section \ref{back} collects the preliminary material required for our work. For completeness and to establish notation, we include a review of key concepts, some of which may be familiar to specialists.

In Section \ref{Levi-section}, we introduce the Levi-Civita product on a class of non-associative anti-commutative superalgebras that we call pre-symplectic. These superalgebras admit a  non-degenerate  anti-symmetric bilinear form $\omega$ (not necessarily closed). 
This product $\bullet$ is the unique bilinear map on the pre-symplectic superalgebra $\A$ that satisfies the following properties: 
\[
[u,v] = u\bullet v - (-1)^{|u||v|} v\bullet u,
\qquad
\omega(u\bullet v,w) = (-1)^{|u||v|}\,\omega(v, u\bullet w),
\]
for all homogeneous $u,v,w\in\A$.  
The product $\bullet$ is also characterized by a Koszul formula (see Eq. \eqref{Koszul}). 
Its existence and uniqueness are established in Prop~\ref{LVCT}.

Further, we study the Levi-Civita product in the case where the pre-symplectic superalgebra is a Lie superalgebra and the form $\omega$ is closed. We call these superalgebras \emph{quasi-Frobenius Lie superalgebras}, following \cite{BBE,BE,BM}. It is worth mentioning that several authors call them \emph{symplectic Lie superalgebras}
(see \cite{BC, BEB, BBNS} to cite a few). We show that the Levi-Civita product associated with a quasi-Frobenius Lie superalgebra $(\g, \omega)$ is always
right-symmetric (see Prop. \ref{al}) and satisfies
\[
\omega(u\bullet v,w) = \omega(u,[v,w]),
\qquad \text{ for all homogeneous }\, u,v,w\in\g.
\]
The opposite product given by
$
u\circ v = -(-1)^{|u||v|} v\bullet u,
$ for all homogeneous $u,v\in \g$,
defines a left-symmetric structure on $\fg$, which superize the classical construction in Eq.\eqref{eq:flat-conn}.

\medskip In addition, we introduce a second product $\star$ on a quasi-Frobenius Lie superalgebra $(\g, [\; , \;], \omega)$, defined by
\[
u\star v - (-1)^{|u||v|} v\star u = [u,v],
\qquad
\omega(u\star v,w) 
= -(-1)^{|v||w|}\,\omega(v, u\star w).
\]
Such a symplectic product exists but is not unique (see Corollary~\ref{No}).  
This contrasts with the uniqueness of the Levi-Civita product on pre-symplectic superalgebras as well as  in the case of pseudo-Euclidean superalgebras (see \cite{BBE2}).

We show nevertheless that in any  quasi-Frobenius Lie superalgebra $(\g, [\; , \;], \omega)$ we can construct a unique symplectic product with the key property that it  depends only on the Lie bracket and the symplectic form $\omega$. This product is explicitly given by:
\begin{equation*}\label{tti}
\omega(u\star v,w)
=
\frac{1}{3}\left(
\omega([u,v],w)
+ (-1)^{|v||w|}\,\omega([u,w],v)
\right),
 \text{ for all homogeneous }\, u,v,w\in\g.
\end{equation*}
We call this product the {\it natural symplectic product} associated with $(\g,\br,\omega)$. When the curvature of the natural  symplectic product $\star$ is zero,  the product is called \emph{flat}, and the triple $(\g, \br, \omega)$ is referred to as a \emph{flat quasi-Frobenius Lie superalgebra}. In fact,  flatness is equivalent to the fact that the product defines a left-symmetric structure on the underlying vector space of $\fg$.

In Section \ref{Sol-nil-sec}, we initiate the study of the structure theory for flat quasi-Frobenius Lie superalgebras. We show that such  algebras are solvable when the form  $\omega$ is even (Theorem \ref{solvable-even}) and nilpotent when the form $\omega$ is odd (Theorem \ref{nilpotent-odd}).

In Section \ref{d-ext}, we introduce the notion of \emph{flat double extension} for flat quasi-Frobenius Lie superalgebras, see Theorem \ref{double-ex1}, Theorem \ref{double-ex2}, Theorem \ref{double-ex3}, and Theorem \ref{double-ex4}. This construction is inspired by the classical double extension, originally introduced by Medina and Revoy \cite{MR} for quadratic Lie algebras, and its symplectic analogs \cite{MR2, DM}. It is worth noticing that this notion differs from the notion of double extensions of quasi-Frobenius Lie superalgebras considered in \cite{BM, BE}.  

In Section \ref{flat-nil-sec}, we show that all non-abelian flat quasi-Frobenius Lie superalgebras have a degenerate center, see Theorem \ref{Centre}. Moreover, every such Lie superalgebra can be obtained via flat double extensions. In particular, every flat orthosymplectic quasi-Frobenius Lie superalgebra can be constructed by a sequence of flat double extensions starting from an abelian Lie superalgebra (see Corollary \ref{start-abel}), while flat periplectic quasi-Frobenius Lie superalgebras can be obtained by a sequence of flat double extensions starting from the trivial algebra $\{0\}$ (see Corollary \ref{ext-triv}). Furthermore, we show that every flat orthosymplectic quasi-Frobenius Lie superalgebra is not only solvable but nilpotent (Theorem \ref{ref-nilp}).

Section \ref{classification} is devoted to classifying all flat quasi-Frobenius Lie superalgebras of total dimension at most five over an algebraically closed field of characteristic zero.

\section{Backgrounds}\label{back}
Let $\mathbb{K}$ be an arbitrary field of characteristic zero.  The group of integers modulo $2$ is denoted by $\Z_{2}$. 

Let $V=V_{\bar 0}\oplus V_{\bar 1}$ be a superspace defined over  ${\mathbb K}$. The parity of a homogeneous element $v\in V_{\bar{i}}$ is denoted by $|v|:=\bar{i}$. The element $v$ is called \textit{even} if $v\in V_{\bar 0}$ and \textit{odd} if $v\in V_{\bar 1}$. 

Throughout the text, all elements in  linear expressions are assumed to be homogeneous, unless otherwise stated. 

Let \(V\) be a vector superspace. For homogeneous elements \(x,y\in V\) we denote by \(x^*,y^*\in V^*\) their duals . For every \(x^*\otimes y^*\in V^*\otimes V^*\) and  \(u\otimes v\in V\otimes V\), we define the pairing
\[
\langle x^*\otimes y^*, u \otimes v\rangle :=(-1)^{|y||u|}x^*(u)\,y^*(v).
\]
The wedge (graded antisymmetrisation) is given by
\[
x^*\wedge y^*=x^*\otimes y^* - (-1)^{|x||y|}\,y^*\otimes x^*.
\]

A linear map $\varphi:V\rightarrow W$ between two superspaces is called \textit{even} if $\varphi(V_{\bar{i}})\subset W_{\bar{i}}$ and \textit{odd} if $\varphi(V_{\bar{i}})\subset W_{\overline{i}+\overline{1}}$.

Let $V$ and $W$ be two superspaces. Let ${\cal B}\in \text{Bil}(V,W)$ be a homogeneous bilinear form. Recall that the Gram matrix $B=(B_{ij})$ associated to $\cal B$ is given by the following formula:
\begin{equation*}\label{martBil}
B_{ij}=(-1)^{p({\cal B})p(v_i)}{\cal B}(v_{i}, w_{j})\text{~~for the basis vectors $v_{i}\in V$ and $w_{j}\in W$ .}
\end{equation*}
This definition allows us to identify a~bilinear form $B(V, W)$ with an element of $\mathrm{Hom}(V, W^*)$. 
Consider the \textit{upsetting} of bilinear forms
$u:\text{Bil} (V, W)\rightarrow \text{Bil}(W, V)$ given by the formula 
\begin{equation*}\label{susyB}
u({\cal B})(w, v)=(-1)^{p(v)p(w)}{\cal B}(v,w)\text{~~for any $v \in V$ and $w\in W$.}
\end{equation*}
In terms of the Gram matrix $B$ of ${\cal B}$, we have 
\begin{equation*}
u(B)=
\left ( \begin{array}{cc}  R^{t} & (-1)^{p(B)}T^{t} \\ (-1)^{p(B)}S^{t} & -U^{t} \end{array} \right ),
\text{ for $B=\left (  \begin{array}{cc}R & S \\ T & U \end{array} \right )$}.
\end{equation*}
From now on, we will assume that $V=W$. Following \cite{L}, we say that: 
\begin{itemize} \item the form
$\cal B$ is  \textit{symmetric} if  and only if $u(B)=B$;

\item the form $\cal B$ is \textit{anti-symmetric} if and only if $u(B)=-B$.

\end{itemize}
Given a linear map $f: V\rightarrow V$. The adjoint of $f$ with respect to ${\cal B}$ is the linear map $f^*$ defined as follows:
\[
{\cal B}(f(v),w)= (-1)^{|f||v|}{\cal B}(v,f^*(w)).
\]
We say that $f$ is ${\cal B}$-symmetric (resp. ${\cal B}$-antisymmetric) if $f^*=f$ (resp. $f^*=-f$).

A  sub-superspace $U\subseteq V$ is called degenerate if and only if $U\cap U^\perp\not = \{0\}$, where $\perp$ stands for the orthogonal with respect to the bilinear form ${\cal B}$.

Throughout this paper, $\omega$ will denote an anti-symmetric bilinear form, which by definition satisfies
\[
\omega(u,v)=-(-1)^{|u||v|} \omega(v,u), \text{ for all homogeneous } u,v\in V.
\]
A non-degenerate, anti-symmetric even form is called {\it orthosymplectic}; an odd one is called {\it periplectic}.

\ssbegin{Lemma}[See \cite{BM}]\label{BM}
Let $V = V_{\bar{0}} \oplus V_{\bar{1}}$ be a finite-dimensional superspace with $V_{\bar{1}} \neq \{0\}$, endowed with a periplectic  form $\omega$. Then  $\dim(V_{\bar{0}}) = \dim(V_{\bar{1}})$.
\end{Lemma}

A superalgebra $(\A, \bullet)$  over a field $\mathbb{K}$ consists of a superspace $\A=\A_{\bar 0}\oplus \A_{\bar 1}$ and a binary operation satisfying: $\A_\alpha \bullet \A_\beta \subseteq \A_{\alpha+\beta}, \forall \alpha, \beta \in \mathbb{Z}_2$. 

 For every $u \in \A$, we denote the left and right multiplication operators by $\Ll_u$ and $\Rr_u$, respectively. More precisely, $\Ll_u$ and $\Rr_u$ are defined as  follows: 
 \[
 \Ll^\bullet_u(v):=u \bullet v \text{  and }\Rr^\bullet_u(v):=(-1)^{|u||v|}v\bullet u, \text{ for all $v \in \A$}.
 \]
Let us define:
$$ N_\ell(\A, \bullet) := \left\{ u \in \A, \mathrm{L}^\bullet_u = 0 \right\}, \;  N_r(\A, \bullet) := \left\{ u \in \A, \mathrm{R}_u^\bullet = 0 \right\}, \; \text{and} \; N(\A, \bullet):= N_\ell(\A, \bullet)\cap  N_r(\A, \bullet). $$
We refer to these as the {\it left normalizer}, {\it right normalizer}, and (two-sided) {\it normalizer}, respectively.

A superalgebra $(\cal A, [\; , \; ])$ is said to be a {\it Lie superalgebra} if the following conditions hold for all homogeneous elements 
$x,y,z \in {\cal A}:$ 
\begin{itemize}
\item (Super anti-commutativity): 
$$
[x,y]=-(-1)^{|x||y|}[y,x].
$$
\item (Super Jacobi identity):
\[
[[x,y],z]=[x,[y,z]] 
- (-1)^{|x||y|}\,[y,[x,z]].
\]
\end{itemize}
We recall that the super Jacobi identity generalizes the classical Jacobi identity by incorporating the parity signs arising from the 
${\mathbb Z}_2$-grading. In particular, when all elements are even, it reduces to the usual Jacobi identity. The identity ensures that the adjoint action defines a representation of the Lie superalgebra.

Throughout this paper, a general superalgebra will be denoted by $
({\cal A},\bullet)$, whereas a Lie superalgebra will be denoted by 
$(\mathfrak{g},[\;, \;])$. For a thorough exposition of Lie superalgebras, we refer the reader to \cite{L,S}.

Let $(\A, \bullet)$ be a non-associative superalgebra. On the underlying vector superspace $\A$, we
define the corresponding commutator 
$$[u, v]_\bullet :=u\bullet v -(-1)^{|u||v|} v\bullet u,\; \text{ for all homogeneous }  u,v\in\A.$$ 
The algebra $(\A, [\; , \; ]_\bullet)$ will be denote by  $\A^-$. The algebra $(\A, \bullet)$ is called Lie-admissible
superalgebra, if $\A^-$ is a Lie superalgebra.  In this case,  $\A^-$ is called the sub-adjacent Lie superalgebra of $(\A, \bullet)$ and the latter is called a compatible (Lie-admissible) superalgebra structure on
the Lie superalgebra $\A^-$.

A superalgebra $(\A, \bullet)$ is called left-symmetric (resp. right-symmetric) if and only if 
 \[
\mathrm{as}(a,b,c) = (-1)^{|a||b|} \, \mathrm{as}(b,a,c), \quad (\text{resp. } \mathrm{as}(a,b,c) = (-1)^{|b||c|} \, \mathrm{as}(a,c,b)),
 \]
 where 
 $\mathrm{as}(a, b, c) =(a \bullet b) \bullet c - a \bullet (b \bullet c)$  is the associator. 
 
It is well-known that such superalgebras are Lie-admissible. 

\ssbegin{Example}
Consider the superalgebra ${\cal A}_n$   defined on the homogeneous basis
\(
(e_1,\dots,e_n \mid f_1,\dots,f_n),
\)
whose non-zero products are given by
\begin{align*}
e_1 \bullet e_i &= e_{i+1},& e_1 \bullet f_j &= f_{j+1}, 
& 2 \le i \le n-1,\, 1 \le j \le n-1, \\
e_1\bullet e_1&= (-1)^{n+1}e_2, & e_i \bullet f_{\,n-i} &
= 
f_{\,n-i} \bullet e_i 
=
(-1)^{n-i}\, f_n,
& 2 \le i \le n-1,
\end{align*}
A straightforward computation shows that $({\cal A}_{n},\bullet)$ is a left-symmetric superalgebra.

The sub-adjacent Lie superalgebra ${\cal A}_n^-$ is given by 
\begin{align*}
[e_1, e_i] = e_{i+1}, 
\quad 
[e_1, f_i] = f_{i+1},\quad 
1 \le i \le n-1.
\end{align*}
\end{Example}

A  pre-symplectic superalgebra $( \A, \cdot, \om)$ is a  non-associative  superalgebra $(\A, \cdot)$ endowed with a
non-degenerate anti-symmetric  bilinear form $\om$.

\section{The Levi-Civita product and the Symplectic product on quasi-Frobenius Lie superalgebras}\label{Levi-section}

\subsection{The Levi-Civita product}
In \cite{BBE2}, the authors introduced the Levi-Civita product associated with a non-associative superalgebra equipped with a non-degenerate {\it symmetric} bilinear form.  Here, we introduce its analogue in the case where the form is  {\it anti-symmetric}. 
\sssbegin{Proposition}\label{LVCT} Let $(\A,\br,\om)$ be a pre-symplectic anti-commutative superalgebra. There exists a unique product $\bullet$ on $\A$ that satisfies \textup{(}for all  $u,v,w\in \A)$:
	\begin{eqnarray} \label{torsion}
		[u, v]_\bullet&=&[u,v], 
	\\ \label{compatibl} (\mathrm{L}_u^\bullet)^* & = &
    \mathrm{L}_u^\bullet.
    \end{eqnarray}
More precisely, the product $\bullet$ is defined by
	\begin{equation}\label{Koszul}
		2\om( u\bullet v,w)=\om( [u, v],w) -(-1)^{|u||v|+|w||u|}\om( [v, w],  u)+(-1)^{|v||w|+|u||w|}\om([w, u],  v).
	\end{equation}
\end{Proposition}
The product $\bullet$ is called the {\it Levi-Civita} product associated with  $(\A,\br,\om)$.	

\begin{proof}
 \textbf{(Uniqueness)}  
 Suppose that there exist two products $\bullet$ and $\bullet'$ satisfying conditions \eqref{torsion},  \eqref{compatibl}. For any  $u,v,w \in \g$, we have
 \[
  \begin{array}{lcl}
  \omega(u \bullet v, w)
      &=& \omega([u,v], w) + (-1)^{|u||v|} \omega(v \bullet u, w) \\[2mm]
      &= &\omega([u,v], w) - (-1)^{|u||v| + |w||u|} \omega(v \bullet w, u) \\[2mm]
      &=& \omega([u,v], w) - (-1)^{|u||v| + |w||u|} \omega([v,w], u)
        - (-1)^{|u||v| + |w||u| + |v||w|} \omega(w \bullet v, u) \\[2mm]
      &=& \omega([u,v], w) - (-1)^{|u||v| + |w||u|} \omega([v,w], u)
        + (-1)^{|v||w| + |u||w|} \omega(w \bullet u, v) \\[2mm]
      &=& \omega([u,v], w) - (-1)^{|u||v| + |w||u|} \omega([v,w], u)
         + (-1)^{|v||w| + |u||w|} \omega([w,u], v)\\[2mm]
         &&
        - \omega(u \bullet v, w).
\end{array}
\]
The same computation by means of $\bullet'$, implies that $
 2\,\omega(u \bullet v, w) = 2\,\omega(u \bullet' v, w).$ Since $\omega$ is non-degenerate, it follows that $\bullet = \bullet'$.
 
 \bigskip
 \textbf{(Existence)}  
 We define the product $\bullet$ by the following relation:
 \[
 2\,\omega(u \bullet v, w)
     = \omega([u,v], w)
       - (-1)^{|u||v| + |w||u|}\omega([v,w], u)
       + (-1)^{|v||w| + |u||w|}\omega([w,u], v),
 \]
for all $u,v,w \in \A$. It is straightforward to verify that this product satisfies
 \[
 \omega(\mathrm{L}_u^\bullet (v), w)
     = -(-1)^{|u||v|} \omega(v, \mathrm{L}_u^\bullet( w)), \text{ and } 
 [u,v]_\bullet = [u,v] .
 \]
 Hence, the product $\bullet$ exists and is unique.
\end{proof}

Let $(\g, \br, \om)$ be a pre-symplectic Lie superalgebra, and let $\bullet$ denote its Levi-Civita product. 


Eqns. \eqref{torsion} and \eqref{compatibl} imply that the multiplication operator $\mathrm{L}^\bullet_u$, is $\omega$-symmetric, for every $u \in \g$, and that $\ad_u = \mathrm{L}^\bullet_u - \mathrm{R}^\bullet_u$, where $\ad_u: \g \to \g$ is the adjoint operator. Let us denote by $Z(\g)$ the center of $\fg$. 
The following relations hold (where $\g \bullet \g = \{ u \bullet v \mid u, v \in \g \}$):
\begin{equation}
[\g, \g]^\perp = \{u \in \g, \mathrm{R}_u^\bullet = -(\mathrm{R}_u^{\bullet})^*\} \quad \text{and} \quad (\g \bullet \g)^\perp = N_r(\g, \bullet), \text{ also } N(\g, \bullet)\subseteq Z(\g).
\label{eq1}\end{equation} 
\sssbegin{Definition}
Let $(\mathfrak{g}, \cdot, \omega)$ be a pre-symplectic superalgebra. The form $\omega$  is said to be closed if, for all elements $u, v, w \in \g$, the following identity holds:
\begin{equation}
(-1)^{|u||w|}\omega( u \cdot v, w ) + (-1)^{|v||u|} \omega( v \cdot w, u )+  (-1)^{|w||v|}\omega( w \cdot u, v ) = 0.
\label{cyclic}
\end{equation}
\end{Definition}
Following \cite{BBE, BE, BM}, if a pre-symplectic {\bf Lie} superalgebra $(\g, [\; , \;], \omega)$ satisfies the condition \eqref{cyclic}, then it is called a {\it quasi-Frobenius} Lie superalgebra. Cohomologically speaking, this condition means that $\omega \in Z^2(\fg; \K)$. In the case where $\omega$ is a coboundary,  we call $(\fg, [\;,\;],\omega)$ a {\it Frobenius} Lie superalgebra, see \cite{E,O}. Furthermore, if $\omega$ is orthosymplectic (resp. periplectic), we call $(\g, \br, \omega)$ an orthosymplectic (resp. periplectic) quasi-Frobenuis Lie superalgebra.
\sssbegin{Proposition}\label{pr3} 
Let $(\g, \br, \om)$ be a quasi-Frobenius Lie superalgebra and $\bullet$ be its Levi-Civita product. The following properties  hold.
\begin{enumerate}
    \item[$(i)$] For all   $u, v, w \in \g$, 
    \begin{equation}
    \om( \mathrm{L}_u^\bullet(v), w ) = \om( u, \ad_v(w) ). \label{Lev-cy}
    \end{equation}
    \item[$(ii)$] We have 
    $$[\g, \g]^\perp = N_\ell(\g, \bullet) = \{ u \in \g, \ad_u^* = -\ad_u \}.$$
    \item[$(iii)$] We have $$(\g \bullet \g)^\perp \subseteq [\g,\g]^\perp, \text{ and } Z(\g) = (\g \bullet \g)^\perp = N_\ell(\g, \bullet) \cap N_r(\g, \bullet).$$
\end{enumerate}
In particular, $[\g, \g]^\perp$ is an abelian subsuperalgebra and $Z(\g) \subseteq [\g, \g]^\perp$.
\end{Proposition}

\begin{proof}
Part (i) follows from the fact the $\omega$ is closed and Eq. \eqref{Koszul}. 
    
Let us prove Part (ii). For all $u \in [\mathfrak{g}, \mathfrak{g}]^{\perp}$ and $v,w \in \mathfrak{g}$, Part (i) implies that
$$
\omega(\mathrm{L}_u^\bullet (v), w) =  \omega(u, [v,w]) = 0.
$$ Since $\omega$ is non-degenerate, it follows that $
\Ll^\bullet_u  = 0,
$ which shows that $ [\mathfrak{g}, \mathfrak{g}]^{\perp} \subseteq N_\ell(\g, \bullet).
$ The converse inclusion follows similarly. On the other hand, from Eq. \eqref{eq1} we know that $\Rr^\bullet_u$ is $\omega$-anti-symmetric; consequently, $\operatorname{ad}_u = -\Rr^\bullet_u$
is also $\om$-anti-symmetric. 

Let us prove Part (iii). Since $[\cdot , \cdot]=[\cdot, \cdot]_\bullet$, it is clear that $[\g, \g]\subseteq \g\bullet \g$ and hence $(\g\bullet \g)^\bot\subseteq [\g,\g]^\bot$. In addition, from Eq. \eqref{eq1} we obtain $\Rr^\bullet_u=0$, for all $u\in (\g\bullet \g)^\bot$. Now, we also have $\mathrm{L}^\bullet_u = 0$, for all $u\in (\g \bullet \g)^\perp$, as seen in Part (i). From Eq. \eqref{eq1} we deduce $(\g \bullet \g)^\perp \subseteq N_\ell(\g, \bullet) \cap N_r(\g, \bullet)$.  Since $(\g\bullet \g)^\perp=N_r(\g, \bullet)$ from Eq. \eqref{eq1}, we finally deduce that $(\g \bullet \g)^\perp = N_\ell(\g, \bullet) \cap N_r(\g, \bullet)$.
    
On the other hand, we have $\om( v \bullet u, w ) = \om( v, [u,w] ) = 0$, for all $u \in Z(\g)$ and $v, w \in \g$. It follows that $\Rr^\bullet_u=\Ll^\bullet_u=0$. Hence, $Z(\g) = (\g \bullet \g)^\perp$. \end{proof}
\sssbegin{Proposition}\label{pr9}
Let $(\g, \br, \om)$ be a quasi-Frobenius nilpotent Lie superalgebra. Then,  $[\g, \g]$ and $Z(\g)$ are both degenerate.
\end{Proposition}
\begin{proof}
Since $(\g, \br)$ is nilpotent, then $[\g, \g] \cap Z(\g) \neq \{0\}$. From Prop. \ref{pr3}, we have $Z(\g) \subseteq [\g, \g]^\perp$. It follows that $\{0\} \neq [\g, \g] \cap Z(\g) \subseteq [\g, \g] \cap [\g, \g]^\perp$, and $\{0\} \neq [\g, \g] \cap Z(\g) \subseteq Z(\g)^\bot \cap Z(\g)$,  which imply that $[\g, \g]$ and $Z(\g)$ are degenerate. 
\end{proof}
\sssbegin{Remark}
Quasi-Frobenius Lie (super)algebras with a degenerate center have been studied in \cite{BE, Fi}. 
\end{Remark}
\sssbegin{Proposition}\label{al}
Let $(\g, \br, \om)$ be a quasi-Frobenius Lie superalgebra, and $\bullet$ be its Levi-Civita product. Then $(\g, \bullet)$ is a right symmetric  superalgebra.
\end{Proposition}

\begin{proof}
Using Prop. \ref{pr3}, we have (for any $u, v, w,z  \in \g):$
\[
\begin{aligned}
\om( & (w \bullet u) \bullet v - w \bullet (u \bullet v) -(-1)^{|u||v|} (w \bullet v) \bullet u + (-1)^{|u||v|}w \bullet (v \bullet u), z) \\
&= \om( w, -[u \bullet v -(-1)^{|u||v|} v \bullet u, z]+[u, [v, z]] - (-1)^{|u||v|}[v, [u, z]]  ) \\
&= -\om( w, [[u, v], z] - [u, [v, z]] + (-1)^{|u||v|}[v, [u, z]]) \\
&= -\om( w, [[u, v], z] +(-1)^{|u||v|+|u||z|} [[v, z], u] + (-1)^{|v||z|+|u||z|}[ [z, u], v])\\&=0.
\end{aligned}
\]
Since $\om$ is non-degenerate, it follows that $(\g, \bullet)$ is a right symmetric superalgebra.
\end{proof}

Let $(\g, \bullet)$ be a right-symmetric superalgebra. Consider the opposite product $\circ$ defined by
$$
u \circ v = -(-1)^{|u||v|} v \bullet u, \quad  \text{for all $u,v\in \g$.}
$$ It is clear that $(\g, \circ)$ is a left-symmetric superalgebra. Moreover, Prop. \ref{al} implies that 
\begin{equation}
\omega(\mathrm{L}_u^\circ (v), w) = -(-1)^{|u||v|} \omega(v, \ad_u(w)), \quad \text{for all $u,v,w \in \g$.}
\label{opp LV}\end{equation}
Therefore, for every quasi-Frobenius  Lie superalgebra $(\g, \br, \omega)$, there exists a unique compatible left-symmetric  structure.
\sssbegin{Proposition} \label{pr4} 
Let $(\A, \circ, \om)$ be a closed pre-symplectic superalgebra. Then the following conditions are equivalent:  
\begin{enumerate}  
    \item[$(i)$] $(\A, \circ)$ is a left symmetric superalgebra such that 
    \begin{equation}
    \om( \Rr_v^\circ (u), w )= (-1)^{|v||w|}\om( u, \Rr_v^\circ (w) ), \quad \text{for all } u, v, w \in \A.
  \label{sym}  \end{equation}  
    \item[$(ii)$] $(\A^-, \om)$ is a quasi-Frobenius Lie superalgebra such that  
    \[
   \omega(\mathrm{L}_u^\circ (v), w) = -(-1)^{|u||v|} \omega(v, \ad_u(w)), \quad \text{for all } u, v, w \in \A.
    \]  
\end{enumerate}  
\end{Proposition}  
\begin{proof}
Assume that $(i)$ holds. Since $(\A, \circ)$ is a left-symmetric superalgebra, it follows that $\A^-$ is a Lie superalgebra. Moreover, since $\omega$ is closed  with respect to $\circ,$ it follows that $\omega$ is also closed with respect to $[\cdot , \cdot ]_\circ$. Hence, $(\A^-, \omega)$ is a quasi-Frobenius Lie superalgebra. 

On the other hand, we have (for every  $u,v,w \in \A$):
\begin{align*}
\omega(\mathrm{L}_u^\circ (v), w)
\overset{\eqref{cyclic}}{=}&
-(-1)^{|u||w| + |u||v|}\,\omega(\mathrm{L}_v^\circ( w), u)
-(-1)^{|u||w| + |w||v|}\,\omega(\mathrm{L}_w^\circ( u), v)
\\\overset{\eqref{sym}}{=}&
-(-1)^{|u||v|}\,\omega(v, \mathrm{L}_u^\circ (w))
-(-1)^{|u||w| + |u||v|}\,\omega(v, \mathrm{L}_w^\circ (u))
\\=& -(-1)^{|u||v|}\,\omega(v, \ad_u(w)).
\end{align*}

Conversely, assume that $(ii)$ holds. Eq. \eqref{opp LV} implies that  the product $\circ$ defines a left-symmetric structure on $\A$. Moreover, for all  $u,v,w \in \A$, we have  
\begin{align*}
\omega(\mathrm{L}_u^\circ (v), w)
&= -(-1)^{|u||v|}\,\omega(v, \ad_u(w))
\\&= (-1)^{|u||v| + |u||w|}\,\omega(v, \ad_w(u))
\\&= -(-1)^{|u||v| + |u||w| + |v||w|}\,\omega(\mathrm{L}_w^\circ( v), u)
\\&= (-1)^{|v||w|}\,\omega(u, \mathrm{L}_w^\circ (v)).
\end{align*}
This completes the proof.
\end{proof}

\subsection{The symplectic product}
In \cite{BCGRS}, the authors established the existence and non-uniqueness of the symplectic product in the case of Lie algebras. Here we extend the result to the super setting. The Levi-Civita product is uniquely determined by a pre-symplectic structure. In contrast, we show in Corollary \ref{No} that the symplectic product on a quasi-Frobenius Lie superalgebra is not unique.
\sssbegin{Definition}(Symplectic product)
Let \((\g, \br, \om)\) be a quasi-Frobenuis Lie superalgebra. A bilinear map  
\[
\star : \g \times \g \to \g
\]
is called a \emph{symplectic product} if it satisfies the following conditions for all \( u, v \in \g \): 
\begin{eqnarray}
    [u, v]_\star &=& [u,v], \label{eq:cyclic1} \\
    (\mathrm{L}_u^\star)^*&=&-\mathrm{L}_u^\star. \label{eq:cyclic2}
\end{eqnarray}
\end{Definition}
\sssbegin{Example}
Following \cite{BM}, consider the periplectic quasi-Frobenuis Lie superalgebra $(D^5, [\;, \;], \omega)$.  In the basis: 
\(\{e_1, e_2  \mid f_1, f_2 \}\) (even $\mid$ odd), it is defined as follows: 
\[
[e_1, f_1] = f_1, 
\qquad 
[e_1, f_2] = f_2, 
\qquad 
[e_2, f_2] = f_1, 
\quad \text{ and } \quad \omega= e_1^*\wedge f_1^*+e_2^*\wedge f_2^*.
\]
A direct computation shows that the following product is symplectic:
\[
\begin{array}{lcrlcrlcrlcrlcrlcr}
e_1 \star e_1 &=& -\tfrac{1}{3} e_1, 
& e_1 \star e_2 &=& -\tfrac{1}{3} e_2, 
& e_1 \star f_1 &=& \tfrac{1}{3} f_1, 
& e_1 \star f_2 &=& \tfrac{1}{3} f_2,\\[3pt]
e_2 \star e_1 &=& -\tfrac{1}{3} e_2, 
& e_2 \star f_2 &=& \tfrac{1}{3} f_1, 
& f_1 \star e_1 &=& -\tfrac{2}{3} f_1, 
& f_2 \star e_1 &=& -\tfrac{2}{3} f_2, 
\\[3pt] 
f_2 \star e_2 &=& -\tfrac{2}{3} f_1.
\end{array}
\]
\end{Example}

\sssbegin{Example} Following \cite{BM}, consider the orthosymplectic quasi-Frobenuis Lie superalgebra $(C^3+A, [\;, \;], \omega)$. In the basis: 
\(\{e_1, e_2  \mid f_1, f_2 \}\) (even $\mid$ odd), it is defined as follows: 
\[
[e_1, f_1] = f_2, 
\qquad 
[f_1, f_1] = e_2, 
\quad
\text{ and }
\quad 
\omega=2 e_1^*\wedge e_2^*-f_1^*\wedge f_2^*.
\]
A direct computation shows that the following product is symplectic:
\[
f_1 \star e_1 = -f_2, 
\quad f_1 \star f_1 = \tfrac{1}{2} e_2.
\]
\end{Example}

\sssbegin{Remark}
We adopt the notations  
$D^5$  and $C_3+A$ from Backhouse \cite{B}, who classified all real Lie superalgebras of dimension 
4. The reader is referred to \cite{B} for an explanation of this nomenclature.
\end{Remark}




We will need the following  lemma later.

\sssbegin{Lemma} \label{pr10} Let $(\A,\star)$ be a  Lie-admissible superalgebra equipped  with a non-degenerate anti-symmetric bilinear form $\omega$ such that  $\Ll^{\star}_u$ is $\om$-anti-symmetric, for any $u\in\A$.  Then the Lie superalgebra $(\A^-,[\; , \;]_\star, \om)$ is quasi-Frobenius.
\end{Lemma}
\begin{proof} For any $u,v,w\in\A$, we have 
	
	     \[\om( [u,v],w)=\om( u\star v,w) -(-1)^{|u||v|}\om( v\star u,w).\]By using the relation
	$\om( \mathrm{L}_u^\star v,w)=-(-1)^{|u||v|}\om( v, \mathrm{L}_u^\star w)$, one can see easily that
	\begin{align*}\oint_{u,v,w}(-1)^{|u||w|}\om(  [u,v],w)&=0\qed 
    . \end{align*}
    \noqed 
\end{proof}

\sssbegin{Remark}
The product constructed in Lemma~\ref{pr10}, whenever it exists, is a symplectic product on $(\A^-,[\, , \,]_\star, \om)$.
\end{Remark}

\sssbegin{Theorem}\label{existence de produit}
On a quasi-Frobenius Lie superalgebra  $(\g, \br, \om)$, there always exists a symplectic product. 
\end{Theorem}

\begin{proof}
    Suppose that $(\g, \br, \om)$ is a quasi-Frobenius Lie superalgebra. Then, there exists a bilinear product $\bullet$ on $\g$ compatible with the Lie bracket $\br$; for instance, by Prop.~\ref{pr3}, one can consider the Levi-Civita product associated with $(\g, \br, \om)$. Define an even bilinear map $M : \g \times \g \to \g$ by the identity:
$$
\om( u \bullet v, w ) -(-1)^{|u||v|}\om( v, u \bullet w ) = \om( M(u,v), w), \quad \text{for all } u,v,w \in \g.
$$
Since $\om$ is anti-symmetric, we have $\om( M(u,v), w ) = (-1)^{|u||v|} \om( v, M(u,w))$. Furthermore, the cyclicity of $\om$ gives
\begin{align*}
\oint_{u,v,w} (-1)^{|u||w|}\om( [u,v], w )
&= \om( u \bullet v, w)- (-1)^{|u||v|} \om( v \bullet u, w ) + \om( v \bullet w, u) \\ &\quad-(-1)^{|v||w|}\om( w \bullet v, u ) + \om( w \bullet u, v ) -(-1)^{|w||u|} \om( u \bullet w, v ) \\
&= (-1)^{|u||w|}\om( M(u,v), w ) + (-1)^{|v||u|}\om( M(v,w), u )\\&\quad +(-1)^{|w||v|} \om( M(w,u), v ).
\end{align*}
Thus, we obtain the identity $
\oint_{u,v,w}(-1)^{|u||w|} \om( M(u,v), w ) = 0.
$

Now, define a new product $\star$ on $\g$ by  
$$
u \star v := u \bullet v + \frac{1}{3}(M(u,v) +(-1)^{|u||v|} M(v,u)).
$$
We check that $\star$ is a compatible product for the Lie superalgebra $(\g, \br, \om)$:
$$
[u, v]_\star = [u, v]_\bullet + \frac{1}{3}(M(u,v) + (-1)^{|u||v|}M(v,u) -(-1)^{|u||v|} M(v,u) - M(u,v)) =[u,v].
$$
Hence, $(\g, \star)$ is also compatible with the Lie superalgebra structure.

To show that $\star$ is $\omega$-anti-symmetric, we have
\begin{align*}
&\om( u \star v, w ) +(-1)^{|u||v|} \om( v, u \star w )
= \om( u \bullet v, w ) +(-1)^{|u||v|} \om( v, u \bullet w ) - \frac{1}{3}\om( M(u,v), w )\\&\quad -\frac{1}{3}(-1)^{|u||v|} \om( M(v,u), w ) -\frac{1}{3}(-1)^{|u||v|}\om( v, M(u,w)) -\frac{1}{3}(-1)^{|u||v||u||w|} \om(v,  M(w,u)) \\
&= \om( M(u,v), w ) - \frac{2}{3}\om (M(u,v), w ) -\frac{1}{3}\om( u, M(v,w)) -\frac{1}{3}(-1)^{|u||v|+|u||w|} \om(v, M(w,u))\\&= \om( M(u,v), w ) - \frac{2}{3}\om (M(u,v), w ) -\frac{1}{3}(-1)^{|u||w|+|v||w|}\om( w, M(u,v))\\
&= (1-\frac{2}{3}-\frac{1}{3}) \om( M(u,v), w )=0.
\end{align*}
because $\om( M(u,v), w ) = (-1)^{|u||v|} \om( v, M(u,w) )$ and  $M$ is cyclic on $\om$. Therefore, $\star$ is $\om$-anti-symmetric, as claimed.
\end{proof}

\sssbegin{Lemma}\label{exists}
Let $E$ be a vector superspace equipped with an orthosymplectic \textup{(}resp. periplectic\textup{)} form $\om$. Then there exists an even \textup{(}resp. odd\textup{)} symmetric bilinear map $f : E \times E \to E$ such that
$$
\omega(f(u, v), w) = -(-1)^{|u||v|}\omega(v, f(u, w)), \quad \textup{for all $u, v, w \in E$}.
$$
\end{Lemma}
\begin{proof}
Let $(E, \omega)$ be a even (resp. odd) symplectic superspace. Consider $T$ as an even (resp. odd) symmetric trilinear form on $E$. For $u, v \in E$, define the map
$$
T(u,v, \cdot): E \to \mathbb{K}
$$
as an even (resp. odd) linear form. Since $\omega$ is non-degenerate, there exists a unique vector $f_T(u,v) \in E$ such that
$$
T(u,v,w) = \omega(f_T(u,v), w) \quad \text{for all } w \in E.
$$
Since $T$ is even (resp. odd) and symmetric on $E$, and $\om$ is even (resp. odd), it is straightforward to show that $f_T: E \times E \to E$ is an even symmetric bilinear map. For any $u,v,w \in E$, we have
$$
\begin{array}{lcl}
\omega(f_T(u,v), w) & = &  T(u,v,w) =(-1)^{|v||w|} T(u,w,v) = (-1)^{|v||w|}\omega(f_T(u,w), v) \\[2mm]
& = & -(-1)^{|u||v|}\omega(v, f_T(u,w)).
\end{array}
$$
This completes the proof.
\end{proof}

\sssbegin{Proposition}\label{autrproduct}
Let $(\g, \br, \om)$ be a quasi-Frobenius Lie  superalgebra, and let $\star$ be a symplectic product on $(\g, \br, \om)$. Consider a bilinear map $f: \g \times \g \to \g$. The product defined by  
\[
u \bar{\star} v = u \star v + f(u, v), \quad \text{for all } u, v \in \g
\]  
is a symplectic product on $(\g, \br, \om)$ if and only if $f$ is symmetric and satisfies the identity  
\[
\om( f(u,v), w)=-(-1)^{|u||v|}\om( v, f(u,w)),\; \text{for all  $u,v,w\in \g.$}
\]
\end{Proposition}

\begin{proof}
Assume that $\star$ is a symplectic product on $(\g, \br, \om)$. For any $u,v\in \g$, we have 
\[\begin{aligned}
\relax [u, v]_{\bar{\star}} &= u \star v + f(u, v) -(-1)^{|u||v|} v \star u -(-1)^{|u||v|} f(u, v) \\
&= [u,v] + f(u, v)- (-1)^{|u||v|} f(v, u).
\end{aligned}
\]
We therefore obtain $[u, v]_{\bar{\star}}=[u,v]$ if and only if $f(u, v) = (-1)^{|u||v|}f(v, u)$. 

On the other hand, for any $u, v, w \in \g$, we have
\[
\begin{aligned}
\om( u \bar{\star} v, w ) + (-1)^{|u||v|}\om( v, u\bar{\star} w) &= \om( u \star v + f(u, v), w ) +(-1)^{|u||v|} \om(v, u \star w + f(u, w)) \\
&= \om( f(u, v), w ) +(-1)^{|u||v|} \om(v, f(u, w)).
\end{aligned}
\]
Thus, 
\[\om( u \bar{\star} v, w ) + (-1)^{|u||v|}\om( v, u \bar{\star} w)= 0 \text{ if and only if } \om( f(u, v), w ) =-(-1)^{|u||v|} \om(v, f(u, w)).
\]
\end{proof}
\sssbegin{Corollary}\label{No}
On a quasi-Frobenius Lie superalgebra  there exist many symplectic products.
\end{Corollary}
\begin{proof}  
According to Prop.  \ref{existence de produit}, there exists a symplectic product $\star$. Moreover, by Lemma \ref{exists}, there exists an even symmetric bilinear map $f: \g \times \g \to \g$ satisfying the condition  
\[
\om( f(u,v),w) = -(-1)^{|u||v|}\om( v, f(u,w)), \quad \text{ for all } u, v, w \in \g.
\]  
Thus, by Prop.  \ref{autrproduct}, the product defined by  
\[
u \bar{\star} v = u \star v + f(u, v), \quad \text{ for all } u, v \in \g,
\]  
is also a symplectic product on $(\g, \br, \omega )$.
\end{proof}

In \cite{BB}, the authors introduced a symplectic connection that depends on the Lie bracket and the symplectic form. In this work, we introduce a natural symplectic product for quasi-Frobenius Lie superalgebras that depends on both the Lie bracket and the form, and we establish its uniqueness in Prop.~\ref{connexion}.

Let $(\g,\br,\omega)$ be a quasi-Frobenius Lie superalgebra. Following the same procedure as in Prop.~\ref{existence de produit}, and using Eq. \eqref{opp LV},  we define a product $\circ$ on $\mathfrak{g}$ by
\[
\omega(\mathrm{L}_u^\circ (v), w) = (-1)^{|u||v|}\,\omega(v,\ad_u(w)), \qquad \text{for all } u,v,w \in \g.
\]
This product satisfies $[u,v] = [u, v]_\circ.$ Moreover, since $\om$ is closed, we have
\[
\omega(u \circ v, w) + (-1)^{|u||v|}\,\omega(v, u \circ w)
  = (-1)^{|u||v|}\,\omega(v \circ u, w).
\]
From the proof of Theorem \ref{existence de produit}, we consider the map
$
M(u,v) = -(-1)^{|u||v|}\, v \circ u.
$
Using this, we define the product $\star$ on $\g$ by
\[
u \star v = u \circ v - \tfrac{1}{3}\Big( (-1)^{|u||v|}\, v \circ u + u \circ v \Big)
= \tfrac{2}{3}\,u \circ v - \tfrac{1}{3}(-1)^{|u||v|}\, v \circ u.
\]
One can verify that this product satisfies
\[
\omega(u \star v, w)
   = \tfrac{2}{3}\,\omega(u \circ v, w)
     - \tfrac{1}{3}(-1)^{|u||v|}\, \omega(v \circ u, w)
   = -(-1)^{|u||v|}\,\tfrac{2}{3}\,\omega(v,[u,w])
     - \tfrac{1}{3}\,\omega(u,[v,w]).
\]
Since $\omega$ is closed, it follows that 
\begin{equation}
\omega(u \star v, w)
   = \tfrac{1}{3}\,\omega([u,v],w)
     + \tfrac{1}{3}(-1)^{|v||w|}\,\omega([u,w],v).
\label{Pr-Symplectic}\end{equation}
Therefore, we obtain a natural symplectic product on $\g$ that depends only on the form $\omega$ and the Lie bracket.

\sssbegin{Proposition}\label{connexion}
Let \((\g, \br, \om)\) be a quasi-Frobenius Lie superalgebra. There exists a unique natural symplectic product that depends on the product  \(\br\) and the form \(\om\). This product is given by \textup{(}for all \(u, v, w \in \g\)\textup{)}:
\begin{equation}
\om( u \star v, w ) = \frac{1}{3} \left( \om( [u, v], w ) +(-1)^{|v||w|} \om( [u, w], v ) \right).
\label{cyclic13}\end{equation}\end{Proposition}
\begin{proof}
 Suppose there exists a symplectic product on \((\g, \br, \om)\) such that:
\[
\om( u \star v, w ) = a \om( [u,v], w ) + b (-1)^{|v||w|} \om( [u,w], v ) + c (-1)^{|v||w|}\om( [v,w], u ),
\]
for some scalars \(a,b,c \in \mathbb{K}\). The  condition \(\om( u \star v, w ) = -(-1)^{|u||v|}\om( v, u \star w )\)  implies that \(b = a\) and \(c = 0\).

Using the identity \([u, v] = [u,v]_\star\), we have
\begin{align*}
\om( [u, v], w ) &= \om(u \star v -(-1)^{|u||v|} v \star u, w) \\&= a \om( [u,v], w ) + (-1)^{|v||w|}a \om( [u,w], v )  -a(-1)^{|u||v|} \om( [v,u], w ) \\&- (-1)^{|u||w|+|u||v|}a \om( [v,w], u ).
\end{align*}
Since the form $\om$ is closed, this simplifies to
$
\om( [u,v], w ) = 3a \om ([u, v], w ),
$
and thus \(a = \frac{1}{3}\).
Therefore, the unique symplectic product compatible with the bracket $\br$ and the form $\omega$ is:
\[
\om( u \star v, w ) = \frac{1}{3} \left( \om( [u,v], w ) + (-1)^{|v||w|}\om( [u,w], v ) \right), \quad \text{for all \(u, v, w \in \g\).}
\qed
\]
\noqed
\end{proof}
\sssbegin{Definition}
We call the product defined in Eq. \eqref{cyclic13} the {\it natural symplectic product.}
\end{Definition}
\section{Flat quasi-Frobenius Lie superalgebras}\label{Sol-nil-sec}

We investigate flat quasi-Frobenius Lie superalgebras, characterized by a natural symplectic product with identically zero curvature. This product induces a left-symmetric structure, extending the Lie algebra case studied in \cite{BELL} to the super setting.

\ssbegin{Definition}
Let $(\g,\br, \om)$ be a quasi-Frobenius Lie superalgebra, and let $\star$ be a symplectic product on $(\g,\br, \om)$.  
We define the curvature operator on $\g$ by
\[
\mathcal{R}(u,v) : = \Ll^\star_{[u,v]} - [\Ll^\star_u, \Ll^\star_v], \qquad \text{ for all } u,v \in \g,
\]
where $\Ll^\star_u$ denotes the left multiplication by $u$ with respect to $\star$.  

Moreover, if this curvature operator vanishes identically, i.e.
\[
\mathcal{R}(u,v) = 0, \quad \text{ for all } u,v \in \g,
\]
then $(\g,\br, \om)$ is called a \emph{flat quasi-Frobenius Lie superalgebra}.

\end{Definition}                        
It is easy to see that the curvature $\mathcal{R}$ of the symplectic product $\star$ vanishes identically if and only if $(\g, \star)$ is a left-symmetric superalgebra.

\ssbegin{Example} In this example, we study an important class of nilpotent Lie superalgebras, known as \emph{filiform Lie superalgebras}. We restrict our attention to the simplest ones, denoted by $L^{n,m}$, where $n,m \in \mathbb{N}$. These are defined on the homogeneous basis
\[
(e_1, \dots, e_n \mid f_1, \dots, f_m),
\]
whose only nonzero bracket products are given by
\begin{align*}
[e_1, e_i] &= e_{i+1}, \qquad 2 \le i \le n-1, \\
[e_1, f_j] &= f_{j+1}, \qquad 1 \le j \le m-1.
\end{align*}

\underline{The case where $m=n$.}  Following \cite{BM}, $L^{n,n}$ admits a closed periplectic form given by
\[
\omega =\lambda\, e_1^* \wedge f_n^*+\mu \sum_{i=2}^{n} (-1)^{\,n-i} e_i^* \wedge f_{\,n+1-i}^*, \quad \text{where $\lambda\mu \neq 0$.}
\]
Using Eq~\eqref{cyclic13}, the associated natural symplectic product $\star$ is given by: (non-displayed products are zero)
\begin{align*}
e_1 \star e_i &= \frac{2}{3} e_{i+1}, 
&\qquad 
e_i \star e_1 &= -\frac{1}{3} e_{i+1}, \quad e_i\star f_{n-i}=f_{n-i}\star e_i=  \frac{\mu}{3\la} (-1)^{n+1-i}f_{n}, \\
e_1 \star f_j &= \frac{2}{3} f_{j+1}, 
&\qquad 
f_j \star e_1 &= -\frac{1}{3} f_{j+1}, \quad  e_1\star e_1=  \frac{\la}{3\mu} (-1)^{n+1}e_{2},\\
e_1 \star f_{n-1} &= \frac{1}{3} f_{n}, 
&\qquad 
f_{n-1} \star e_1 &= -\frac{2}{3} f_{n},
\end{align*}
where $2 \le i \le n-1$ and $1 \le j \le n-2$.

A direct computation gives 
\[
\mathcal{R}(e_1,e_i)(e_1)
=
\big(\Ll^\star_{[e_1,e_i]} - [\Ll^\star_{e_1}, \Ll^\star_{e_i}]\big)(e_1)
=-
\frac{1}{9}\, e_{i+2}
\neq 0.
\]
Therefore, $(L^{n,n},[\;,\;], \omega)$ is not flat. 

\underline{The case where $n$ is even and $m$ is odd}. Following \cite{BM}, the Lie superalgebra $L^{n,m}$ admits a closed  orthosymplectic form given by
\[
\omega 
=
\lambda\, e_1^* \wedge e_n^*
+
\mu \sum_{i=2}^{\frac{n}{2}}
(-1)^{\frac{n}{2}-i}
e_i^* \wedge e_{n+1-i}^*
+
(-1)^{\frac{m+1}{2}} \nu
\sum_{i=1}^{\frac{m-1}{2}}
(-1)^i
f_i^* \wedge f_{m+1-i}^*
+
\frac{\nu}{2}
f_{\frac{m+1}{2}}^* \wedge f_{\frac{m+1}{2}}^*,
\]
where $\lambda\mu\nu \neq 0$.

Using again Eq.~\eqref{cyclic13}, one obtains  the associated product $\star$ as in the previous case. In particular, $(L^{n,m},[\;, \;], \omega)$ is not flat. We omit the details. 

\underline{The case where $n$ is even and $m=n-2$}. Following \cite{BM}, the Lie superalgebra $L^{n,n-2}$ admits a closed, antisymmetric, nondegenerate form that is {\bf nonhomogeneous}. This lies beyond the scope of the present paper.




\underline{The case where $n$ is even and $m\not =n-2$.} It was conjectured in \cite{BM} that no closed and nondegenerate form exists.

\end{Example}

As we proceed in this work, we consider the \emph{natural symplectic product} defined by Eq.~ \eqref{cyclic13}. 

Let $(\g,\br,\om)$ be a quasi-Frobenius Lie superalgebra, and let $\star$ be the {\it natural} symplectic product. The left and right multiplications with respect to $\star$ are expressed as (for all $u \in \g$):
\begin{equation}\Ll_u^\star = \tfrac{1}{3}\big(\ad_u - \ad_u^{*}\big), 
\qquad 
\Rr_u^\star =- \tfrac{1}{3}\big(2\ad_u + \ad_u^{*}\big),
\label{LR}\end{equation}
where $\ad_u^{*}$ denotes the adjoint of $\ad_u$ with respect to the form $\omega$.

\ssbegin{Proposition}\label{orth}
Let $(\g,\br,\omega)$ be a quasi-Frobenius Lie superalgebra, and let $\star$ be the natural symplectic product.  Then:
\begin{enumerate}
\item[\textup{(i)}] For all $u \in [\g,\g]^{\perp}$, we have
$
\Ll_u^\star = -2\Rr^\star_u=\tfrac{2}{3}\,\ad_u$. Moreover,
\[
[\g,\g]^{\perp} = \{\, u \in \g \;\mid\; (\Rr_u^{\star})^* = -\Rr^\star_u \,\}.
\]
\item[\textup{(ii)}] The center of $\g$ satisfies
\[
Z(\g) = (\g \star \g)^{\perp} = N_r(\g, \star)= N_r(\g, \star)\cap N_\ell(\g, \star) 
       = N_\ell(\g, \star)\cap [\g,\g]^{\perp}.
\]
\end{enumerate}
\end{Proposition}
\begin{proof}
For part (i),  since $\omega$ is closed, we have  (for all $u,v,w \in \g$):
\[
\begin{array}{lcl}
\omega(\mathrm{L}_u^\star v, w)
  & = & \tfrac{1}{3}\,\omega([u,v],w)
     + \tfrac{1}{3}(-1)^{|v||w|}\,\omega([u,w],v)\\[2mm]
  & =& \tfrac{2}{3}\,\omega([u,v],w)
     + \tfrac{1}{3}(-1)^{|u|(|v|+|w|)}\,\omega([v,w],u).
     \end{array}
\]
Now, if $u \in [\g,\g]^{\perp}$ and for all $v,w \in \g$, we have 
\[
\omega(\mathrm{L}_u^\star (v), w)
   = \tfrac{2}{3}\,\omega([u,v],w).
\]
Since $\omega$ is non-degenerate, it follows that
$
\Ll^\star_u = \tfrac{2}{3}\,\ad_u.
$
Using the relation $\ad_u = \Ll^\star_u - \Rr^\star_u$, we deduce
$
\Rr^\star_u = -\tfrac{1}{3}\,\ad_u.$ Thus, for $u \in [\g,\g]^{\perp}$, we obtain $(\Rr_u^\star)^{*} = -\Rr_u^\star$.

For part (ii), let $u \in Z(\g)$. Since $\ad_u = 0$, Eq.\eqref{LR} implies that $\Ll^\star_u = -\Rr^\star_u.$ Therefore, $\Ll^\star_u = \Rr^\star_u = 0$.  

Now, let $u \in (\g \star \g)^\perp$. For all $v,w \in \g$, we have
\[
0 = \omega(v \star w, u) 
   = -(-1)^{|v||w|}\,\omega(w, v \star u).
\]
Since $\omega$ is non-degenerate, this implies $\Rr^\star_u = 0$. Moreover,
$
\omega([v,w],u) = 0,$ for every $v,w \in \g,$
so that
$
(\g \star \g)^\perp \subset [\g,\g]^\perp.
$

From (i), we know that for $u \in [\g,\g]^\perp$, we have 
$
\Ll^\star_u = -\tfrac{1}{2}\Rr^\star_u.
$
Therefore, for $u \in (\g \star \g)^\perp$, we deduce
$
\Ll^\star_u = \Rr^\star_u = 0.
$ Finally, let $u \in N_\ell(\g, \star) \cap [\g,\g]^\perp$. By Part (i) again, we conclude that
$
\Ll^\star_u = \Rr^\star_u = 0.
$
We thus obtain the desired equalities:
\[
Z(\g) = (\g \star \g)^\perp = N_r(\g, \star)\cap N_\ell(\g, \star) 
      = N_\ell(\g, \star)\cap [\g,\g]^\perp.
      \qed
\]
\noqed
\end{proof}

Let $(\g,\br,\omega)$ be a  quasi-Frobenius Lie superalgebra.  
The vanishing of the curvature, $\mathcal{R}=0$, is equivalent to the relation (for all $u,v\in \g$)
\begin{equation}
\Rr^\star_{u \star v} -(-1)^{|u||v|} \Rr^\star_v \circ \Rr^\star_u = [\Ll^\star_u, \Rr^\star_v]. 
\label{flat-LR}\end{equation}

\ssbegin{Proposition}\label{ort-2}
Let $(\g,\br,\omega)$ be a flat quasi-Frobenius Lie superalgebra. For all  $u,v \in [\g,\g]^\perp$, we have
\[
\ad_u \circ \ad_v = 0; \qquad \text{in particular,} \qquad \ad_u^2=0.
\]
\end{Proposition}

\begin{proof}
Let $u,v \in [\g,\g]^\perp$ and $w\in \g$.  
From the definition of the natural symplectic product \eqref{cyclic13}, we have
\[
\omega(u \star v, w) 
 = \tfrac{1}{3}\,\omega([u,v],w)
   + \tfrac{1}{3}(-1)^{|v||w|}\,\omega([u,w],v).
\]
Since $v \in [\g,\g]^\perp$, it follows that $\omega([u,w],v)=0$ for all $w\in \g$.  
Hence
\[
\omega(u \star v, w)=\tfrac{1}{3}\,\omega([u,v],w).
\]
On the other hand, as $u\in [\g,\g]^\perp$, we have $\ad_u^*=-\ad_u$ and therefore
\[
\omega([u,v],w) = -(-1)^{|u||v|}\,\omega(v,[u,w])=0.
\]
Thus $\omega(u\star v,w)=0$ for all $w$, and by the non-degeneracy of $\omega$, we deduce
$
u\star v = 0, \, \text{ for all } \, u,v \in [\g,\g]^\perp.
$

Using the flatness relation (see Eq. \eqref{flat-LR}), we get
\[
0 = \Rr^\star_{u\star v} - (-1)^{|u||v|}\Rr^\star_v \circ \Rr^\star_u - [\Ll^\star_u, \Rr^\star_v].
\]
Since $u\star v=0$, and the fact that  $\Ll^\star_u = -\tfrac{1}{2}\Rr^\star_u = \tfrac{2}{3}\ad_u$, for $u\in [\g,\g]^\perp$, by Prop.~\ref{orth}, this reduces to
$$
0 = -  (-1)^{|u||v|}\Rr^\star_v \circ \Rr^\star_u - [\Ll^\star_u, \Rr^\star_v]=  -(-1)^{|u||v|}\tfrac{1}{9}\ad_v\circ \ad_u+\tfrac{2}{9}[\ad_u, \ad_v]
.$$ Therefore, 
\[
0 = \tfrac{2}{9}\ad_u\circ \ad_v- \tfrac{1}{3}(-1)^{|u||v|}\ad_v\circ \ad_u.
\]
But since $u,v\in[\g,\g]^\perp$, one checks that $[\ad_u,\ad_v]=\ad_{[u,v]}=0$.  
Thus, the relation reduces to
\[
\ad_v \circ \ad_u = 0.
\] In particular, setting $u=v$, we deduce $ \ad_u^2=0.$ This completes the proof.
\end{proof}

\ssbegin{Theorem}\label{solvable-even}
Let $(\g, \br,\omega)$ be a flat orthosymplectic quasi-Frobenius Lie superalgebra.   Then, the  Lie superalgebra $(\g,\br)$ is solvable.
\end{Theorem}
\begin{proof}
Since $\omega$ is even, the even part $(\g_{\bar 0},\br_{\bar 0},\omega_{\bar 0})$ is a flat quasi-Frobenius Lie algebra.  
It is shown in \cite{BELL} that such a Lie algebra is nilpotent, and hence solvable.  Moreover, according to \cite{Kac, S}, the solvability of $(\g_{\bar 0},\br_{\bar 0})$ implies that the  Lie superalgebra $(\g,\br)$ is solvable.  
Hence, $(\g,\br)$ is a solvable Lie superalgebra.
\end{proof}
\ssbegin{Theorem}\label{nilpotent-odd}
Let $(\g, \br,\omega)$ be a flat periplectic quasi-Frobenius Lie superalgebra. Then, the  Lie superalgebra $(\g,\br)$ is nilpotent.
\end{Theorem}
\begin{proof}
Let $(\g = \g_{\bar{0}} \oplus \g_{\bar{1}}, \br, \omega)$ be a flat periplectic quasi Frobenius Lie superalgebra. According to Lemma \ref{BM}, $\dim(\g_{\bar 0})=\dim(\g_{\bar 1})$. Without loss of generality, we may assume that $\g$ admits a homogeneous basis (even $\mid$ odd):
\[
\{U_1, \dots, U_n \, \mid \,  V_1, \dots, V_n\}.
\]
The Lie bracket and the periplectic form $\omega$ are completely determined by the following relations: 
\[
\begin{cases}
[U_i, U_j] = \sum_{k=1}^n C_{ij}^k U_k, & C_{ij}^k =- C_{ji}^k, \\
[U_i, V_j] = \sum_{k=1}^n E_{ij}^k V_k, & E_{ij}^k =- E_{ji}^k, \\
[V_i, V_j] = \sum_{k=1}^n F_{ij}^k U_k, & F_{ij}^k = F_{ji}^k, \\
\omega(U_i, V_i) =1, & 
\end{cases}
\]
where all structure constants satisfy the super Jacobi identity and $\omega$ is closed.

Now, consider the Lie algebra $\bar{\g}$ and the bilinear form $\bar{\omega}$ obtained by forgetting the super structure.  
The same set
\[
\{\overline{U_1}, \dots, \overline{U_n}, \overline{V_1}, \dots, \overline{V_n}\}
\]
forms a basis of $\bar{\g}$, and the bracket and the symplectic form are defined by
\[
[\overline{U_i}, \overline{U_j}] = \sum_{k=1}^n C_{ij}^k \overline{U_k}, \quad 
[\overline{U_i}, \overline{V_j}] = \sum_{k=1}^n E_{ij}^k \overline{V_k}; \quad 
\bar{\omega}(\overline{U_i}, \overline{V_i}) = 1.
\]
It is straightforward to check that $\bar{\g}$ is a Lie algebra and $\bar{\omega}$ is closed  on $\bar{\g}$.  
In particular, $\bar{\g}$ inherits a $\mathbb{Z}_2$-grading, so $(\bar{\g}, \br, \bar{\omega})$ is a quasi-Frobenius Lie algebra.  

We define the natural symplectic product $\bar{\star}$ on $\bar{\g}$ by
\[
\bar{\omega}(\bar{u} \,\bar{\star}\, \bar{v}, \bar{w}) 
   = \frac{1}{3}\,\bar{\omega}([\bar{u},\bar{v}], \bar{w}) + \frac{1}{3}\,\bar{\omega}([\bar{u},\bar{w}], \bar{v}), 
   \qquad \text{for all }  \bar{u},\bar{v},\bar{w} \in \bar{\g}.
\]
Recall that the structure constants $C_{ij}^k$ and $E_{ij}^k$ completely determine the Lie algebra structure on $\bar{\g}$ and the symplectic form $\bar{\omega}$.  
Let us show that the product $\bar{\star}$ depends only on the Lie bracket of $\bar{\g}$ and the symplectic form $\bar{\omega}$. This  would imply that  $(\bar{\g}, \br, \bar{\omega})$ is a flat quasi-Frobenius Lie algebra. Indeed, 
\begin{align*}
 &\om(U_i\star U_j, V_k)= \frac{1}{3}\om([U_i, U_j], V_k)+\frac{1}{3} \om([U_i, V_k], V_j)\\&=\frac{1}{3} \sum_{l=1}^n C_{ij}^l\om( U_l, V_k)+ \frac{1}{3} \sum_{l=1}^n E_{ik}^l\om( V_l, U_j) 
 = \frac{1}{3} \sum_{l=1}^n C_{ij}^l\om( U_l, V_k)- \frac{1}{3} \sum_{l=1}^n E_{il}^j \om( U_l, V_k) \\&= \frac{1}{3} \sum_{l=1}^n (C_{ij}^l-E_{il}^j)\om( U_l, V_k).
\end{align*}
Therefore, $U_i\star U_j= \frac{1}{3}\sum_{l=1}^n (C_{ij}^l-   E_{il}^j)U_l$, and hence $\bar{U_i}\bar{\star} \bar{U_j}= \frac{1}{3} \sum_{l=1}^n (C_{ij}^l-   E_{il}^j)\bar{U_l}$. As the constant structures $C_{ij}^k$ and $E_{ij}^k$ determine completely the product $\bar{\star}$ and $\star$ is flat, it follows that $(\bar{\g}, \br, \bar{\om})$ is flat.

It is shown in \cite{BELL} that every flat quasi-Frobenius Lie algebra of this type is nilpotent. Therefore, $\bar{\g}_0$ is nilpotent, and so is the adjoint representation of $\bar{\g}_0$ on $\bar{\g}_1$.

We observe that the structure constants $C_{ij}^k$ determine the nilpotent Lie algebra structure on $\bar{\g}_0$, while the constants $E_{ij}^k$ ensure that the representation $\rho_{\bar{\g}}$ is nilpotent.  
Consequently, $\bar{\g}_0$ is nilpotent, and the corresponding adjoint action on $\bar{\g}_1$ is nilpotent as well.  
This implies that the Lie superalgebra $(\g = \g_{\bar 0} \oplus \g_{\bar 1}, [\; , \;])$ is nilpotent \cite{RSS}. 
\end{proof}
\ssbegin{Remark}
The reader may wonder why the argument in Theorem \ref{nilpotent-odd} fails in the orthosymplectic case. The point is that, upon forgetting the superstructure, the Lie algebra 
$(\bar{\fg},[\;, \;], \bar{\omega})$ need not be symplectic. Therefore, the result of \cite{BELL} does not apply, and one cannot deduce that $(\bar{\fg},[\;, \;])$ is nilpotent.
\end{Remark}
\ssbegin{Proposition}\label{degenerate}
Let $(\g, \br,\omega)$ be a flat  quasi-Frobenius Lie superalgebra.  
Then, $[\g,\g]$ is degenerate.
\end{Proposition}
\begin{proof}
If $\omega$ is odd, then according to Theorem~\ref{nilpotent-odd}, $(\g, \br)$ is nilpotent.  
In particular, we have 
$
[\g,\g] \cap Z(\g) \neq \{0\}.
$ 
Since $\omega$ is periplectic, we also have 
$
Z(\g) \subseteq [\g,\g]^\perp,
$ 
and therefore 
$$
[\g,\g] \cap [\g,\g]^\perp \neq \{0\}.
$$

Now, suppose $\omega$ is even, and assume that 
$
\g = [\g,\g] \oplus [\g,\g]^\perp.
$
According to Prop.~\ref{solvable-even}, $(\g,\br)$ is solvable.  
Hence, $[\g,\g]$ is nilpotent, and its center satisfies 
$
Z([\g,\g]) \neq \{0\}.
$ 

We claim that $Z([\g,\g]) \subseteq Z(\g)$. Indeed, let $u \in Z([\g,\g])$. For any $v \in [\g,\g]^\perp$ and $w \in [\g,\g]$, we have
\[
(-1)^{|u||w|} [[u,v], w] = -(-1)^{|u||v|} [[v,w], u] - (-1)^{|v||w|} [[w,u], v] = 0.
\] 
On the other hand, Prop.~\ref{ort-2} implies that $\mathrm{ad}_v \circ \mathrm{ad}_w = 0$, for any $w \in [\g,\g]^\perp$, and hence 
$
[[u,v],w] = 0.
$
Thus, $[u,v] \in Z(\g) \subseteq [\g,\g]^\perp$, which implies $[u,v] = 0$.  This shows that $u \in Z(\g)$. 

Now the fact that $Z([\g,\g]) \subseteq Z(\g)$ and $Z(\g) \subseteq [\g,\g]^\perp$ leads to $Z([\g, \g])\subseteq [\g,\g]\cap [\g,\g]^\perp=\{0\} $, which is a contradiction.
\end{proof}
\ssbegin{Proposition}\label{sided}
Let $(\g, \br, \om)$ be a flat quasi-Frobenius Lie superalgebra. Then,  $N_\ell(\g, \star)$ is a two-sided ideal of $(\g, \star)$ and $[\g, [\g, \g]^\bot]\subseteq N_\ell(\g, \star)$. 
\end{Proposition}

\begin{proof}
It is clear that $N_\ell(\g, \star)$ is a right ideal. Let us now show that it is also a left ideal.  
Let $u \in N_\ell(\g, \star)$ and $v \in \g$. Since $(\g, \br, \omega)$ is flat, we have
\[
(-1)^{|u||v|}\Ll^\star_{v \star u} = \Ll^\star_{[u,v]} = [\Ll^\star_u, \Ll^\star_v] = 0.
\]
It follows that $\Ll^\star_{v \star u} = 0$, and hence $v \star u \in N_\ell(\g, \star)$. 

Now, let $u, v \in \g$. From \eqref{LR}, one has 
$
\Ll^\star_u = \frac{1}{3} (\ad_u - \ad_u^*).
$
Thus,
\begin{eqnarray*}\nonumber
0 &=& \Ll^\star_{[u,v]} - [\Ll^\star_u, \Ll^\star_v] \\ \nonumber
&=& \frac{1}{3}(\ad_{[u,v]} - \ad_{[u,v]}^*) - \frac{1}{9} [\ad_u - \ad_u^*, \ad_v - \ad_v^*] \\ \nonumber
&=& \frac{1}{3}\ad_{[u,v]} - \frac{1}{3}\ad_{[u,v]}^* - \frac{1}{9} [\ad_u, \ad_v] +\frac{1}{9}[\ad_u,  \ad_v^*] +\frac{1}{9}[\ad_u^*,  \ad_v]  -\frac{1}{9}[\ad_u^*,  \ad_v^*]\\
&=& \frac{2}{9}\ad_{[u,v]} - \frac{2}{9}\ad_{[u,v]}^* + \frac{1}{9}[\ad_u^*, \ad_v] + \frac{1}{9}[\ad_u, \ad_v^*].
\end{eqnarray*}
Now, $\ad_v^* = -\ad_v$, for any $u \in \g$ and $v \in [\g, \g]^\perp$. Substituting this into the previous relation, we obtain
\begin{eqnarray*}
0 &=& \frac{2}{9}\ad_{[u,v]} - \frac{2}{9}\ad_{[u,v]}^* + \frac{1}{9}[\ad_u^*, \ad_v] + \frac{1}{9}[\ad_u, \ad_v^*]= \frac{1}{9}\ad_{[u,v]} - \frac{1}{9}\ad_{[u,v]}^*.
\end{eqnarray*}
Therefore, $\ad_{[u,v]} = \ad_{[u,v]}^*$, which implies $\Ll^\star_{[u,v]} = 0$, and hence $[u,v] \in N_\ell(\g, \star)$.  
\end{proof}
\section{Flat  double extensions of flat quasi-Frobenius Lie superalgebras}\label{d-ext}
In this section, we introduce the notion of {\it flat}  double extension for flat quasi-Frobenius Lie superalgebras. It is inspired by the concept of double extension, originally introduced by Medina and Revoy in \cite{MR} for quadratic Lie algebras.  It is worth noticing that this notion is different from the notion of double extensions of quasi-Frobenius Lie superalgebras introduced in \cite{BE, BM}.

\subsection{Central extension of left-symmetric superalgebras}\label{5.1-dia}
Let $(\A, \bullet)$ be a left-symmetric superalgebra,  let $\h := \K e$ be a one-dimensional vector superspace, and let $\mu : \A \times \A \rightarrow \K$ be a homogeneous bilinear map of a given parity.  
We define a new product $\diamond$ on the vector space
$
\widetilde{\A} := \A \oplus \K e
$
as follows:
\[
(u+\alpha e)\diamond(v+\beta e) := u \bullet v + \mu(u,v)e,
\qquad \text{for all } u,v \in \A,\; \alpha, \beta \in \K.
\]
Following \cite{BBE2}, the superalgebra $(\widetilde{\A}, \diamond)$ is left-symmetric if and only if the following condition holds
\begin{equation}
\mu(\mathrm{L}_u^\bullet (v), w) - \mu(u, \mathrm{L}_v^\bullet (w)) 
= (-1)^{|u||v|}\big(\mu(\mathrm{L}_v^ \bullet (u), w) - \mu(v, \mathrm{L}_u^ \bullet (w))\big), \quad \text{for all $u,v,w\in \A$}. 
\label{con-left}
\end{equation}  
In this case, $(\widetilde{\A}, \diamond)$ is called the \emph{central extension} of $(\A, \bullet)$ by means of  $\mu$.
\subsection{Almost semi-direct product of left-symmetric superalgebras}\label{5.2-dia}

Let $(\A, \bullet)$ be a left-symmetric superalgebra. Let $\mathcal{S} := \K d$ be a one-dimensional vector superspace, let $
D : \A \rightarrow \A, \;
\xi : \A \rightarrow \A,
$ be homogeneous endomorphisms of parity $|d|$, and let $b_0 \in \A_{\bar 0}$  and a scalar $\lambda \in \K$.  
On the vector superspace 
$
\overline{\A} := \K d \oplus \A,
$
we define a new product $\bar{\bullet}$ as follows: (here $\lambda=0$ if $d$ is odd)
\[
d \bar{\bullet} d := \lambda d + b_0, \qquad
d \bar{\bullet} u := D(u), \qquad
u \bar{\bullet} d := \xi(u), \qquad
u \bar{\bullet} v := u \bullet v,
\]
for all $u,v \in \A$.  

Then $(\overline{\A}, \bar{\bullet})$ is a left-symmetric superalgebra if and only if the following conditions are satisfied (for all $u,v \in \A$): 
\begin{equation}
\begin{cases}
\xi([u, v]_\bullet)
=\mathrm{L}_ u^\bullet (\xi(v)) - (-1)^{|u||v|} \mathrm{L}_v^\bullet (\xi(u)), \\[0.4em] 
D(\mathrm{L}_u^\bullet (v))= \mathrm{L}_{D(u)}^\bullet (v) + (-1)^{|d||u|} \mathrm{L}_u^ \bullet (D(v))- (-1)^{|d||u|}\mathrm{L}_{\xi(u)}^ \bullet (v),\\[0.4em]
(\xi\circ D-D\circ \xi)(u)
= (-1)^{|d||u|} \big( \xi^2 - \lambda \xi - \Rr^\bullet_{b_0} \big)(u),\\[0.4em]
(1 - (-1)^{|d|}) \big(D^2 - \lambda D - \Ll^\bullet_{b_0}\big) = 0,\;
(1 - (-1)^{|d|}) (\xi - D)(b_0) = 0.
\end{cases}
\label{produit-semi}
\end{equation}
\ssbegin{Definition}
If $(D, \xi, b_0, \lambda)$ satisfies the compatibility conditions \eqref{produit-semi},  
then the left-symmetric superalgebra $(\overline{\A}, \bar{\bullet})$ is called the 
\emph{almost semi-direct product} of the left-symmetric superalgebra $(\A,\bullet)$ 
by the one-dimensional vector superspace $\K d$ by means of $(D, \xi, b_0, \lambda)$.  
In this case, the quadruple 
$
(\mathcal{S}, \A, D, \xi, b_0, \lambda)
$
is referred to as a \emph{context of almost semi-direct product of left-symmetric superalgebras}.

\end{Definition}
\ssbegin{Remark}
Let $
(\mathcal{S}, \A, D, \xi, b_0, \lambda)
$
be a context of almost semi-direct product of left-symmetric superalgebra.  
The endomorphism defined by $
\Delta(u) := (D-\xi)(u)
$, for all $u\in \A$, 
is a derivation on $(\A, [\,,]_{\bullet})$. Indeed, the statement follows from the first and the second conditions of the system \eqref{produit-semi}. 
\end{Remark}

\subsection{Flat  double extensions of flat quasi-Frobenius Lie superalgebras}
In this paragraph, we introduce the notion of a {\it flat double extension} of flat quasi-Frobenius Lie superalgebras.  
We then show that if a flat orthosymplectic quasi-Frobenius  Lie superalgebra is nilpotent,  then any its flat double extension of is also nilpotent.

We explain how the double extension is constructed in the following diagram.

\begin{equation*}
\begin{tikzcd}
(\fb,[\;,\;]_{\fb}, \omega_\fb)
  \arrow[r, "\text{Step 1}"]
&
(\fb,\bullet)
  \arrow[r, "\text{Step 2}"]
&
(\fb\oplus V^* , \diamond)
  \arrow[r, "\text{Step 3}"]
&
(\fg:=V\oplus \fb\oplus V^* , {\bar \bullet})
  \arrow[d, "\text{Step 4}"] \\
&&&
(\fg, [\;,\;], \omega)
\end{tikzcd}
\end{equation*}
In the above diagram, we now describe the construction in detail.

\medskip

\noindent
\textbf{Step 1.}
Select the natural symplectic product associated with 
\((\mathfrak b, [\;,\;]_{\mathfrak b}, \omega_{\mathfrak b})\).

\medskip

\noindent
\textbf{Step 2.}
Construct the central extension of \((\mathfrak b, \bullet)\) (see Section~\ref{5.1-dia}), 
where the cocycle is given by
\[
\mu(\cdot,\cdot) = \omega_{\mathfrak b}(\xi(\cdot), \cdot),
\]
and $\xi$ satisfies the conditions ensuring that equation~\eqref{con-left} holds.

\medskip

\noindent
\textbf{Step 3.}
Form the semidirect product as in Section~\ref{5.2-dia}, using the data
\[
(D = \xi^* - \xi, \xi, b_{\bar 0}, \lambda),
\]
subject to the conditions that guaranty that equation~\eqref{produit-semi} is valid.

\medskip

\noindent
\textbf{Step 4.}
Compute the commutator $[\;,\;]_{\bar \bullet}$ to obtain the Lie bracket on 
$\mathfrak g$, and extend the form $\omega_{\mathfrak b}$ naturally to $\mathfrak g$.

\subsubsection{Even flat double extensions of orthosymplectic quasi-Frobenius Lie superalgebras.}

Let $(\mathfrak{b}, \br_\mathfrak{b}, \omega_\mathfrak{b})$ be a flat orthosymplectic quasi-Frobenius  Lie superalgebra. Denote by $\star$ the natural symplectic product associated with it.  
Let $V = \K d$ be a one-dimensional vector superspace and $V^* = \K e$ be its dual.  

In order to introduce the process of flat double extension by $V$, we consider an even linear map $
\xi : \mathfrak{b} \rightarrow \mathfrak{b},
$ and an even element $b_0 \in \mathfrak{b}_{\bar 0}$. We extend the linear maps $\xi$ and $D := \xi^* - \xi$ to 
$
\widetilde{D}, \widetilde{\xi} : \mathfrak{b} \oplus V^* \rightarrow \mathfrak{b} \oplus V^*,
$
as follows:   
\begin{equation}
\label{eq:double-tilde}
\begin{array}{c}
\widetilde{D}(u + \alpha e)
:= (\xi^* - \xi)(u) + \dfrac{1}{3}\, \omega(b_0, u)\, e, \quad   
\widetilde{\xi}(u + \alpha e)
:= \xi(u) - \dfrac{2}{3}\, \omega(b_0, u)\, e.
\end{array}
\end{equation}
for all $u \in \mathfrak{b}$ and $\alpha \in \K$. We also choose $\widetilde{b_0}
:= \dfrac{1}{3}\, b_0 \in \mathfrak{b} \oplus V^*$.

\medskip

\sssbegin{Lemma} \label{Lemma1}
The superspace $\mathcal{H} := \mathfrak{b} \oplus V^*$ equipped with the product
\[
(u + \alpha e) \diamond (v + \beta e)
:= u \star v + \omega(\xi(u), v)\, e,
\qquad \text{for all } u,v \in \mathfrak{b},\; \alpha, \beta \in \K,
\]
is a left-symmetric superalgebra if and only if
\begin{equation}
\xi([u,v]_{\mathfrak b}) = \mathrm{L}_u^\star (\xi(v)) - (-1)^{|u||v|} \mathrm{L}_v^\star (\xi(u)), \quad \text{for all $u,v \in \mathfrak{b}$.}
\label{eq:xi-condition0}
\end{equation}
\end{Lemma} 

\begin{proof} 
The superalgebra $(\mathcal{H}, \diamond)$ is left-symmetric if and only if the bilinear map 
\[
\mu(u,v) := \omega(\xi(u), v)
\]
satisfies Eq. \eqref{con-left}, which is equivalent to Eq.  \eqref{eq:xi-condition0}. 
\end{proof}

Now, let us also assume that $(V, \mathfrak{b}, \xi^* - \xi, \xi^*, \tfrac{1}{3} b_0, 0)$  is a context of almost semi-direct product 
of left-symmetric superalgebras.

\medskip

\sssbegin{Lemma}\label{Lemma2}
$(V, \mathcal{H} := \mathfrak{b} \oplus V^*, \widetilde{D}, \widetilde{\xi}, \widetilde{b_0}, 0)$ is a context of almost semi-direct product 
of left-symmetric superalgebras if and only if  for all $u,v \in \mathfrak{b}$, the following conditions hold:
\begin{equation}\label{eq:claim2}
\left\{
\begin{array}{lcllcl}
\xi ([u,v]_{\mathfrak b}) & = &  \mathrm{L}_u^\star (\xi(v)) - (-1)^{|u||v|} \mathrm{L}_v^\star (\xi(u)), & 
\xi^* \circ \xi & = &  \tfrac{1}{3}(\Rr^\star_{b_0}+ (\Rr^\star_{b_0})^*), \\[2mm]
D(u\star v) & = & D(u)\star v+u\star D(v)-\xi(u)\star v , & 
[\xi, \xi^*] & = & \xi^2 - \tfrac{1}{3} \Rr^\star_{b_0},\;
b_0\in \ker{D},
\end{array}
\right.
\end{equation}
where $D=\xi^*-\xi$.

\end{Lemma}

\begin{proof}
 Let $u, v \in \mathfrak{b}$ and $\alpha, \beta \in \mathbb{K}$. 
A straightforward computation shows that the equality
\[
\widetilde{\xi}([u,v]_{\mathfrak b}) = \mathrm{L}_u^\diamond (\widetilde{\xi}(v)) - (-1)^{|u||v|} \mathrm{L}_v^ \diamond (\widetilde{\xi}(u))
\]
is equivalent to
\[
\xi ([u,v]_{\mathfrak b}) = \mathrm{L}_u^\star (\xi(v)) - (-1)^{|u||v|} \mathrm{L}_v^\star (\xi(u)), 
\quad \text{and} \quad 
\xi^* \circ \xi = \tfrac{1}{3}(\Rr^\star_{b_0}+ (\Rr^\star_{b_0})^*).
\]
Moreover, the relation
\[
\widetilde{D}(u) \diamond v - \widetilde{D}(u  \diamond v)
= (-1)^{|d||u|} \big(\widetilde{\xi}(u)  \diamond v - u  \diamond \widetilde{D}(v)\big)
\]
is equivalent to
\[
D(u\star v)=D(u)\star v+u\star D(v)-\xi(u)\star v ,
\]
together with
\[
[\xi^*, \xi] = \xi^2 - \tfrac{1}{3} \Rr_{b_0}^\star.
\]
Finally, the identity
\[
[\widetilde{\xi}, \widetilde{D}](u)
= (-1)^{|d||u|} \big(\widetilde{\xi}^2 - \Rr^\diamond_{\widetilde{b_0}}\big)(u)
\]
is equivalent to
\[
[\xi, \xi^*] = \xi^2 - \tfrac{1}{3} \Rr_{b_0}^\star,
\quad \text{and} \quad
b_0 \in \ker (\xi^* - \xi).
\]
Since  $(V, \mathfrak{b}, \xi^* - \xi, \xi, \tfrac{1}{3} b_0, 0)$ is a context of almost semi-direct product of left-symmetric superalgebras, it follows that 
$(V, \mathcal{H} := \mathfrak{b} \oplus V^*, \widetilde{D}, \widetilde{\xi}, \widetilde{b_0}, 0)$ is also a context of almost semi-direct product  of left-symmetric superalgebras if and only if the conditions \eqref{eq:claim2} are satisfied. This proves the lemma.
\end{proof}
\sssbegin{Remark}
Let $(V, \mathcal{H} := \mathfrak{b} \oplus V^*, \widetilde{D}, \widetilde{\xi}, \widetilde{b_0}, 0)$ be a context of an almost semi-direct product of left-symmetric superalgebras.  
The endomorphism defined by $
\Delta(u) := \xi^* - 2 \xi
$, for all $u\in \mathfrak{b}$, 
is a derivation on $(\mathfrak{b}, [\,,]_{\mathfrak{b}})$.  
\end{Remark}

We now define the even flat double extension of a flat orthosymplectic quasi-Frobenius Lie superalgebra by the one-dimensional Lie algebra.

\sssbegin{Theorem}\label{double-ex1}

Let $(\mathfrak{b}, \br_{\mathfrak{b}}, \omega_{\mathfrak{b}})$ be a flat orthosymplectic quasi-Frobenius Lie superalgebra, and let $\star_{\mathfrak{b}}$ be  the natural symplectic  product associated with it. 
Let $V = \mathbb{K}d$ be a one-dimensional vector space and $V^* = \mathbb{K}e$ be its dual. Assume that there exist an even linear map 
$\xi : \mathfrak{b} \to \mathfrak{b}$, 
and an element $b_{0} \in \mathfrak{b}_{\bar 0}$, satisfying the system \eqref{eq:claim2}. Define 
$
\g := \mathbb{K}d \oplus \mathfrak{b} \oplus \mathbb{K}e,
$
equipped with the Lie bracket and even bilinear form $\omega$ as follows:
\begin{equation}\label{Liebrackets1}
[d,u] = (\xi^* - 2\xi)(u) + \omega_{\mathfrak{b}}(b_0, u)e, \quad\text{ and }\quad
[u,v] = [u,v]_{\mathfrak{b}} + \omega_{\mathfrak{b}}\big((\xi + \xi^*)(u), v\big)e,
\end{equation}
for all $u,v \in \mathfrak{b}$ and
$$
\omega|_{\mathfrak{b} \times \mathfrak{b}} = \omega_{\mathfrak{b}},
\qquad \omega(e,d) = -\omega(d,e) = 1,
\qquad \omega(d,\mathfrak{b}) = \omega(e,\mathfrak{b}) = 0.
$$
Then, $(\g, \br, \omega)$ is a flat orthosymplectic quasi-Frobenius Lie superalgebra.

Furthermore,  the natural symplectic product $\star$  associated with $(\g,\br, \omega)$,  is given by
\begin{equation}\label{Produit1-13}
\begin{array}{lll} \displaystyle 
e \star u = u \star e =0,& e \star d = d \star e = 0,& \displaystyle  
d \star u = (\xi^* - \xi)(u) + \tfrac{1}{3}\omega_{\mathfrak{b}}(b_0, u)e, \\[0.3em]
\displaystyle  u \star d = \xi(u) - \tfrac{2}{3}\omega_{\mathfrak{b}}(b_0, u)e, &d \star d = \tfrac{1}{3} b_0, & 
u \star v = u \star_{\mathfrak{b}} v + \omega_{\mathfrak{b}}(\xi(u), v)e, 
\end{array}
\end{equation}
for all $ u,v \in \mathfrak{b}$.
\end{Theorem}
The flat orthosymplectic quasi-Frobenius Lie superalgebra $(\g, \br, \omega)$ is called the even flat double extension of $(\mathfrak{b}, \br_{\mathfrak{b}}, \omega_{\mathfrak{b}})$ by means of $(\xi, b_0).$
\begin{proof}
It is easy to check that $(\mathfrak{g},\br, \omega)$ is an orthosymplectic quasi-Frobenius Lie superalgebra.  Now, we show that it is  flat. From Eq. \eqref{cyclic13}, we derive  the natural symplectic product given by \eqref{Produit1-13}. Hence
\[
(u + \alpha e) \star (v + \beta e) = u \star_{\mathfrak{b}} v + \omega_{\mathfrak{b}}(\xi(u), v)e,\quad \text{for all $u,v \in \mathfrak{b}$ and $\alpha, \beta \in \mathbb{K}$.}
\]
 By Lemma \ref{Lemma1}, the first condition of system~\eqref{eq:claim2}, for all $u,v\in \mathfrak{b}$ 
\[
\xi([u,v]) = \mathrm{L}_u^{\star_\mathfrak{b}} (v)-(-1)^{|u||v|} \mathrm{L}_v^{\star_\mathfrak{b} }(u),
\]
implies that the above product defines on the  vector superspace  $\mathcal{H} = \mathfrak{b} \oplus \mathbb{K}e$ a left-symmetric structure which is a central extension of $\mathfrak{b}$ by $V$ by means of $\mu$, where 
$\mu(u,v) = \omega_{\mathfrak{b}}(\xi(u), v)
$, for all $u,v\in \mathfrak{b}$.

Now,  since $(\xi, b_0)$ satisfies the system \eqref{eq:claim2}, it follows from Lemma \ref{Lemma2}, that 
$
(V, \mathcal{H} := \mathfrak{b} \oplus V^*, \widetilde{D}, \widetilde{\xi}, \widetilde{b_0}, 0)
$
is a context of an almost semi-direct product of left-symmetric superalgebras. 
Hence, $(\mathfrak{g}, \star)$ is a left-symmetric superalgebra, 
which implies that $(\mathfrak{g}, [\, , \,], \omega)$ is a flat orthosymplectic quasi-Frobenius Lie superalgebra.
\end{proof}

\sssbegin{Lemma}\label{latrace}
Let $(\g, \br, \omega)$ be an orthosymplectic quasi-Frobenius nilpotent Lie superalgebra. If $a \in \g_{\bar{0}}$, then $\tr(\Rr^{\star}_a) = 0$.
\end{Lemma}

\begin{proof}
Let $a \in \g_{\bar{0}}$. Then, we have $
\omega(\Ll^{\star}_a(u), v) = -\,\omega(u, \Ll^{\star}_a(v)),$ for all  $u,v \in \g$. Choose a homogeneous basis (even $\mid$ odd)
$
\{X_1, Y_1, \ldots, X_n, Y_n \mid f_1, \ldots, f_m\}
$
of $\g$ such that 
\[
\omega(X_i, Y_i) = -\omega(Y_i, X_i) = -1, \quad \text{and} \quad \omega(f_i, f_i) = \epsilon_i, \ \text{where } \epsilon_i = \pm 1.
\]
We compute the (ordinary) trace of $\Ll^{\star}_a$:
\[
\begin{aligned}
\tr(\Ll^{\star}_a) &= \sum_{i=1}^n \omega(\Ll^{\star}_a(X_i), Y_i)
           - \sum_{i=1}^n \omega(\Ll^{\star}_a(Y_i), X_i)
           + \sum_{i=1}^m \epsilon_i \,\omega(\Ll^{\star}_a(f_i), f_i) \\
&= -\sum_{i=1}^n \omega(X_i, \Ll^{\star}_a(Y_i))
   + \sum_{i=1}^n \omega(Y_i, \Ll^{\star}_a(X_i))
   - \sum_{i=1}^m \epsilon_i \,\omega(f_i, \Ll^{\star}_a(f_i)) \\
&= \sum_{i=1}^n \omega(\Ll^{\star}_a(Y_i), X_i)
   - \sum_{i=1}^n \omega(\Ll^{\star}_a(X_i), Y_i)
   - \sum_{i=1}^m \epsilon_i \,\omega(\Ll^{\star}_a(f_i), f_i) \\
&= -\,\tr(\Ll^{\star}_a).
\end{aligned}
\]
Hence $\tr(\Ll^{\star}_a) = 0$.  
Since $\ad_a = \Ll^{\star}_a - \Rr^{\star}_a$, we have
\[
\tr(\ad_a) = \tr(\Ll^{\star}_a) - \tr(\Rr^{\star}_a) = -\tr(\Rr^{\star}_a).
\]
But $(\g, \br)$ is nilpotent, so $\tr(\ad_a) = 0$. Therefore $\tr(\Rr^{\star}_a) = 0$.
\end{proof}
\sssbegin{Lemma}\label{Lenilpotent}
Let $(\g, \br, \om)$ be a flat orthosymplectic quasi-Frobenius nilpotent Lie  superalgebra. If $a \in  \g_{\bar 0}$, then the right multiplication operator $\Rr_a$ is nilpotent.
\end{Lemma}

\begin{proof}
Since $(\g, \br, \om)$ is a flat, then according to relation \eqref{LR}, we have
\begin{equation}\label{R12}
\Rr^{\star}_{u \star v} - (-1)^{|u||v|} \Rr^{\star}_v \circ \Rr^{\star}_u = [\Ll^{\star}_u, \Rr^{\star}_v].
\end{equation}
Now, let $a \in \g_0$. We define inductively 
$
x_0 = a, \; x_{n+1} = \Ll^{\star}_{x_n}(a), \quad \text{for any } n \ge 0.
$ We shall show by induction that, for any $n \ge 1$, we have 
\begin{equation}\label{RR}
\Rr^{\star}_{x_n} = (\Rr^{\star}_a)^{\,n+1} + (\Rr^{\star}_a)^{\,n-2} \circ [\Ll^{\star}_a, \Rr^{\star}_a] + (\Rr^{\star}_a)^{\,n-3} \circ [\Ll^{\star}_{x_1}, \Rr^{\star}_a] + \dots + [\Ll^{\star}_{x_{n-1}}, \Rr^{\star}_a].
\end{equation}

The relation is clearly true for $n=1$. Assume it holds for some $n \ge 1$. Then, using Eq. \eqref{R12}, we have
\begin{equation*}
\begin{aligned}
\Rr^{\star}_{x_{n+1}} &= \Rr^{\star}_{x_n \star a} = \Rr^{\star}_a \circ \Rr^{\star}_{x_n} + [\Ll^{\star}_{x_n}, \Rr^{\star}_a] \\
&= \Rr^{\star}_a \circ \Big( (\Rr^{\star}_a)^{\,n+1} + (\Rr^{\star}_a)^{\,n-2} \circ [\Ll^{\star}_a, \Rr^{\star}_a] + (\Rr^{\star}_a)^{\,n-3} \circ [\Ll^{\star}_{x_1}, \Rr^{\star}_a] + \dots + [\Ll^{\star}_{x_{n-1}}, \Rr^{\star}_a] \Big)\\
& \quad + [\Ll^{\star}_{x_n}, \Rr^{\star}_a] \\
&= (\Rr^{\star}_a)^{\,n+2} + (\Rr^{\star}_a)^{\,n-1} \circ [\Ll^{\star}_a, \Rr^{\star}_a] + (\Rr^{\star}_a)^{\,n-2} \circ [\Ll^{\star}_{x_1}, \Rr^{\star}_a] + \dots + [\Ll^{\star}_{x_n}, \Rr^{\star}_a].
\end{aligned}
\end{equation*}
This completes the induction step.

 On the other hand, since $(\g, \br)$ is nilpotent,  Lemma \ref{latrace} implies that  $\tr(\Rr^{\star}_b)=0,$ for all $b\in \g_0$. Now, we put $b=x_{n+1} \in \mathfrak{g}_0$, we get  
 $$
\operatorname{tr}(\Rr^{\star}_{x_{n+1}}) = 0 \quad \text{and} \quad
\operatorname{tr}((\Rr^{\star}_a)^k \circ [\Ll^{\star}_a, \Rr^{\star}_a]) = 0 \quad \text{for all } n\, ,k \in \mathbb{N},
$$
since,  
 \begin{align*}
      \operatorname{tr}((\Rr^{\star}_a)^k \circ [\Ll^{\star}_a, \Rr^{\star}_a]) &=  \operatorname{tr}((\Rr^{\star}_a)^k \circ (\Ll^{\star}_a\circ\Rr^{\star}_a-\Rr^{\star}_a\circ \Ll^{\star}_a))\\
      & =\operatorname{tr}((\Rr^{\star}_a)^k \circ (\Ll^{\star}_a\circ\Rr^{\star}_a))-\operatorname{tr}((\Rr^{\star}_a)^{k+1}\circ \Ll^{\star}_a))\\&=\operatorname{tr}((\Rr^{\star}_a)^{k+1} \circ \Ll^{\star}_a)-\operatorname{tr}((\Rr^{\star}_a)^{k+1}\circ \Ll^{\star}_a)=0.\end{align*}
From Eq. \eqref{RR}, we deduce by induction that
$$
\operatorname{tr}((\Rr^{\star}_a)^{n+1}) = 0 \quad \text{for all } n \in \mathbb{N},
$$
and hence $\Rr^{\star}_a$ is nilpotent. 
\end{proof}
\sssbegin{Theorem}\label{nil1}
Let $(\mathfrak{b}, \br_\mathfrak{b}, \om_\mathfrak{b})$ be a flat orthosymplectic quasi-Frobenius Lie superalgebra and let $(\g, \br, \om)$ be the even flat double extension of $(\mathfrak{b}, \br_\mathfrak{b}, \om_\mathfrak{b})$ by means of $(\xi, b_0)$. If $(\mathfrak{b}, \br_\mathfrak{b})$ is nilpotent, then $(\g, \br)$ is nilpotent as well.
\end{Theorem}
\begin{proof}
First, we show that, for any $k \in \mathbb{N}^*$ and $a \in \mathfrak{b}_0$,
\begin{equation}\label{ee1}
\operatorname{tr}\big(\xi^k \circ \Rr_a^{\star_\mathfrak{b}}\big) = \operatorname{tr}\big(\Rr^{\star_\mathfrak{b}}_{\xi^k(a)}\big).
\end{equation}
Indeed, for any $a \in \mathfrak{b}_{\bar 0}$ and $u \in \mathfrak{b}$, according to system \eqref{eq:claim2}, we have
$$
\xi([u,a]) = u \star \xi(a) - a \star \xi(u),
$$
which implies $$
\xi \circ \Rr^{\star_\mathfrak{b}}_a - \xi \circ \Ll^{\star_\mathfrak{b}}_a = \Rr^{\star_\mathfrak{b}}_{\xi(a)} - \Ll^{\star_\mathfrak{b}}_a \circ \xi.
$$
Taking the trace gives $\operatorname{tr}(\xi \circ \Rr^{\star_\mathfrak{b}}_a) = \operatorname{tr}(\Rr^{\star_\mathfrak{b}}_{\xi(a)})$. Hence, Eq. \eqref{ee1} holds for $k=1$. Assume the property is true for some $k \ge 1$. Then for $k+1$, we have
$$
\xi^{k+1}([a,u]) = \xi^{k+1} \circ \Ll^{\star_\mathfrak{b}}_a(u) - \xi^{k+1} \circ \Rr^{\star_\mathfrak{b}}_a(u).
$$
On the other hand, for any $k\in \mathbb{N}$, we have 
$$
\xi^{k+1}([a,u]) = \xi^k \circ \xi([a,u]) = \xi^k \big(a \star \xi(u) - u \star \xi(a)\big) = \xi^k \circ \Ll^{\star_\mathfrak{b}}_a \circ \xi(u) - \xi^k \circ \Rr^{\star_\mathfrak{b}}_{\xi(a)}(u).
$$
Comparing the two expressions and taking the trace, we obtain
$$
\operatorname{tr}(\xi^{k+1} \circ \Rr^{\star_\mathfrak{b}}_a) = \operatorname{tr}(\xi^k \circ \Rr^{\star_\mathfrak{b}}_{\xi(a)}) = \operatorname{tr}(\Rr^{\star_\mathfrak{b}}_{\xi^{k+1}(a)}),
$$
which proves the property for $k+1$.


We set \(\Delta := \xi^* - 2\xi\). From \eqref{eq:claim2}, we have
$
[\xi, \Delta] = \xi^2 - \frac{1}{3} \Rr^{\star_\mathfrak{b}}_{b_0}.
$
By induction, one can show that for any \(k \in \mathbb{N}^*\), we have 
$$
[\xi^k, \Delta] = k \xi^{k+1} - \frac{1}{3} \sum_{p=0}^{k-1} \xi^p \circ \Rr^{\star_\mathfrak{b}}_{b_0} \circ \xi^{k-1-p}.
$$
It follows that 
$
\operatorname{tr}(\xi^{k+1}) = \frac{1}{3} \operatorname{tr}\big(\Rr^{\star_\mathfrak{b}}_{b_0} \circ \xi^{k-1}\big).
$
Then, from Eq. \eqref{ee1}, we deduce that
$
\operatorname{tr}(\xi^{k+1}) = \frac{1}{3} \operatorname{tr}\big(\Rr^{\star_\mathfrak{b}}_{\xi^{k-1}(b_0)}\big)$ for all  $k \in \mathbb{N}^*.
$ Since \((\mathfrak{b}, \br_\mathfrak{b})\) is nilpotent, Lemma \ref{Lenilpotent} implies that \(\operatorname{tr}(\Rr^{\star_\mathfrak{b}}_a) = 0\) for any \(a \in \mathfrak{b}_{\bar 0}\). Therefore, since $\xi^{k-1}(b_0)\in \mathfrak{b}_{\bar 0}$, then \(\operatorname{tr}(\xi^{k+1}) =\frac{1}{3} \operatorname{tr}\big(\Rr^{\star_\mathfrak{b}}_{\xi^{k-1}(b_0)}\big)= 0\) for all \(k \in \mathbb{N}^*\), which implies that \(\xi\) is nilpotent.

From the system \eqref{eq:claim2}, we have $
\xi^* \circ \xi = \frac{1}{3} (\Rr^{\star_\mathfrak{b}}_{b_0} + (\Rr^{\star_\mathfrak{b}}_{b_0})^*), $ and $ [\xi, \xi^*] = \xi^2 - \frac{1}{3} \Rr^{\star_\mathfrak{b}}_{b_0}.
$
It follows that
$
\xi \circ \xi^* = \xi^2 + \frac{1}{3} (\Rr^{\star_\mathfrak{b}}_{b_0})^*,$ and  $ (\xi^*)^2 = [\xi, \xi^*] + \frac{1}{3} (\Rr_{b_0}^{\star_\mathfrak{b}})^*$.
   
Replacing these into the expression
$\Delta^2 = (\xi^* - 2\xi)^2,$ 
we get
\[
\De^2 = 3 \xi^2 - 3 \xi^* \circ \xi = -3 \De \circ \xi - 3 \xi^2.
\]
By induction, one can show  that
$
\De^k = a_k  \De \circ \xi^{k-1} + b_k  \xi^k$ for all $k \ge 2$, where $a_k, b_k$ are some scalars. Since $\xi$ is nilpotent, it follows that $\De$ is nilpotent.

Now, if $(\mathfrak{g},\br, \omega)$ is an even flat double extension of $(\mathfrak{b}, \br_{\mathfrak{b}}, \omega_{\mathfrak{b}})$, the Lie bracket in $\mathfrak{g}$ is given by
$$
[d,u] = \De(u) + \omega_{\mathfrak{b}}(b_0, u)e, \qquad [u,v] = [u,v]_{\mathfrak{b}} + \omega_{\mathfrak{b}}\big((\xi+\xi^*)(u), v\big) e,
$$
for all $u,v \in \mathfrak{b}$.  By induction, one can show
that for any $k\in \mathbb{N}$, we have 
$$
\ad^k_d(u)=\De^k(u)+\omega_{\mathfrak{b}}(b_0, \De^{k-1}(u))e,\text{ and  }
\ad_u^k(v)=(\ad_u^\mathfrak{b})^k(v)+ \omega_{\mathfrak{b}}((\xi+\xi^*)(u), (\ad_u^\mathfrak{b})^{k-1}(v))e,
$$
for all $u,v\in \g$. Since $(\mathfrak{b}, \br_{\mathfrak{b}})$ is nilpotent and $\De$ is nilpotent, it follows that $(\mathfrak{g}, \br)$ is nilpotent. This completes the proof.
\end{proof}
\subsubsection{Odd flat double extensions of  orthosympletic Lie superalgebras}
We now introduce the notion of  {\it odd flat double extension} of flat orthosymplectic quasi-Frobenius Lie superalgebras by a purely odd one-dimensional vector superspace.

Let $(\mathfrak{b}, \br_\mathfrak{b}, \omega_\mathfrak{b})$ be a   flat orthosymplectic quasi-Frobenius Lie superalgebra, and denote by $\star$ the natural symplectic product  associated with it.  
Let $V = \K d$ be a purely odd one-dimensional vector superspace and $V^* = \K e$ be its dual.  

In order to introduce the process of the odd flat double extension, we consider an odd linear map  $
\xi : \mathfrak{b} \rightarrow \mathfrak{b},
$ and an even element $b_0 \in \mathfrak{b}_{\bar 0}$. We define the linear map
\begin{equation}\label{Bj} D(u) =- (-1)^{|u|}(\xi^* + \xi)(u)\quad \text{For all  $u\in {\mathfrak{b}}$.} \end{equation} 
We extend the two maps $D$ and $\xi$ to 
$
\widetilde{D}, \widetilde{\xi} : \mathfrak{b} \oplus V^* \longrightarrow \mathfrak{b} \oplus V^*, 
$
as follows: 
\begin{equation}
\widetilde{D}(u + \alpha e)
= -(-1)^{|u|}(\xi^* + \xi)(u)+ \, \omega(b_0, u)\, e, \quad \text{ and }\quad 
\widetilde{\xi}(u + \alpha e)
= \xi(u),
\label{eq:double-tilde1}
\end{equation}
for all $u \in \mathfrak{b}$ and $\alpha \in \K$. Put $
\widetilde{b_0}
=  b_0.$

\medskip

\sssbegin{Lemma}\label{Lemma3} 
The superspace $\mathcal{H} := \mathfrak{b} \oplus V^*$ equipped with the product
\[
(u + \alpha e) \diamond (v + \beta e)
:= u \star v +(-1)^{|v|} \omega(\xi(u), v)\, e,
\qquad \text{ for all } u,v \in \mathfrak{b},\; \alpha, \beta \in \K,
\]
is a left-symmetric superalgebra if and only if
\begin{equation}
\xi([u,v]_{\mathfrak b}) = \mathrm{L}_u^\star( \xi(v)) - (-1)^{|u||v|} \mathrm{L}_v^\star (\xi(u)), \quad \text{for all $u,v \in \mathfrak{b}$.}
\label{eq:xi-condition}
\end{equation}

\end{Lemma}

\begin{proof}
The algebra $(\mathcal{H}, \diamond)$ is a left-symmetric superalgebra if and only if the bilinear map 
\[
\mu(u,v) := \omega(\xi(u), v)
\]
satisfies condition \eqref{con-left}, which is equivalent to relation \eqref{eq:xi-condition}.\end{proof}

Now, let us also assume that $(V, \mathfrak{b}, D, \xi, b_0, 0)$ is a context of almost semi-direct product 
of left-symmetric superalgebras, where  $D(u)=-(-1)^{|u|}(\xi^* + \xi)(u)$, for all $u\in \mathfrak{b}$.

\medskip

\sssbegin{Lemma}\label{Lemma4}
$(V, \mathcal{H} := \mathfrak{b} \oplus V^*, \widetilde{D}, \widetilde{\xi}, \widetilde{b_0}, 0)$ is a context of almost semi-direct product 
of left-symmetric superalgebras if and only if 
\begin{equation}\label{eq:claim4}
\left\{
\begin{array}{rclrcl}
\xi([u,v]_{\mathfrak b}) & = &  \mathrm{L}_u^\star (\xi(v)) - (-1)^{|u||v|} \mathrm{L}_v^\star (\xi(u)),&\left[\xi , \xi^*\right]
& = &  \Rr^\star_{b_0}-3\xi^2,\\[2mm]
 D(u\star v)
& = & D(u)\star v+(-1)^{|u|}u\star D(v)+\xi(u)\star v,& 
\Ll^\star_{b_0} & = & -(\xi+\xi^*)^2 ,\\[2mm]
(2\xi+\xi^*)(b_0)&=&0,\end{array}
\right.
\end{equation}
for all $u, v \in \mathfrak{b}$. 
\end{Lemma}
\begin{proof}
Let $u, v \in \mathfrak{b}$ and $\alpha, \beta \in \mathbb{K}$. 
A straightforward computation shows that the equality
\[
\widetilde{\xi}([u,v]) = \mathrm{L}_u^\diamond (\widetilde{\xi}(v)) - (-1)^{|u||v|} \mathrm{L}_v^\diamond (\widetilde{\xi}(u))
\]
is equivalent to
\[
\xi([u,v]_{\mathfrak b}) = \mathrm{L}_u^\star (\xi(v)) - (-1)^{|u||v|} \mathrm{L}_v^\star (\xi(u)).
\]
Moreover, the relation
\[
\widetilde{D}(u) \diamond v - \widetilde{D}(u \diamond v)
= (-1)^{|d||u|} \big(\widetilde{\xi}(u) \diamond v - u \diamond \widetilde{D}(v)\big)
\]
is equivalent to
\[
-(-1)^{|u|}(\xi^* + \xi)(u) \star v 
+ (-1)^{|u|+|v|} (\xi^* + \xi)(u \star v)
= (-1)^{|u|} \xi(u) \star v 
+ (-1)^{|u|+|v|} u \star (\xi^* + \xi)(v),
\]
together with
\[
-(-1)^{|u|+|v|}\,\omega(\xi(\xi+\xi^*)(u), v)
- \omega(b_0, u \star v)
= (-1)^{|u|+|v|}\,\omega(\xi^2(u), v)
+ \omega((\xi+\xi^*)(\xi(u)), v).
\]
Since both $\xi \circ (\xi+\xi^*)$ and $\xi^2$ are even, and $\omega$ is even, it follows that when $|u| + |v| = \bar{1}$, the equality reduces to the trivial identity $0 = 0$. Consequently, the only nontrivial case occurs when $|u| + |v| = \bar{0}$, which yields the relation:
\[
-\omega(\xi(\xi+\xi^*)(u), v)
+ \omega(u \star b_0, v)
= \omega(\xi^2(u), v)
+ \omega((\xi+\xi^*)(\xi(u)), v).
\]
Since $\omega$ is nondegenerate, we deduce that
\[
\xi \circ \xi^* + \xi^* \circ \xi
= -3\xi^2 + \Rr^\star_{b_0}.
\]
The identity
\[
(\widetilde{\xi}\circ  \widetilde{D}- \widetilde{D}\circ \widetilde{\xi})(u)
= (-1)^{|d||u|} \big(\widetilde{\xi}^2 - \Rr^\diamond_{\widetilde{b_0}}\big)(u)
\]
is equivalent to
\[
\xi \circ \xi^* + \xi^* \circ \xi= -3\xi^2 + \Rr^\star_{b_0}, \text{ and  } \om(b_0, \xi(u))=(-1)^{|u|}\om(b_0, \xi(u)).
\]
Since both $b_0$ and $\om$ are even, and $\xi$ is odd,  we deduce that $ (-1)^{|u|}\om(b_0, \xi(u))=\om(b_0, \xi(u))$.

In addition, the relation
\[
(1 - (-1)^{|d|}) \big(\widetilde{D}^2 - \lambda \widetilde{D} - \Ll^\diamond_{b_0}\big) = 0,
\]
is equivalent to
\[
(2\xi+\xi^*)(b_0) = 0, 
\quad \text{and} \quad 
(\xi+\xi^*)\circ (\xi+\xi^*) = -\Ll^\star_{b_0}.
\]
Finally, the identity 
$
(\widetilde{\xi} - \widetilde{D})(\widetilde{b_0}) = 0
$
is equivalent to
$
(2\xi+\xi^*)(b_0) = 0.
$

Since $(V, \mathfrak{b},D, \xi^*,  b_0, 0)$ is a context of almost semi-direct product of left-symmetric superalgebras, it follows that 
$(V, \mathcal{H} := \mathfrak{b} \oplus V^*, \widetilde{D}, \widetilde{\xi}, \widetilde{b_0}, 0)$ is also a context of almost semi-direct product of left-symmetric superalgebras if and only if the conditions \eqref{eq:claim4} are satisfied. This proves the lemma.
\end{proof}

 Now, we are in a position to introduce the odd flat double extensions of   flat orthosymplectic quasi-Frobenius Lie  superalgebras. 


\sssbegin{Theorem}\label{double-ex2}
Let $(\mathfrak{b}, \br_{\mathfrak{b}}, \omega_{\mathfrak{b}})$ be a   flat orthosymplectic quasi-Frobenius Lie superalgebra, and let $\star_{\mathfrak{b}}$ denote the natural symplectic  product associated with it.
Let $V = \mathbb{K}d$ be a purely odd one-dimensional
vector superspace and $V^* = \mathbb{K}e$ its dual. Assume that there exist an odd homogeneous linear map 
$\xi : \mathfrak{b} \to \mathfrak{b}$, 
and an even element $b_0 \in \mathfrak{b}_0$, satisfying the system \eqref{eq:claim4}. Define the vector superspace
$
\g := \mathbb{K}d \oplus \mathfrak{b} \oplus \mathbb{K}e,
$
equipped with the Lie bracket and even bilinear form $\omega$ as follows: 
\begin{equation}\label{Liebrackets2}
\begin{array}{lcllcl}
 [d,u] & = & - (-1)^{|u|}(\xi^* + 2\xi)(u) + \omega_{\mathfrak{b}}(b_0, u)e, & [d,d]& = & 2b_0,\\[2mm] 
[u,v] &=& [u,v]_{\mathfrak{b}} + \big((-1)^{|v|}\omega_{\mathfrak{b}}(\xi(u), v) +(-1)^{|u|}\omega_{\mathfrak{b}}(\xi^*(u),  v)\big)e,\end{array}
\end{equation}
for all $u,v \in \mathfrak{b}$ and
$$
\omega|_{\mathfrak{b} \times \mathfrak{b}} = \omega_{\mathfrak{b}},
\qquad \omega(e,d) = \omega(d,e) = 1,
\qquad \omega(d,\mathfrak{b}) = \omega(e,\mathfrak{b}) = 0.
$$
Then, $(\g, \br, \omega)$ is a flat orthosymplectic quasi-Frobenius Lie superalgebra.

Furthermore,  the natural symplectic product $\star$  associated with $(\g,\br, \omega)$,  is given by
\begin{equation}\label{Produit3-13}
\begin{array}{lll}
e \star u = u \star e =0, & e \star d = d \star e = 0, & 
d \star u = -(-1)^{|u|}(\xi^* + \xi)(u) + \omega_{\mathfrak{b}}(b_0, u)e, \\[0.3em]
u \star d = \xi(u), & d \star d = b_0, & 
u \star v = u \star_{\mathfrak{b}} v + (-1)^{|v|}\omega_{\mathfrak{b}}(\xi(u), v)e,  
\end{array}
\end{equation}
for all $ u,v \in \mathfrak{b}$.
\end{Theorem}
The   flat orthosymplectic quasi-Frobenius Lie superalgebra $(\g, \br, \omega)$ is called the odd flat double extension of $(\mathfrak{b}, \br_{\mathfrak{b}}, \omega_{\mathfrak{b}})$ by means of $(\xi, b_0)$.
\begin{proof}
The proof follows the same lines as that of Theorem~\ref{double-ex1}, using Lemmas \ref{Lemma3} and \ref{Lemma4}.
\end{proof}
\sssbegin{Proposition}\label{nil2}
Let $(\mathfrak{b}, \br_\mathfrak{b}, \om_\mathfrak{b})$ be a   flat orthosymplectic quasi-Frobenius Lie superalgebra and let $(\g, \br, \om)$ be the odd flat double extensionn of $(\mathfrak{b}, \br_\mathfrak{b}, \om_\mathfrak{b})$ by means of $(\xi, b_0)$. If $(\mathfrak{b}, \br_\mathfrak{b})$ is nilpotent, then $(\g, \br)$ is nilpotent as well.
\end{Proposition}
\begin{proof}
Let us define the operator $\De$ by
\[
\De(u): = (-1)^{|u|}(2\xi + \xi^*)(u), \quad \text{for all $u\in \mathfrak{b}$}.
\]
Then
\[
\De^2(u) 
= (-1)^{|u|+|d|}( 2\xi+\xi^* )\big((-1)^{|u|}(2\xi+\xi^*)(u)\big)
= -\big(4\xi^2+2\xi\circ \xi^*+2\xi^*\circ \xi+(\xi^*)^2\big)(u).
\]
According to the system \eqref{eq:claim4}, we have 
\[
\xi \circ \xi^* + \xi^* \circ \xi+3\xi^2=\Rr^{\star_{\mathfrak{b}}}_{b_0}, \text{ and }
\Ll^{\star_{\mathfrak{b}}}_{b_0}=-(\xi^2+\xi\circ\xi^*+\xi^*\circ\xi+(\xi^2)^*).
\]
It follows that $\Delta^2=\Ll^{\star_{\mathfrak{b}}}_{b_0}-\Rr^{\star_{\mathfrak{b}}}_{b_0}=\ad_{b_0}^\mathfrak{b}$. Since $(\mathfrak{b}, [\, , \,]_{\mathfrak{b}})$ is nilpotent, it follows that $\Delta^2$ is nilpotent as well.

Now, if $(\mathfrak{g},\br, \omega)$ is an odd flat double extension of $(\mathfrak{b}, \br_{\mathfrak{b}}, \omega_{\mathfrak{b}})$, the Lie brackets in $\mathfrak{g}$ are given by
$$
[d,u] = \De(u) + \omega_{\mathfrak{b}}(b_0, u)e, \qquad [u,v] = [u,v]_{\mathfrak{b}}  + \big((-1)^{|v|}\omega_{\mathfrak{b}}(\xi(u), v) +(-1)^{|u|}\omega_{\mathfrak{b}}(\xi^*(u),  v)\big)e,
$$
for all $u,v \in \mathfrak{b}$. By induction, one can show
that for any $k\in \mathbb{N}$, we have 
$$
\ad^k_d(u)=\De^k(u)+\omega_{\mathfrak{b}}(b_0, D^{k-1}(u))e,$$
and 
$$
\ad_u^k(v)=(\ad_u^\mathfrak{b})^k(v)+ \big((-1)^{\al}\omega_{\mathfrak{b}}(\xi(u), (\ad_u^\mathfrak{b})^{k-1}(v)) +(-1)^{|u|}\omega_{\mathfrak{b}}(\xi^*(u),  (\ad_u^\mathfrak{b})^{k-1}(v))\big)e,
$$
where $\al=|(\ad_u^\mathfrak{b})^{k-1}(v)|$. Since $(\mathfrak{b}, \br_{\mathfrak{b}})$ is nilpotent, $\De$ is nilpotent and $e\in Z(\g)$, it follows that $(\mathfrak{g}, \br)$ is nilpotent as well.  This completes the proof.
\end{proof}
\subsubsection{Even Flat double extension of flat periplectic Lie superalgebras}
Following the approach of Theorem~\ref{double-ex1}, we now define the even flat double extension for  periplectic quasi-Frobenius Lie superalgebras. The proof, which relies on Lemmas \ref{Lemma1} and \ref{Lemma2}, yields the following result.
\sssbegin{Theorem}\label{double-ex3}

Let $(\mathfrak{b}, \br_{\mathfrak{b}}, \omega_{\mathfrak{b}})$ be a flat periplectic  quasi-Frobenius Lie superalgebra, and let $\star_{\mathfrak{b}}$ denote the natural symplectic  product associated with it.
Let $V = \mathbb{K}d$ be a one-dimensional vector space and  and $\Pi(V^*) = (\mathbb{K} e)_{\bar{1}}$, where $V^*$ its dual. Assume that there exist an even linear map 
$\xi : \mathfrak{b} \to \mathfrak{b}$, 
and an even element $b_0 \in \mathfrak{b}_{\bar 0}$, satisfying the system \eqref{eq:claim2}. Define the vector superspace
$
\g :=V \oplus \mathfrak{b} \oplus \Pi(V^*),
$
equipped with the Lie bracket and odd bilinear form $\omega$ as follows:
\begin{equation}\label{Liebrackets3}
[d,u] = (\xi^* - 2\xi)(u) + \omega_{\mathfrak{b}}(b_0, u)e, \quad\text{ and }\quad
[u,v] = [u,v]_{\mathfrak{b}} + \omega_{\mathfrak{b}}\big((\xi + \xi^*)(u), v\big)e,
\end{equation}
for all $u,v \in \mathfrak{b}$, and
$$
\omega|_{\mathfrak{b} \times \mathfrak{b}} = \omega_{\mathfrak{b}},
\qquad \omega(e,d) = -\omega(d,e) = 1,
\qquad \omega(d,\mathfrak{b}) = \omega(e,\mathfrak{b}) = 0.
$$
Then, $(\g, \br, \omega)$ is a flat periplectic  quasi-Frobenius Lie superalgebra.

Furthermore,  the natural symplectic product $\star$  associated with $(\g,\br, \omega)$,  is given by
\begin{equation}\label{Produit_odd1}
\begin{array}{lll} \displaystyle 
e \star u = u \star e =0, & e \star d = d \star e = 0,& \displaystyle  
d \star u = (\xi^* - \xi)(u) + \tfrac{1}{3}\omega_{\mathfrak{b}}(b_0, u)e, \\[0.3em]
\displaystyle  u \star d = \xi(u) - \tfrac{2}{3}\omega_{\mathfrak{b}}(b_0, u)e, &d \star d = \tfrac{1}{3} b_0,& 
u \star v = u \star_{\mathfrak{b}} v + \omega_{\mathfrak{b}}(\xi(u), v)e, 
\end{array}
\end{equation}
for all $ u,v \in \mathfrak{b}$.
\end{Theorem}
The flat periplectic quasi-Frobenius Lie superalgebra $(\g, \br, \omega)$ is called the even flat double extension of $(\mathfrak{b}, \br_{\mathfrak{b}}, \omega_{\mathfrak{b}})$ by means of $(\xi, b_0)$.

\begin{proof} The proof is similar to that of Theorem~\ref{double-ex1}, using Lemmas~\ref{Lemma1} and~\ref{Lemma2}.
\end{proof}
\subsubsection{Odd Flat double extension of flat periplectic Lie superalgebras}
We now introduce the notion of  odd flat double extension of flat periplectic quasi-Frobenius Lie superalgebras by a purely odd one-dimensional vector superspace.

Let $(\mathfrak{b}, \br_\mathfrak{b}, \omega_\mathfrak{b})$ be a   flat periplectic quasi-Frobenius Lie superalgebra, and let $\star_{\mathfrak{b}}$ denote the natural symplectic  product associated with it. Let $V = \K d$ be a purely odd one-dimensional vector superspace and $\Pi(V^*) = \K e$, where $V^*$ its dual.  

To introduce the process of the odd flat double extension by $V$, we consider an odd linear map $
\xi : \mathfrak{b} \rightarrow \mathfrak{b},$ and an even element $b_0 \in \mathfrak{b}_{\bar 0}$. We introduce the map $
D : \mathfrak{b} \rightarrow \mathfrak{b}$  defined as 
\[
D(u)=(-1)^{|u|}( \xi^* - \xi)(u), \quad \text{for all $u\in \mathfrak{b}$}.
\]
These two maps can be extended as 
$
\widetilde{D}, \widetilde{\xi} : \mathfrak{b} \oplus \Pi(V^*) \rightarrow \mathfrak{b} \oplus \Pi(V^*),
$ as follows:  
\begin{equation}
\widetilde{D}(u + \alpha e)
= (-1)^{|u|}(\xi^* - \xi)(u)+ \, \omega(b_0, u)\, e, \quad 
\widetilde{\xi}(u + \alpha e)
= \xi(u), 
\label{eq:double-tilde2}
\end{equation}
for all $u \in \mathfrak{b}$ and $\alpha \in \K$. Put $\widetilde{b_0}
=  -b_0$. 

\sssbegin{Lemma} \label{Lemma5}
The superspace $\mathcal{H} := \mathfrak{b} \oplus \Pi(V^*)$ equipped with the product
\[
(u + \alpha e) \diamond (v + \beta e)
:= u \star v +(-1)^{|v|} \omega(\xi(u), v)\, e,
\qquad \text{for all } u,v \in \mathfrak{b},\; \alpha, \beta \in \K,
\]
is a left-symmetric superalgebra if and only if
\begin{equation}
\xi([u,v]_{\mathfrak b}) = \mathrm{L}_u^\star( \xi(v)) - (-1)^{|u||v|} \mathrm{L}_v^\star (\xi(u)),
\text{ for all $u,v \in \mathfrak{b}$.} \label{eq:xi-condition2}\end{equation}

\end{Lemma}

\begin{proof}
The proof is similar to that of Lemma \ref{Lemma1}. \end{proof}

Now, let us also assume that $(V, \mathfrak{b}, D, \xi,- b_0, 0)$ is a context of almost semi-direct product 
of left-symmetric superalgebras, where  $D(u)=(-1)^{|u|}(\xi^* - \xi)(u)$, for all $u\in \mathfrak{b}$.

\sssbegin{Lemma}\label{Lemma6}
$(V, \mathcal{H} := \mathfrak{b} \oplus \Pi(V^*), \widetilde{D}, \widetilde{\xi}, \widetilde{b_0}, 0)$ is a context of almost semi-direct product 
of left-symmetric superalgebras if and only if 
\begin{equation}\label{eq:claim6}
\left\{
\begin{array}{rclrcl}
\xi([u,v]_{\mathfrak b})&=&  \mathrm{L}_u^\star (\xi(v)) - (-1)^{|u||v|} \mathrm{L}_v^\star( \xi(u)),& [\xi ,\xi^*]  &=&  3\xi^2 + \Rr^\star_{b_0},\\[2mm] D(u\star v)&=& D(u)\star v-(-1)^{|u|}u\star D(v)+
(-1)^{|u|} \xi(u) \star v, &
\Ll^\star_{b_0}  &=& (\xi-\xi^*)^2,  \\[2mm] 
2\xi(b_0)&=&\xi^*(b_0),\end{array}
\right.
\end{equation}
for all $u, v \in \mathfrak{b}$. 
\end{Lemma}

\begin{proof}
Let $u, v \in \mathfrak{b}$ and $\alpha, \beta \in \mathbb{K}$. 
A straightforward computation shows that the equality
\[
\widetilde{\xi}([u,v]) = \mathrm{L}_u^\diamond (\widetilde{\xi}(v)) - (-1)^{|u||v|} \mathrm{L}_v^\diamond (\widetilde{\xi}(u))
\]
is equivalent to
\[
\xi([u,v]_{\mathfrak b}) = \mathrm{L}_u^\star (\xi(v)) - (-1)^{|u||v|} \mathrm{L}_v^\star (\xi(u)).
\]
Moreover, the relation
\[
\widetilde{D}(u) \diamond v - \widetilde{D}(u \diamond v)
= (-1)^{|d||u|} \big(\widetilde{\xi}(u) \diamond v - u \diamond \widetilde{D}(v)\big)
\]
is equivalent to
\[
(-1)^{|u|}(\xi^* - \xi)(u) \star v 
- (-1)^{|u|+|v|} (\xi^* - \xi)(u \star v)
= (-1)^{|u|} \xi(u) \star v 
- (-1)^{|u|+|v|} u \star (\xi^*- \xi)(v),
\]
together with
\[
(-1)^{|u|+|v|}\,\omega(\xi(\xi^*-\xi)(u), v)
- \omega(b_0, u \star v)
= (-1)^{|u|+|v|}\,\omega(\xi^2(u), v)
+ \omega((\xi^*-\xi)(\xi(u)), v).
\]
Since both $\xi \circ (\xi^*-\xi)$ and $\xi^2$ are even, and $\omega$ is odd, it follows that when $|u| + |v| = \bar{0}$, the equality reduces to the trivial identity $0 = 0$. Consequently, the only nontrivial case occurs when $|u| + |v| = \bar{1}$, leading to the following relation:
\[
-\omega(\xi(\xi^*-\xi)(u), v)
+ \omega(u \star b_0, v)
= -\omega(\xi^2(u), v)
+ \omega((\xi^*-\xi)(\xi(u)), v).
\]
Since $\omega$ is nondegenerate, we deduce that
\[
\xi \circ \xi^* + \xi^* \circ \xi
= 3\xi^2 + \Rr^\star_{b_0}.
\]
The identity
\[
(\widetilde{\xi}\circ \widetilde{D}-\widetilde{D}\circ \widetilde{\xi})(u)
= (-1)^{|d||u|} \big(\widetilde{\xi}^2 - \Rr^\diamond_{\widetilde{b_0}}\big)(u),
\]
is equivalent to
\[
\xi \circ \xi^* + \xi^* \circ \xi= 3\xi^2 + \Rr^\star_{b_0}, \text{ and  } \om(b_0, \xi(u))=(-1)^{|u|}\om(b_0, \xi(u)).
\]
Since  $b_0$ is even, $\om$ is odd, and $\xi$ is odd, then the relation $\om(b_0, \xi(u))=(-1)^{|u|}\om(b_0, \xi(u))$ is therefore always satisfied. In addition, the relation
\[
(1 - (-1)^{|d|}) \big(\widetilde{D}^2 - \lambda \widetilde{D} - \Ll^\diamond_{b_0}\big) = 0,
\]
is equivalent to
\[
(2\xi-\xi^*)(b_0) = 0, 
\quad \text{and} \quad 
(\xi^*-\xi)\circ (\xi^*-\xi) = \Ll_{b_0}.
\]
Finally, the identity 
$
(\widetilde{\xi} - \widetilde{D})(\widetilde{b_0}) = 0
$
is equivalent to
$
(2\xi-\xi^*)(b_0) = 0.
$

Since $(V, \mathfrak{b},D, \xi^*,  b_0, 0)$ is a context of almost semi-direct product of left-symmetric superalgebras, it follows that 
$(V, \mathcal{H} := \mathfrak{b} \oplus \Pi(V^*), \widetilde{D}, \widetilde{\xi}, \widetilde{b_0}, 0)$ is also a context of almost semi-direct product of left-symmetric superalgebras if and only if the conditions \eqref{eq:claim4} are satisfied. This proves the Lemma.
\end{proof}

\sssbegin{Theorem}\label{double-ex4}
Let $(\mathfrak{b}, \br_{\mathfrak{b}}, \omega_{\mathfrak{b}})$ be a   flat periplectic quasi-Frobenius Lie superalgebra, and let $\star_{\mathfrak{b}}$ denote the natural symplectic  product associated with it. Let $V=V_{\bar{1}} = (\mathbb{K} d)_{\bar{1}}$ be a purely odd one-dimensional vector superspace  and $\Pi(V^*) = \mathbb{K} e$ where $V^*$ its dual. Assume that there exist an odd homogeneous linear map 
$\xi : \mathfrak{b} \to \mathfrak{b}$, 
and an even element $b_0 \in \mathfrak{b}_0$, satisfying the system \eqref{eq:claim6}. Define the vector superspace
$
\g := \mathbb{K}d \oplus \mathfrak{b} \oplus \mathbb{K}e,
$
equipped with the Lie bracket and odd bilinear form $\omega$ given as follows: 
\begin{equation}\label{Liebrackets4}
\begin{array}{lcllcl}
[d,u] & = & (-1)^{|u|}(\xi^* -2\xi)(u) + \omega_{\mathfrak{b}}(b_0, u)e, & [d,d] & = &-2b_0,\\ [2mm]
[u,v] & = &  [u,v]_{\mathfrak{b}} + \big((-1)^{|v|}\omega_{\mathfrak{b}}(\xi(u), v) +(-1)^{|u|}\omega_{\mathfrak{b}}(\xi^*(u),  v)\big)e,\end{array}
\end{equation}
for all $u,v \in \mathfrak{b}$ and
$$
\omega|_{\mathfrak{b} \times \mathfrak{b}} = \omega_{\mathfrak{b}},
\qquad \omega(e,d) = -\omega(d,e) = 1,
\qquad \omega(d,\mathfrak{b}) = \omega(e,\mathfrak{b}) = 0.
$$
Then $(\g, \br, \omega)$ is a   flat periplectic  quasi-Frobenius Lie superalgebra.

Furthermore,  the natural symplectic product associated with $(\g,\br, \omega)$,  is given by
\begin{equation}\label{Produit2-13}
\begin{array}{lll}
e \star u = u \star e =0, &  e \star d = d \star e = 0, & 
d \star u = (-1)^{|u|}(\xi^* - \xi)(u) + \omega_{\mathfrak{b}}(b_0, u)e, \\[0.3em]
u \star d = \xi(u) , & d \star d =- b_0, & 
u \star v = u \star_{\mathfrak{b}} v + (-1)^{|v|}\omega_{\mathfrak{b}}(\xi(u), v)e. 
\end{array}
\end{equation}
for all $ u,v \in \mathfrak{b}$.
\end{Theorem}
The   flat periplectic  quasi-Frobenius Lie superalgebra $(\g, \br, \omega)$ is called the odd flat double extension of $(\mathfrak{b}, \br_{\mathfrak{b}}, \omega_{\mathfrak{b}})$ by means of $(\xi, b_0)$.

\begin{proof}
The structure of this proof is analogous to that of Theorem~\ref{double-ex1}; the key difference is the use of Lemmas~\ref{Lemma5} and~\ref{Lemma6} instead of the lemmas used there.
\end{proof}

\section{Flat quasi-Frobenius Lie superalgebras are nilpotent}\label{flat-nil-sec}
 In this section, we show that all non-abelian flat quasi-Frobenius Lie superalgebras have a degenerate center.
Moreover, every such Lie superalgebra can be obtained via flat double extensions. In particular, every flat orthosymplectic quasi-Frobenius Lie superalgebra can be constructed by a sequence of flat double extensions starting from an an abelian Lie superalgebra,  while flat periplectic quasi-Frobenius Lie superalgebras can be obtained by a sequence of flat double extensions starting from the trivial algebra $\{0\}$. Furthermore, we show that every flat orthosymplectic quasi-Frobenius Lie superalgebra is not only solvable but nilpotent

\ssbegin{Theorem}\label{Centre}
Let $(\g, \br, \om)$ be a flat quasi-Frobenius Lie superalgebra. Then $Z(\g)\neq\{0\}$. Moreover, if $(\g, [\; , \;])$ is a non abelian Lie superalgebra, then  $Z(\g)$ is degenerate. 
\end{Theorem}
\begin{proof} We distinguish two cases.

\underline{The case where $\omega$ is odd.}  If $(\mathfrak{g}, \br)$ is abelian, then it is clear that $Z(\mathfrak{g}) \neq 0.$  If $(\mathfrak{g}, \br)$ is non-abelian then according to Theorem~\ref{nilpotent-odd}, $(\g, \br)$ is nilpotent. It follows that   $ Z(\g)\cap [\g, \g] \neq \{0\},$ and since $Z(\g)\subseteq [\g, \g]^\bot$, it follows that $Z(\g)\cap Z(\g)^\bot \neq \{0\}$ and we are done. 

\underline{The case where $\omega$ is even}. If $(\mathfrak{g}, \br)$ is abelian, then it is clear that $Z(\mathfrak{g}) \neq 0.$ 
Now, we assume that  $(\mathfrak{g}, \br)$ is non-abelian. We need the following two lemmas.

\medskip
\ssbegin{Lemma}
\label{Claim 1.}
The following implication holds:
\begin{equation}\label{imp}
Z(\g)\cap[\g, \g] = \{0\} \;\Longrightarrow\; Z([\g,\g]) \cap N_\ell(\g, \star) = \{0\}.
\end{equation}
\end{Lemma}

\begin{proof} Suppose that $Z(\g)\cap[\g, \g] = \{0\}$ and let $u \in Z([\g,\g]) \cap N_\ell(\g, \star)$.  
Then, for any $w \in \g$ and any $v \in [\g,\g]$, and since $\Ll^\star_u = 0$ together with \eqref{LR}, we have $\ad_u^* = \ad_u$. Hence,
\[
0 = \omega([u,v],w) = (-1)^{|u||v|}\omega(v,[u,w]).
\]
It follows that $[u,w] \in [\g,\g]^\perp$. But we have seen in Prop.~\ref{sided} that $[u,w] \in N_\ell(\g, \star)$, thus
$
[u,w] \in [\g,\g]^\perp \cap N_\ell(\g, \star) = Z(\g) \text{ by Prop. \ref{orth}}$, and consequently $ [u,w] \in Z(\g)\cap[\g, \g] = \{0\} $
which shows that $u \in Z(\g)\cap[\g, \g]$, and hence $u = 0$. \hfill $\blacksquare$
\noqed
\end{proof}

With this observation in mind, let us prove our result by contradiction.  
According to Prop.~\ref{solvable-even}, $(\g, \br)$ is solvable and $Z(\g)\cap[\g, \g] = \{0\}$.  
We consider the sequence of vector subspaces defined by
\[
C_0 = [\g,\g] \cap [\g,\g]^\perp, 
\qquad 
C_k = [[\g,\g], C_{k-1}], \quad k \geq 1.
\]
By virtue of Prop. ~\ref{degenerate}, we have $C_0 \neq \{0\}$, and according to Prop.~\ref{sided}, 
$C_k \subset N_\ell(\g, \star)$ for any $k \geq 1$.
\medskip

\ssbegin{Lemma}\label{Claim 2.}
For any $k \geq 1$, if $C_k = \{0\}$ then $C_{k-1} = \{0\}$. 
\end{Lemma}
\begin{proof}
Suppose that $C_k = \{0\}$. Then, for any $u \in C_{k-1}$, any $v \in [\g,\g]$, and any $w \in \g$, we have
\[
0 = \omega([u,v],w) = \pm (-1)^{|u||v|}\omega(v,[u,w]),
\]
since $\ad_u^* = -\ad_u$ if $k=1$ and $\ad_u^* = \ad_u$ if $k \geq 2$. Thus, $[u,w] \in [\g,\g]^\perp \cap N_\ell(\g, \star) = Z(\g)$ and hence $[u,w] \in Z(\g)\cap [\g, \g]$ which shows that $u = 0$.  
This completes the proof of the lemma.
\hfill $\blacksquare$
\noqed
\end{proof}

\medskip

Since $C_0 \neq \{0\}$, we deduce that for all $k \geq 1$, $C_k \neq \{0\}$. But $[\g,\g]$ is nilpotent, and therefore there exists $k_0$ such that $C_{k_0} \subseteq Z([\g,\g])$.  
We also have $C_{k_0} \subseteq N_\ell(\g, \star)$, according to \eqref{imp}, we obtain that $Z(\g)\cap [\g, \g]\neq{0}$. Since $[\g,\g]\subseteq  Z(\g)^\perp$, we deduce that
$
Z(\g)\cap [\g, \g]\subseteq Z(\g)\cap Z(\g)^\perp \neq \{0\}
$.\end{proof}


\ssbegin{Remark}\label{exemple-dim2}
Every flat quasi-Frobenius Lie superalgebra of dimension $2$ is abelian. Indeed, Theorem~\ref{Centre} implies  $Z(\mathfrak{g}) \neq 0$. Hence, there exists a non-zero element $a \in Z(\mathfrak{g})_{\bar{i}}$.

If $\omega$ is odd, then $\dim(\mathfrak{g}_{\bar{0}})=\dim( \mathfrak{g}_{\bar{1}})$ by Lemma \ref{BM}.  Assume $a \in Z(\mathfrak{g})_{\bar{0}}$. We can find an element $b\in \g_{\bar 1}$ such that $\omega(a, b) = 1$ and  $[b, b] = \alpha a$, where $\alpha \in \mathbb{K}$.    
As $\omega$ is closed, we obtain
\[
3\,\omega([b, b], b) = 0 \quad \Rightarrow \quad 3\alpha = 0,
\]
which implies $\alpha = 0$, and therefore $b$ is also central. Hence,  $(\mathfrak{g}, \br)$ is abelian. In the case where \(a \in Z(\mathfrak{g})_{\bar{1}}\), any element $b\in \g_{\bar 0}$ commutes with $a$ and is therefore central. 

If $\omega$ is even, then  $(\mathfrak{g}, [\; , \;])$ is either purely even or purely odd, as $\omega_{|_{\fg_{\bar 0}}}$ is also nondegenerate and therefore exists only if $\fg_{\bar 0}$ is even-dimensional. In the case where $(\mathfrak{g}, [\; , \;])$ is purely odd, the Lie superalgebra is abelian by definition. Let us assume that $(\mathfrak{g}, [\; , \;])$ is purely even. Any 2-dimensional Lie algebra with a non-zero center must be abelian.  
\end{Remark}

\ssbegin{Proposition}\label{reduit}  
Let $(\g, [\; ,\;], \omega)$ be a   flat quasi-Frobenius Lie superalgebra.  
Let $0\not = a\in Z(\g) \cap Z(\g)^{\perp}$,  and denote  
$
I := \mathbb{K} a,$ and   $\mathfrak{b} := I^{\perp} / I.
$
Let $\pi_{\mathfrak{b}} : I^{\perp} \to \mathfrak{b}$ be the canonical projection. We have
\begin{enumerate}
    \item[$(i)$] $I$ and $I^\perp$ are totally isotropic graded two-sided ideals of $(\g, \star)$. 
    \item[$(ii)$] The quotient $\mathfrak{b} = I^{\perp} / I$ admits a canonical structure of a   flat  quasi-Frobenius Lie superalgebra, defined by
    \[
        [\pi_{\mathfrak{b}}(u), \pi_{\mathfrak{b}}(v)]_{\mathfrak{b}} = \pi_{\mathfrak{b}}([u,v]) \quad \text{and} \quad 
        \omega_{\mathfrak{b}}(\pi_{\mathfrak{b}}(u), \pi_{\mathfrak{b}}(v)) = \omega(u, v), \quad \text{ for all $u,v \in I^{\perp}$. }
    \]

\end{enumerate}
\end{Proposition}

\begin{proof}
 
Let us proof Part (i). Since $Z(\g)$ is a homogeneous ideal, the element $a$ can be chosen to be homogeneous. Since $a \in Z(\g)$, we have $a \star u = u \star a = 0$ for all $u \in \g$, which implies that $I$ is a totally isotropic graded two-sided ideal of $(\g, \star)$.

    Now, let $v \in I^{\perp}$ and $w \in \g$. We have
    \[
        \omega(w \star v, a) = (-1)^{|w||v|} \omega(v, w \star a) = 0, 
        \qquad 
        \omega(v \star w, a) = (-1)^{|w||v|} \omega(w, v \star a) = 0,
    \]
    hence $w \star v,\, v \star w \in I^{\perp}$.
    
    If $\omega$ is even, then $I^{\perp}$ is a graded two-sided ideal of $(\g, \star)$.

If $\omega$ is odd, we distinguish two cases.  
If $a$ is even, and since $\omega$ is odd, we have $\g_{\bar{0}} \subseteq I^{\perp}$, 
and therefore $(I^{\perp})_{\bar{0}} = \g_{\bar{0}}$. 
If $a$ is odd, then $\g_{\bar{1}} \subseteq  I^{\perp}$, and hence 
$(I^{\perp})_{\bar{1}} = \g_{\bar{1}}$. 
In both cases, we conclude that $I^{\perp}$ is also a graded two-sided ideal of $(\g, \star)$.
 
    Since $I^{\perp}$ is a graded two-sided ideal of $(\g, \star)$, the quotient $\mathfrak{b} = I^{\perp} / I$ naturally inherits the structure of a Lie superalgebra and of a quasi-Frobenius. Indeed,  
    for $u, v \in I^{\perp}$, we define
    \[
        [\pi_{\mathfrak{b}}(u), \pi_{\mathfrak{b}}(v)]_{\mathfrak{b}} = \pi_{\mathfrak{b}}([u,v]), 
        \qquad 
        \omega_{\mathfrak{b}}(\pi_{\mathfrak{b}}(u), \pi_{\mathfrak{b}}(v)) = \omega(u, v).
    \]
    These operations are well-defined since $I$ is totally isotropic, i.e., $\omega(u,a)=0$ for all $u \in I^{\perp}$ and $a \in I$.

    Moreover, since $(\g, [\; ,\;], \om)$ is  quasi-Frobenius then $(\mathfrak{b}, [\; ,\;]_{\mathfrak{b}}, \omega_{\mathfrak{b}})$ is quasi-Frobenius as well. Furthermore, the product defined on $\mathfrak{b}$ is given by 
	\begin{align*}
\omega_{\mathfrak{b}}\left(\pi_{\mathfrak{b}}(u)\star_\mathfrak{b}\pi_{\mathfrak{b}}(v),\pi_{\mathfrak{b}}(w)\right)&=\frac{1}{3}\left(\om_{\mathfrak{b}}\left([\pi_{\mathfrak{b}}(u),\pi_{\mathfrak{b}}(v)],\pi_{\mathfrak{b}}(w)\right)+\om_{\mathfrak{b}}\left([\pi_{\mathfrak{b}}(u),\pi_{\mathfrak{b}}(w)],\pi_{\mathfrak{b}}(v)\right)\right)\\&=\frac{1}{3}\left(\om\left([u,v],\w\right)+\om\left([u,w],v\right)\right)\\&=\om(u\star v,w)
	\end{align*}

Thus, we deduce that $(\mathfrak{b}, \star_\mathfrak{b})$ is a  left-symmetric superalgebra and hence $(\mathfrak{b}, [\; ,\;]_{\mathfrak{b}}, \omega_{\mathfrak{b}})$ is a flat quasi-Frobenius Lie superalgebra.
\end{proof}

\ssbegin{Definition}
The Lie superalgebra $\mathfrak{b}$ obtained in Prop.~\ref{reduit} is called the \emph{flat quasi-Frobenius Lie superalgebra deduced from $I^{\perp}$ by means of $I$}.
\end{Definition}

We now prove the converse of Theorems \ref{double-ex1} and \ref{double-ex2}. 

\ssbegin{Theorem}\label{central}
Let $(\g, [\;,\;], \omega)$ be a   flat orthosymplectic quasi-Frobenius non-abelian Lie superalgebra.
Then $(\g, [\;,\;], \omega)$ is either:
\begin{enumerate}
    \item[$(i)$] an even flat double extension of a   flat orthosymplectic quasi-Frobenius Lie superalgebra $(\mathfrak{b}, [\;,\;]_{\mathfrak{b}}, \omega_{\mathfrak{b}})$ by means of $(\xi, b_0)$, or
    \item[$(ii)$] an odd flat double extension of a   flat orthosymplectic quasi-Frobenius Lie superalgebra $(\mathfrak{b}, [\;,\;]_{\mathfrak{b}}, \omega_{\mathfrak{b}})$ by means of $(\xi, b_0)$.
\end{enumerate}
\end{Theorem}
\begin{proof}
Theorem \ref{Centre} implies \( Z(\g) \cap Z(\g)^\perp \neq 0 \).
By Prop. \ref{reduit}, there exists a totally isotropic two-sided graded ideal \( I = \K e \) of \( (\g, \star) \), and its orthogonal \( I^\perp \) is also a two-sided graded ideal of \( (\g, \star) \).  We now distinguish two cases.

\underline{The case where  \( e \) is even.} Since \( \omega \) is even and non-degenerate, there exists an element \( d \in \g_{\bar{0}} \setminus \{0\} \) such that
$
\omega(e, d) = -\omega(d, e) = 1.
$

  Since $I\subseteq I^\perp$, there exists a subsuperspace $\mathfrak{B}$ such that  $I^\perp = I \oplus \mathfrak{B}$ and  the restriction  $\omega_{\mathfrak{B}}=\omega|_{\mathfrak{B}\times \mathfrak{B}}$  
 is non-degenerate.  Let us write
$
\g = \mathbb{K} e \oplus {\mathfrak B} \oplus \mathbb{K} d,
$
and since \( I^\perp = I \oplus {\mathfrak B} \) is a two-sided graded ideal of \( (\g, \star) \), for any  elements \( u, v \in {\mathfrak B} \), we have
\[
u \star v = u \star_{\mathfrak B} v + \mu(u, v)e,
\]
where \( \mu : {\mathfrak B} \times {\mathfrak B} \to \K \) is an even bilinear form, and 
\( \star_{\mathfrak B} : {\mathfrak B} \times {\mathfrak B} \to {\mathfrak B} \) is an even bilinear map.  

It follows that \( ({\mathfrak B}, \star_{\mathfrak B}) \) is a left-symmetric superalgebra,
\( \mu \) satisfies the compatibility condition \eqref{con-left},
and the restriction \( \omega_{\mathfrak B} = \omega|_{{\mathfrak B} \times {\mathfrak B}} \) is  orthosymplectic and closed on \( ({\mathfrak B}, \br_{\mathfrak B}) \). Moreover, the canonical projection
$
\pi : {\mathfrak B} \rightarrow I^\perp / I=\mathfrak{b}$, $ \pi(u) = \pi_{\mathfrak{b}}(u)$,  for all $  u \in {\mathfrak B},$ is an isomorphism of left-symmetric superalgebras. Thus, we can identify  \( {\mathfrak B} \cong I^\perp / I=\mathfrak{b} \), and hence
$
\g = \mathbb{K} e \oplus \mathfrak{b} \oplus \mathbb{K} d,
$
where \( (\mathfrak{b}, \br_{\mathfrak{b}}, \omega_{\mathfrak b}) = (I^\perp / I, \br_{\mathfrak{b}}, \om_\mathfrak{b}) \)
is a   flat orthosymplectic quasi-Frobenius Lie superalgebra.

Since \( I^\perp =  \mathbb{K} e \oplus \mathfrak{b} \) is a two-sided graded ideal of  \( (\g, \star) \), it is also a graded ideal of \( (\g, [\; , \;]) \).
Then, for all \( u,v \in \mathfrak{b} \), the Lie brackets on \( \g \) are given by
\[
[d, u] = D(u) + T(u)e, \qquad [u, v] = [u, v]_{\mathfrak b} + \big(\mu(u,v) - (-1)^{|u||v|}\mu(v,u)\big)e,
\]
where \( D \in \mathrm{End}(\mathfrak{b})_{\bar{0}} \), and \( T : \mathfrak{b} \to \K \) is an even linear map.  

Since \( \omega \) is even and non-degenerate, there exists \( b_0 \in \mathfrak{b}_{\bar{0}} \)
such that \( T(u) = \omega_\mathfrak{b}(b_0, u) \).
Moreover, the product associated with \( (\g, [\, ,\,], \omega) \) is given by
\[
\begin{aligned}
u \star v &= u \star_\mathfrak{b} v + \mu(u,v)e, \;
&u \star d = \xi(u) + f(u)e, \\
d \star u &= \rho(u) + g(u)e, \;
& d \star d = \lambda d + c_0 + \alpha e,
\end{aligned}
\]
where \( \xi, \rho \in \mathrm{End}(\mathfrak{b})_{\bar{0}} \),  \( f, g : \mathfrak{b} \to \K \) are even linear maps, and $\al, \la\in \K$.

Using the natural symplectic product  \eqref{cyclic13}, for any $u,v\in \mathfrak{b}$ we have
\[
\mu(u,v) = \omega(u \star v, d) = -(-1)^{|u||v|}\omega(v, u \star d)
= \omega_\mathfrak{b}(\xi(u), v),
\]
and hence
\[
\omega_\mathfrak{b}(\xi(u), v)=\om(u\star v, d) =\tfrac{1}{3}\om([u, v],d)+\tfrac{1}{3}\om([u,d], v)= \tfrac{1}{3}\omega_\mathfrak{b}((\xi + \xi^* - D)(u), v),
\]
which gives \( D = \xi^* - 2\xi \).
Similarly,
\[
g(u) = \omega(d \star u, d) = \tfrac{1}{3}\omega_\mathfrak{b}(b_0, u),
\]
and since \( [d, u] = d \star u - u \star d \), we obtain \( \rho(u) = (\xi^* - \xi)(u) \)
and \( f(u) = -\tfrac{2}{3}\omega_\mathfrak{b}(b_0, u) \).

Moreover,
\[
\omega(d \star d, e) = -\omega(d, d \star e) = \lambda = 0,
\quad
\omega(d \star d, d) =\alpha= \om([d,d], d)  = 0,
\]
and
\[
\omega_\mathfrak{b}(c_0, u) = \omega(d \star d, u)
= -\omega(d, d \star u)
= \tfrac{1}{3}\omega_\mathfrak{b}(b_0, u),
\]
hence \( c_0 = \tfrac{1}{3}b_0 \).

Therefore, the Lie bracket and the associated product take the simplified form:
\[
[d, u] = (\xi^* - 2\xi)(u) + \omega_\mathfrak{b}(b_0, u)e, 
\qquad
[u, v] = [u, v]_\mathfrak{b} + \omega_\mathfrak{b}((\xi + \xi^*)(u), v)e,
\]
and
\[
\begin{aligned}
u \star v &= u \star_\mathfrak{b} v + \omega_\mathfrak{b}(\xi(u), v)e,\; \quad\quad\quad u \star d = \xi(u) -\tfrac{2}{3} \omega_\mathfrak{b}(b_0, u)e, \\
d \star u &= (\xi^* - \xi)(u) +\tfrac{1}{3} \omega_\mathfrak{b}(b_0, u)e, \quad\;\;  d \star d = \tfrac{1}{3}b_0.
\end{aligned}
\]

Now, it is easy to show that \( (\g, \star) \) is a left-symmetric superalgebra if and only if the pair \( (\xi, b_0) \) satisfies the system \eqref{eq:claim2}. Hence, \( (\g, [\; ,\;], \omega) \) is the even flat double extension of the   flat orthosymplectic quasi-Frobenius Lie superalgebra \( (\mathfrak{b}, [\; ,\;]_{\mathfrak{b}}, \omega_{\mathfrak{b}}) \) by means of \( (\xi, b_0) \).

\underline{ The case where \( e \) is odd.}  Since \( \omega \) is even and non-degenerate, there exists an element \( d \in {\mathfrak b}_{\bar{1}} \setminus \{0\} \) such that
$
\omega(e, d) = \omega(d, e) = 1.
$

$
$
$
$
Since $I\subseteq I^\perp$, there exists a subsuperspace $\mathfrak{b}$ such that  $I^\perp = I \oplus \mathfrak{B}$ and  the restriction  $\omega_{\mathfrak{B}}=\omega|_{\mathfrak{B}\times \mathfrak{B}}$  
 is non-degenerate.  Let us write
$
\g = \mathbb{K} e \oplus {\mathfrak B} \oplus \mathbb{K} d,
$
and since \( I^\perp = I \oplus {\mathfrak B} \) is a two-sided graded ideal of \( (\g, \star) \), for any  elements \( u, v \in {\mathfrak B} \), we have
\[
u \star v = u \star_{\mathfrak B} v + \mu(u, v)e,
\]
where \( \mu : {\mathfrak B} \times {\mathfrak B} \to \K \) is an even bilinear form, and 
\( \star_{\mathfrak B} : {\mathfrak B} \times {\mathfrak B} \to {\mathfrak B} \) is an even bilinear map.  

It follows that \( ({\mathfrak B}, \star_{\mathfrak B}) \) is a left-symmetric superalgebra,
\( \mu \) satisfies the compatibility condition \eqref{con-left},
and the restriction \( \omega_{\mathfrak B} = \omega|_{{\mathfrak B} \times {\mathfrak B}} \) is  orthosymplectic and closed on \( ({\mathfrak B}, \br_{\mathfrak B}) \). Moreover, the canonical projection
$
\pi : {\mathfrak B} \rightarrow I^\perp / I=\mathfrak{b}$, $ \pi(u) = \pi_{\mathfrak{b}}(u)$,  for all $  u \in {\mathfrak B},$ is an isomorphism of left-symmetric superalgebras. Thus, we can identify  \( {\mathfrak B} \cong I^\perp / I=\mathfrak{b} \), and hence
$
\g = \mathbb{K} e \oplus \mathfrak{b} \oplus \mathbb{K} d,
$
where \( (\mathfrak{b}, \br_{\mathfrak{b}}, \omega_{\mathfrak b}) = (I^\perp / I, \br_{\mathfrak{b}}, \om_\mathfrak{b}) \)
is a   flat orthosymplectic quasi-Frobenius Lie superalgebra.

Since \( I^\perp =  \mathbb{K} e \oplus \mathfrak{b} \) is a two-sided graded ideal of  \( (\g, \star) \), it is also a graded ideal of \( (\g, \br) \).
Then, for all \( u,v \in \mathfrak{b} \), the Lie brackets on \( \g \) are given by
\[
[d,d]=c_0,\quad [d, u] = D(u) + T(u)e, \qquad [u, v] = [u, v]_B + \big(\mu(u,v) - (-1)^{|u||v|}\mu(v,u)\big)e,
\]
where \( D \in \mathrm{End}(\mathfrak{b})_{\bar{1}} \), and \( T : \mathfrak{b} \to \K \) is an  even linear map.  

Since \( \omega \) is odd and non-degenerate, there exists \( b_0 \in \mathfrak{b}_{\bar{0}} \)
such that \( T(u) = \omega_\mathfrak{b}(b_0, u) \).
Moreover, the product associated with \( (\g, [\, ,\,], \omega) \) is given by
\[
\begin{aligned}
u \star v &= u \star_\mathfrak{b} v + \mu(u,v)e, \;
u \star d = \xi(u) + f(u)e, \\
d \star u &= \rho(u) + g(u)e, 
\quad\quad d \star d =  a_0,
\end{aligned}
\]
where \( \xi, \rho \in \mathrm{End}(\mathfrak{b})_{\bar{1}} \), and \( f, g : \mathfrak{b} \to \K \) are even linear maps.

Using the natural symplectic product  \eqref{cyclic13}, for any $u,v\in \mathfrak{b}$, we have
\[
\mu(u,v) = \omega(u \star v, d) =- (-1)^{|u||v|}\omega(v, u \star d)
= (-1)^{|v|}\omega_\mathfrak{b}(\xi(u), v),
\]
and hence
\begin{align*}
(-1)^{|v|}\,\omega_{\mathfrak{b}}(\xi(u), v)
&= \omega(u \star v, d)
= \tfrac{1}{3}\om([u, v], w) + \tfrac{1}{3}(-1)^{|v|}\omega([u, d], v) \\
&= \tfrac{1}{3}(-1)^{|v|}\omega_{\mathfrak{b}}(\xi(u), v)
+ \tfrac{1}{3}(-1)^{|u|}\omega_{\mathfrak{b}}(\xi^{*}(u), v)
- \tfrac{1}{3}(-1)^{|u|+|v|}\omega_{\mathfrak{b}}(D(u), v).
\end{align*}
Hence,
\[
\omega_{\mathfrak{b}}(D(u), v)
= -2(-1)^{|u|}\omega_{\mathfrak{b}}(\xi(u), v)
+ (-1)^{|u|+|v|}\omega_{\mathfrak{b}}(\xi^{*}(u), v).
\]
Since both $\xi$ and $D$ are odd, and $\omega$ is even, it follows that when $|u| + |v| = \bar{0}$, the equality reduces to the trivial identity $0 = 0$. Consequently, the nontrivial case occurs when $|u| + |v| = \bar{1}$, leading to the following relation
\[
\omega_{\mathfrak{b}}(D(u), v)
=- (-1)^{|u|}\omega_{\mathfrak{b}}((2\xi+\xi^*)(u), v).
\]
It follows that $D(u)=- (-1)^{|u|}(2\xi+\xi^*)(u)$. 
Similarly,
\[
f(u) = \omega(u \star d, d) = 0,
\]
and since \( [d, u] = d \star u - u \star d \), we obtain $ g(u)=\om_\mathfrak{b}(b_0, u)$ and $\rho(u) = (-1)^{|u|}(\xi^* + \xi)(u)$.

Moreover,
\[
\omega_\mathfrak{b}(a_0, u) = \omega(d \star d, u)
= \omega(d, d \star u)
= \omega_\mathfrak{b}(b_0, u),
\]
hence \( a_0 = b_0 \) and since $[d,d]=2d\star d=c_0=2 b_0$.

Therefore, the Lie bracket and the associated product take the simplified form:
\begin{equation*}
\begin{array}{lcllcl}
[d, u] & = &  -(-1)^{|u|}(  2\xi+\xi^*)(u) + \omega_\mathfrak{b}(b_0, u)e,& [d,d] & = &2 b_0,\\[2mm]
[u, v] & = & [u, v]_\mathfrak{b} + ((-1)^{|v|}\omega_\mathfrak{b}(\xi(u), v) +(-1)^{|u|}\omega_\mathfrak{b}( \xi^*(u), v))e,
\end{array}
\end{equation*}
and
\[
\begin{aligned}
u \star v &= u \star_\mathfrak{b} v + (-1)^{|v|}\omega_\mathfrak{b}(\xi(u), v)e,\quad\quad
u \star d = \xi(u), \\
d \star u &= -(-1)^{|u|}(\xi^* + \xi)(u) + \omega_\mathfrak{b}(b_0, u)e, \;
d \star d = b_0.
\end{aligned}
\]

Now, it is easy to show that \( (\g, \star) \) is a left-symmetric superalgebra if and only if the pair \( (\xi, b_0) \) satisfies the system \eqref{eq:claim4}. Hence, \( (\g, [\; ,\;], \omega) \) is the odd flat double extension of the   flat orthosymplectic quasi-Frobenius Lie superalgebra \( (\mathfrak{b}, [\; ,\;]_{\mathfrak{b}}, \omega_{\mathfrak{b}}) \) by means of \( (\xi, b_0) \).\end{proof}

\ssbegin{Corollary}\label{start-abel}
An orthosymplectic quasi-Frobenius  Lie superalgebra $(\g, \br, \om)$ is flat if and only if it can be obtained by a sequence 
of even or odd flat double extensions of flat orthosymplectic quasi-Frobenius Lie superalgebras, starting from the abelian superalgebra.
\end{Corollary}
\begin{proof}
If $(\g, \br)$ is non-abelian Lie superalgebra, then according to Theorem \ref{central}, $(\g, \br)$ is even or odd flat double extension of a flat orthosymplectic quasi-Frobenius Lie superalgebra
$(\mathfrak{b}_1,\br_{\mathfrak{b}_1},  \om_{\mathfrak{b}_1})$ by means of  $(\xi_1, b_1)$. If $(\mathfrak{b}_1,\br_{\mathfrak{b}_1})$ is a non-abelian Lie superalgebra, then it is itself an even or odd flat double extension of a flat orthosymplectic quasi-Frobenius Lie superalgebra
$(\mathfrak{b}_2,\br_{\mathfrak{b}_2},  \om_{\mathfrak{b}_2})$ by means of  $(\xi_2, b_2)$. This process has to stop since $\g$ is finite-dimensional. 
\end{proof}
The following theorem offers a sharpening of Theorem \ref{solvable-even}:
\ssbegin{Theorem}\label{ref-nilp}
Every flat orthosymplectic quasi-Frobenius Lie superalgebra $(\g, \br, \omega)$ is nilpotent.
\end{Theorem}
\begin{proof}

We proceed by induction on $\dim \g$. If $\dim \g = 1$ or $\dim \g = 2$, the result is trivial since in this case the superalgebra is abelian, and hence nilpotent (see Remark~\ref{exemple-dim2}).

Assume that the statement holds for every flat orthosymplectic quasi-Frobenius Lie superalgebra of total dimension $\leq n$, and let us prove it for dimension $n+1$. 

Let $(\g, \br, \omega)$ be a flat orthosymplectic quasi-Frobenius Lie superalgebra of total dimension $n+1$. According to Theorem~\ref{central}, $(\g, \br, \omega)$ is either an even or an odd flat double extension of a flat orthosymplectic quasi-Frobenius Lie superalgebra $(\mathfrak{b}, \br_{\mathfrak{b}}, \omega_{\mathfrak{b}})$ of total dimension $n-1$. 

By the induction hypothesis, $(\mathfrak{b}, \br_{\mathfrak{b}})$ is nilpotent. Then, by Props. ~\ref{nil1} and~\ref{nil2}, we deduce that $(\g, \br, \omega)$ is nilpotent as well.
\end{proof}

Now, we give the converse of Theorem \ref{double-ex3} and Theorem \ref{double-ex4}. 

\sssbegin{Theorem}\label{central1}
Let $(\g, [\;,\;], \omega)$ be a   flat periplectic  quasi-Frobenius Lie superalgebra.
Then $(\g, [\;,\;], \omega)$ is either:
\begin{enumerate}
    \item[$(i)$] an even flat double extension of a   flat periplectic  quasi-Frobenius Lie superalgebra $(\mathfrak{b}, [\;,\;]_{\mathfrak{b}}, \omega_{\mathfrak{b}})$ by means of $(\xi, b_0)$, or
    \item[$(ii)$] an odd flat double extension of a   flat periplectic  quasi-Frobenius Lie superalgebra $(\mathfrak{b}, [\;,\;]_{\mathfrak{b}}, \omega_{\mathfrak{b}})$ by means of $(\xi, b_0)$.
\end{enumerate}
\end{Theorem}
\begin{proof}
The proof is similar to that of Theorem \ref{central}.
\end{proof}

\sssbegin{Corollary}\label{ext-triv}
A periplectic quasi-Frobenius  Lie superalgebra $(\g, \br, \om)$ is flat if and only if it can be obtained by a sequence 
of even or odd flat double extensions of flat periplectic quasi-Frobenius Lie superalgebras, starting from the trival superalgebra $\{0\}$.
\end{Corollary}
\begin{proof}
 Theorem \ref{central1} implies that $(\g, \br)$ is an even or odd flat double extension of a flat  periplectic quasi-Frobenius Lie superalgebra
$(\mathfrak{b}_1,\br_{\mathfrak{b}_1},  \om_{\mathfrak{b}_1})$ by means of  $(\xi_1, b_1)$. The latter is itself an even or an odd flat double extension of a flat  periplectic quasi-Frobenius Lie superalgebra
$(\mathfrak{b}_2,\br_{\mathfrak{b}_2},  \om_{\mathfrak{b}_2})$ by means of  $(\xi_2, b_2)$. Since the superdimension of $\g$ is finite and even (see Lemma \ref{BM}), there exists $k\in \mathbb{N}^*$ such that $(\mathfrak{b}_k,\br_{\mathfrak{b}_k})$ is trivial Lie superalgebra, i.e,. $\mathfrak{b}_k=\{0\}$.
\end{proof}
\section{Classification of flat  quasi-Frobenius Lie superalgebras of low dimensions}\label{classification}
Here, we provide a classification of flat quasi-Frobenuis Lie superalgebras of dimension $\leq 5$. This classification is carried out over a field of characteristic zero that is algebraically closed. It is interesting to investigate the classification problem over any field. 
\subsection{Four-dimensional flat quasi-Frobenius Lie superalgebras}
Here, we classify flat quasi-Frobenius Lie superalgebras of total dimension four.  
Theorems~\ref{central} and~\ref{central1} imply that every four-dimensional flat quasi-Frobenius Lie superalgebra $(\mathfrak{g},\br, \omega)$ is an odd or even flat double extension of a two-dimensional abelian one $(\mathfrak{b},\br_{\mathfrak{b}}, \omega_{\mathfrak{b}})$, defined by the pair $(\xi, b_0)$.  
Indeed, every two-dimensional quasi-Frobenius Lie superalgebra is abelian (see Remark~\ref{exemple-dim2}).
\sssbegin{Definition} Let $(\mathfrak{b}, [\; , \;]_{\mathfrak{b}}, \omega_{\mathfrak{b}})$ be a flat quasi-Frobenius Lie superalgebra. If $\omega_{\mathfrak{b}}$ is even, we say that a pair $(\xi, b_0) \in \operatorname{End}(\mathfrak{b}) \times \mathfrak{b}$ is \emph{even admissible} (resp. \emph{odd admissible}) if $(\xi, b_0)$ satisfies System~\eqref{eq:claim2} (resp. System~\eqref{eq:claim4}). If $\omega_{\mathfrak{b}}$ is odd, we say that a pair $(\xi, b_0) \in \operatorname{End}(\mathfrak{b}) \times \mathfrak{b}$ is \emph{even admissible} (resp. \emph{odd admissible}) if $(\xi, b_0)$ satisfies System~\eqref{eq:claim2} (resp. System~\eqref{eq:claim6}). \end{Definition}

\sssbegin{Remark}
Every flat orthosymplectic quasi-Frobenius Lie superalgebra of dimension three is abelian. 
Indeed, we suppose the Lie superalgbera is not abelian. Since $\omega$ is even, the even part $\mathfrak{g}_0$ has even dimension. It follows from Theorem \ref{central} that this  $(\mathfrak{g}, [\; , \;], \omega)$ is an even flat double extension of a   one-dimensional flat orthosymplectic quasi-Frobenius Lie superalgebra $(\mathfrak{b}, \br_\mathfrak{b},  \omega_{\mathfrak{b}})$ by means of  $(\xi, b_0)$. Since $\om_\mathfrak{b}$ is even, it follows that  $\mathfrak{b}_{\bar{0}}=0$. Moreover, as both $\xi$ and $b_0$ are even, it follows that $\xi = 0$ and $b_0 = 0$. Hence, from \eqref{Liebrackets1}, all brackets vanish, and therefore $(\mathfrak{g}, [\, , \,], \omega)$ is abelian, a contradiction.  
\end{Remark}

\sssbegin{Proposition}\label{solution admissible dim2}
Let $(\mathfrak{b}, [\; , \;]_{\mathfrak b},\omega_{\mathfrak{b}})$ be a two-dimensional abelian orthosymplectic quasi-Frobenius Lie superalgebra.
\begin{itemize}
    \item[(i)] The pair $(\xi, b_0)$ is \emph{even admissible} if and only if  there exists a basis $\mathbb{B} = \{e_1, e_2\}$ of $\mathfrak{b}$ in which both elements are even \textup{(}resp. both odd\textup{)}, such that  the following conditions are satisfied:
    \begin{itemize}
        \item If $\{e_1, e_2\}$ are both even, then we either have
\[
\xi = 0, \quad b_0 = \alpha e_1 + \beta e_2, \quad \omega = e_1^* \wedge e_2^*, \quad \alpha, \beta \in \mathbb{K}, \text{ or }
\]
  \[
        \xi = 
        \begin{pmatrix}
        0 & a \\
        0 & 0
        \end{pmatrix}, \quad
        b_0 = \alpha e_1, \quad
        \omega = e_1^* \wedge e_2^*, \quad \text{ where $\alpha \in \mathbb{K}$ and $a \neq 0$.}
        \]
       \\
 \item If $\{e_1, e_2\}$ are both odd, then we either have
\[
\xi = 0, \quad b_0 =0, \quad \omega = e_1^* \wedge e_1^* + \epsilon\, e_2^*\wedge e_2^*, \quad \epsilon = \pm 1.
\]
or \textup{(}where  $a \neq 0$\textup{)}
        \[
        \xi = 
        \begin{pmatrix}
        0 & a \\
        0 & 0
        \end{pmatrix}, \quad
        b_0 = 0, \quad
        \omega = -e_1^* \wedge e_2^*.
        \]

    \end{itemize}

    \item[(ii)] The pair $(\xi, b_0)$ is \emph{odd admissible} if and only if there exists a basis $\mathbb{B} = \{e_1, e_2\}$ of $\mathfrak{b}$ in which both elements are even \textup{(}resp. both odd\textup{)}, such \\
\begin{itemize} \item If $\{e_1, e_2\}$ are both even, we have
\[
\xi = 0, \quad b_0 = \alpha e_1 + \beta e_2, \quad \omega = e_1^* \wedge e_2^*, \quad \alpha, \beta \in \mathbb{K}.
\]
\item If $\{e_1, e_2\}$ are both odd, we have
\[
\xi = 0, \quad b_0 =0, \quad \omega = e_1^* \wedge e_1^* + \epsilon\, e_2^* \wedge e_2^*, \quad \epsilon = \pm 1.
\]
\end{itemize}       
\end{itemize}
\end{Proposition} 
\begin{proof}
Let us only prove Part (i). If the pair $(\xi, b_0)$ is even admissible, then 
\[
\xi^* \circ \xi = 0, \qquad \xi \circ \xi^* = \xi^2, \qquad \text{and} \qquad b_0 \in \ker(\xi^* - \xi).
\]
Hence, $\xi$ is nilpotent, and therefore $\xi = 0$ or $\xi^2 = 0$ with $\xi \neq 0$.  
Since $\omega$ is even and $\mathfrak{b}$ is two-dimensional, we must have $\mathfrak{b} = \mathfrak{b}_{\bar 0}$ or $\mathfrak{b} = \mathfrak{b}_{\bar 1}$. 

\medskip
\noindent
\underline{The case where $\mathfrak{b} = \mathfrak{b}_{\bar 0}$, and $\xi \neq 0$.}  
We can choose a symplectic basis $\{e_1, e_2\}$ of $\mathfrak{b}$ such that 
\[
\omega = e_1 \wedge e_2, \quad 
\xi =
\begin{pmatrix}
0 & a \\
0 & 0
\end{pmatrix},
\text{ and }
b_0 = \alpha e_1 + \beta e_2,
\]
where $a \neq 0$ and $\alpha, \beta \in \mathbb{K}$.  
The condition $(\xi^* - \xi)(b_0) = 0$ implies $\beta = 0$.

\medskip
\noindent
\underline{The case where $\mathfrak{b} = \mathfrak{b}_{\bar 1}$, and $\xi \neq 0$.}  
We can choose a basis $\{f_1, f_2\}$ of $(\mathfrak{b}, \omega)$ such that
\[
\omega =- f_1^* \wedge f_2^*, \quad 
\xi =
\begin{pmatrix}
0 & a \\
0 & 0
\end{pmatrix},
\text{ and }
b_0 = \alpha f_1 + \beta f_2,
\]
where $a \neq 0$ and $\alpha, \beta \in \mathbb{K}$.  
Again, $(\xi^* - \xi)(b_0) = 0$ implies $\beta = 0$. If $\{f_1, f_2\}$ are both odd, then $b_0$ being even forces $b_0 = 0$.

\noindent \underline{The case where $\xi = 0$.} The conclusion follows immediately. 
\end{proof}

\sssbegin{Proposition} \label{solution admissible dim2 odd}
Let $(\mathfrak{b}, [\; , \;]_{\mathfrak b}, \omega_{\mathfrak{b}})$ be the two-dimensional periplectic quasi-Frobenius abelian Lie superalgebra.
\begin{itemize}
    \item[(i)] The pair $(\xi, b_0)$ is \emph{even admissible} if and only if there exists a basis $\mathbb{B} = \{e_1\mid  f_1\}$ of $\mathfrak{b}$ such that the following conditions are satisfied:
 $$
         \xi = 0,\; b_0 = \alpha e_1 ,\quad \text{and}\quad \omega = e_1^* \wedge f_1^*,\quad \text{where}\quad \alpha \in \mathbb{K}.$$

   \item[(ii)]  
The pair $(\xi, b_0)$ is \emph{odd admissible} if and only if there exists 
a homogeneous basis $\mathbb{B}=\{e_1 \mid f_1\}$ of $\mathfrak{b}$ such that  
$\omega = e_1^{*} \wedge f_1^{*}$, and one of the following holds:

\[
\begin{array}{lll}
\textup{(1)} & \xi = 0, 
& b_0 = \alpha\, e_1, \quad \al \in \mathbb{K},
\\[4pt]
\textup{(2)} & 
\xi = 
\begin{pmatrix}
0 & 0 \\[2pt]
a & 0
\end{pmatrix},
& b_0 = 0, \qquad a \neq 0,
\\[10pt]
\textup{(3)} &
\xi = 
\begin{pmatrix}
0 & a \\[2pt]
0 & 0
\end{pmatrix},
& b_0 = \alpha\, e_1, \qquad a \neq 0,\; \al \in \mathbb{K}.
\end{array}
\]

\end{itemize}
\end{Proposition}

\begin{proof}
Let us only prove Part (i). If the pair $(\xi, b_0)$ is even admissible, then 
\[
\xi^* \circ \xi = 0, \qquad 
\xi \circ \xi^* = \xi^2, \qquad 
b_0 \in \ker(\xi^* - \xi).
\]
Hence, $\xi$ is nilpotent.   Since $\omega$ is odd and $\xi$ is even, there exists a homogeneous basis $\{e_1\mid f_1\}$ of $\mathfrak{b}$ such that
\[
\omega = e_1^* \wedge f_1^*,
\qquad 
\xi(e_1) = a\, e_1, \qquad 
\xi(f_1) = b\, f_1, \quad b_0=\al e_1, 
\]
for some scalars $a, b, \al \in \mathbb{K}$.  Since $\xi$ is nilpotent, we must have $a = b = 0$, and therefore $\xi = 0$.
\end{proof}

\sssbegin{Proposition}\label{cls-dim4}
Any four-dimensional flat quasi-Frobenius non-abelian Lie superalgebra is isomorphic to one of the following Lie superalgebras:

\begin{itemize}
\item[$(i)$]  $(\mathbb{K}\oplus \h_3, \br, \om)$ \textup{(}in the basis $\{x_1, x_2, x_3\}$\textup{)}:   
$$[x_1, x_2]=x_3, \quad \om= x_1^*\wedge x_4^*+x_2^* \wedge x_3^*,$$

\item[$(ii)$]  $(\g^2, \br, \om)$ \textup{(}in the basis $\{x_1, x_2\mid  y_1, y_2\}$\textup{)}:
$$
	[x_1,y_1]=y_2,\quad [y_1, y_1]= x_2\mbox{ and }\om=2x_1^*\wedge x_2^*-y_1^*\wedge y_2^*,
	$$
\item[$(iii)$]  $(\g^3, \br, \om)$ \textup{(}in the basis $\{x_1, x_2\mid  y_1, y_2\}$\textup{)}:
$$[x_1, y_1]=y_2, \quad \om= x_1^*\wedge y_2^*+x_2^* \wedge y_1^*,$$
\item[$(iv)$]  $(\g^4, \br, \om)$ \textup{(}in the basis $\{x_1, x_2\mid  y_1, y_2\}$\textup{)}:
$$[y_1, y_1]=x_1,\quad [y_1, y_2]=x_2, \quad \om= -2x_1^*\wedge y_2^*+x_2^* \wedge y_1^*,$$
\end{itemize}
\end{Proposition}

\begin{proof}
	Let $(\mathfrak{g},[\;,\;], \om)$ be  a flat  quasi-Frobenius non-abelian Lie superalgebra of total dimension $4$. Theorems~\ref{central} and~\ref{central1}, implies that such a superalgebra is an even or an odd flat double extension of the two-dimensional abelian quasi-Frobenius Lie superalgebra $(\mathfrak{b},[\;,\;]_{\mathfrak{b}}, \omega_{\mathfrak{b}})$ by means of  $(\xi, b_0)$. Following Props. \ref{solution admissible dim2} and \ref{solution admissible dim2 odd}, we argue as follows. 
    
\underline{The case where  $\omega_{\mathfrak{b}}$ is even.} We have several cases:
\begin{itemize}
\item If $(\xi, b_0)$ is an even admissible pair, we distinguish several cases: \begin{itemize} \item 
If $\xi=0$,  $b_0=\alpha e_{1}+\beta e_2$, and $\om=e_1^*\wedge e^*_2$, where $\alpha,\beta\in \mathbb{K}$,  then Eq. \eqref{Liebrackets1}, implies that there exists a basis $\{e,d,e_{1},e_2\}$ of $\mathfrak{g}$ such that the Lie brackets are given by 
		\begin{equation*}
	\left[d,e_1\right]=-\beta e,\;\left[d,e_2\right]=\alpha e, \quad \om= e^*\wedge d^*+e_1^* \wedge e_2^*.
		\end{equation*} 
		Since $(\g, \br)$ is non-abelian, then $\left (\alpha,\beta \right )\neq(0,0)$. If for instance $\beta\neq0$, then we choose   $$(x_1,x_2,x_3,x_4)=\left (-\frac{1}{\beta}e_1, d, -e, -\beta e_2-\alpha e_1\right ).$$ 
		It follows that $(\mathfrak{g},[\; , \;], \om)$ is symplectomorphic to $(\mathbb{K}\oplus \h_3, \br, \om)$.
\item 
If $\xi\neq 0$ and $b_0=\alpha e_{1}$, where $\alpha\in \mathbb{K}$,  then Eq. \eqref{Liebrackets1}, implies that there exists a basis $\{e,d,e_{1},e_2\}$ of $\mathfrak{g}$ such that the Lie brackets are given by 
$$\left[d,e_2\right]=\alpha e-3ae_1, \,\, \text{ where } a\neq 0.$$
		Put  $$(x_1,x_2,x_3,x_4)=\left(3a d+\alpha e_2,d+\frac{\alpha+1}{3a}e_2, \alpha e-3ae_1,-\frac{\alpha+1}{3a}e+e_1\right).$$ 
		In this new basis, the Lie brackets and the form are reduced to
		\[
		[x_1,x_2]=x_3\mbox{ and }\om=x_1^*\wedge x_4^*+x_2^*\wedge x_3^*,
		\]
		which implies that $(\mathfrak{g},[\; , \;], \om)$ is symplectomorphic to $(\mathbb{K}\oplus \h_3, \br, \om)$.

       \item  If $\xi\neq 0$ and $b_0=0$,  then Eq.  \eqref{Liebrackets1} implies that there exists a basis $\{e,d\mid  e_{1},e_2\}$  of $\mathfrak{g}$ such that the Lie brackets are given by 
$$\left[d,e_2\right]=-ae_1,\quad [e_2, e_2]=2 a e, \,\, a\neq 0.$$
		Put   $$(x_1, x_2,  y_1,y_2)=\left(-\frac{1}{a}d, 2a e, - e_2, -e_1\right).$$ 
		In this new basis, the Lie brackets and the form are reduced to 
		\[
		[x_1,y_1]=y_2,\quad [y_1,y_1]=x_2\mbox{ and }\om=2x_1^*\wedge x_2^*-y_1^*\wedge y_2^*,
		\]
		which implies that $(\mathfrak{g},[\;, \;], \om)$ is symplectomorphic to $(\g^2, \br, \om)$.

        \end{itemize}
        \item If $(\xi, b_0)$ is an odd admissible pair, it follows that   $\xi=0$, $b_0=\alpha e_{1}+\beta e_2$, and $\om=e_1^*\wedge e^*_2$, where $\alpha,\beta\in \mathbb{K}$. Eq  \eqref{Liebrackets2}implies that there exists a basis $\{e_{1},e_2\mid e, d\}$  of $\mathfrak{g}$ such that the Lie brackets and the form are given by 
		\begin{equation*}
[d,d]= 2\al e_1+2\beta e_1,\;\;	\left[d,e_1\right]=-\beta e,\;\left[d,e_2\right]=\alpha e, \quad \om= -e^*\wedge d^*+e_1^* \wedge e_2^*.
\end{equation*} 
Since $(\g, \br)$ is non-abelian, then $(\alpha,\beta)\neq(0,0)$. If for instance $\beta\neq0$, then we put  $$(x_1, x_2, y_1,y_2)=\left (\frac{1}{\beta}e_1, 2\beta e_2+2\alpha e_1, d,  e\right ).$$ 
		Thus, $(\mathfrak{g},[\; ,\;], \om)$ is symplectomorphic to $(\g^2, \br, \om)$.
\end{itemize}

\underline{The case where  $\omega_{\mathfrak{b}}$ is odd.} We have several cases:
\begin{itemize}
\item If $(\xi, b_0)$ is an even admissible pair, then   $\xi=0$ and $b_0=\alpha e_{1}$, and $\om=e_1^*\wedge f^*_1$, where $\alpha\in \mathbb{K}$. Eq.  \eqref{Liebrackets3} implies that there exists a basis $\{d,e_{1}\mid e, f_1\}$  of $\mathfrak{g}$ such that the Lie brackets are given by 
		\begin{equation*}
	\left[d,f_1\right]=\alpha e, \quad \om= e^*\wedge d^*+e_1^* \wedge f_1^*.
		\end{equation*} 
		Since $(\g, \br)$ is non-abelian, then $\alpha\neq 0$. We put  $$(x_1,x_2, y_1,y_2)=\left (d, -\al e_1,- \frac{1}{\al}f_1, -e \right ).$$ 
		Thus,  $(\mathfrak{g},[\; , \;], \om)$ is symplectomorphic to $(\g^3, \br, \om)$.

\item If $(\xi, b_0)$ is an odd admissible pair, we distinguish several cases:
\begin{itemize} \item 
If $\xi=0$,  $b_0=\alpha e_{1}$, and $\om=e_1^*\wedge e^*_2$, where $\alpha\in \mathbb{K}$,  then  Eq.\eqref{Liebrackets4}, implies that there exists a basis $\{e,e_{1}\mid d, f_1\}$ of $\mathfrak{g}$ such that the Lie brackets are given by 
		\begin{equation*}
	\left[d,d\right]=-2\al e_1,\quad \left[d,f_1\right]=\alpha e, \quad \om= e^*\wedge d^*+e_1^* \wedge f_1^*.
		\end{equation*} 
		Since $(\g, \br)$ is non-abelian, then $\alpha\neq 0$. We put  $$(x_1,x_2,  y_1,y_2)=\left (-2\al e_1,  e,d, \frac{1}{\al}f_1, \right ).$$ 
		Thus,  $(\mathfrak{g}, [\;,\;],\om)$ is symplectomorphic to $(\g^4, \br, \om)$.
    \item If $\xi\neq 0$,  $b_0=0$, and $\om=e_1^*\wedge e^*_2$,   then  Eq.\eqref{Liebrackets4}, implies that there exists a basis $\{e,e_{1}\mid d,f_1\}$  of $\mathfrak{g}$ such that the Lie brackets are given by 
		\begin{equation*}
	 \left[d,e_1\right]=a  f_1, \quad \om= e^*\wedge d^*+e_1^* \wedge f_1^*.
		\end{equation*} 
		Put  $$(x_1,x_2,   y_1,y_2)=\left (e_1,-a  e,-\frac{1}{a}d, f_1, \right ).$$ 
		Thus,  $(\mathfrak{g},[\; ,\;], \om)$ is symplectomorphic to $(\g^3, \br, \om)$. 

\item If $\xi\neq 0$,  $b_0=\al e_1$, and $\om=e_1^*\wedge e^*_2$, where $\al\in \mathbb{K}$   then  Eq.\eqref{Liebrackets4}, implies that there exists a basis $\{e,e_{1}\mid d,f_1\}$  of $\mathfrak{g}$ such that the Lie brackets are given by 
		\begin{equation*}
	 \left[d,d\right]=-2\al e,\quad  \left[d,f_1\right]=a  e_1+\al e_0,\quad  \left[f_1,f_1\right]=-2 a e  \quad \om= e^*\wedge d^*+e_1^* \wedge f_1^*.
		\end{equation*} 
		 If $\al=0$, then we  put  $$ (x_1,x_2,  y_1,y_2)=\left (-\frac{2}{a} e,a e_1,\frac{1}{a}f_1, a d \right ).$$ 
		But if $\alpha\neq 0$, put
\[
(x_1,x_2, y_1,y_2)=
\left(
-\lambda  e_1
-
\mu e,\;
-\lambda \, e_1
+
\mu\, e,\;
\frac{1}{2}\big(\lambda^{-1}\, d
-
\mu^{-1}\, f_1\big),\;
\lambda^{-1}\, d
+
\mu^{-1}\, f_1
\right),
\]
where $\lambda=\frac{\sqrt[3]{\alpha^{2}}}{\sqrt[3]{a}}$ and $\mu=\frac{\sqrt[3]{a^{2}}}{\sqrt[3]{\alpha}}$. 

Thus, $(\mathfrak{g},[\;, \;], \om)$ is symplectomorphic to $(\g^4, \br, \om)$.
        \end{itemize}
   \qed     
\end{itemize}
\noqed
\end{proof}
\sssbegin{Remark}
The Lie superalgebras $\fg^2$, $\fg^3$, and $\fg^4$ correspond to $C^3+A$, $\mathbb{K}\oplus C^3$, and $(2A_{1,1} + 2A)_{\frac12}^3$ of \cite{B}, respectively. See also \cite{BM, BR} for more details. 
\end{Remark}
\subsection{Five-dimensional flat quasi-Frobenius Lie superalgebras}
Here, we classify flat quasi-Frobenius Lie superalgebras of total dimension five. Nilpotent Lie superalgebras of total dimension 5 over the reals were classified by Hegazi in \cite{H}.

\sssbegin{Proposition}\label{solution admissible dim3}
Let $(\mathfrak{b},[\;,\;]_\fb,  \omega_{\mathfrak{b}})$ be the three-dimensional orthosymplectic quasi-Frobenius abelian Lie superalgebra.
\begin{itemize}
    \item[(i)] The pair $(\xi, b_0)$ is \emph{even admissible} if and only if  there exists a basis $\mathbb{B} = \{e_1, e_2\mid f_1\}$ of $\mathfrak{b}$ such that  the following conditions are satisfied
    \begin{itemize}
        \item Either 
\[
\xi = 0, \quad b_0 = \alpha e_1 + \beta e_2, \quad \omega = e_1^* \wedge e_2^*+ f_1^* \wedge f_1^*, \quad \alpha, \beta \in \mathbb{K}.
\]
\item Or\[
\xi =
\begin{pmatrix}
0 & a & 0 \\
0 & 0 & 0 \\
0 & 0 & 0
\end{pmatrix},\; b_0=\al e_1,\, \omega = e_1^* \wedge e_2^*+f_1^* \wedge f_1^*, 
\quad 0\not=a, \al\in \mathbb{K}.
\]\end{itemize}
\begin{enumerate}
\item[(ii)] The pair $(\xi, b_0)$ is \emph{odd admissible} if and only if there exists a basis $\mathbb{B} = \{e_1, e_2 \mid f_1\}$ of $\mathfrak{b}$  such  that 
\[
\xi = 0, \quad b_0 = \alpha e_1 + \beta e_2, \quad \omega = e_1^* \wedge e_2^*+ f_1^* \wedge f_1^*, \quad \alpha, \beta \in \mathbb{K}.
\]

\end{enumerate}

\end{itemize}
\end{Proposition}

\begin{proof}
We distinguish two cases. 
 
\noindent\underline{The case where the pair $(\xi, b_0)$ is even admissible.} It follows that 
\[
\xi^* \circ \xi = 0, \qquad \xi \circ \xi^* = \xi^2, \qquad \text{and} \qquad b_0 \in \ker(\xi^* - \xi).
\]
Hence, $\xi$ is nilpotent, and therefore either $\xi = 0$ or $\xi \neq 0$.  Since $\omega_{\mathfrak{b}}$ is even, then  $\dim(\mathfrak{b}_{\bar{1}})=1$. Because $\xi$ is even and nilpotent, we have $\xi|_{\mathfrak{b}_{\bar{1}}} = 0$.   Choose a basis $\{e_1, e_2 \mid f_1\}$ of $\mathfrak{g}$, such that 
\[
\omega_{\mathfrak{b}} = e_1^* \wedge e_2^* +f_1^* \wedge f_1^*.
\]
Since $b_0$ is even, we can write
$
b_0 = \alpha e_1 + \beta e_2$, where $   \alpha, \beta \in \mathbb{K}.$

Finally, since $\xi$ is even and nilpotent, its matrix with respect to $\{e_1, e_2 \mid f_1\}$ is
\[
\xi =
\begin{pmatrix}
0 & a & 0 \\
0 & 0 & 0 \\
0 & 0 & 0
\end{pmatrix},
\qquad a \neq 0.
\]
Since $b_0 \in \ker(\xi^* - \xi)$, we deduce that $\beta = 0$.

\noindent\underline{The case where the pair $(\xi, b_0)$ is odd admissible.} It follows that 
\[
\xi^* \circ \xi = -\xi \circ \xi^*,\;  \xi^2 = 0, 
\qquad \text{and} \qquad b_0 \in \ker(2\xi+\xi^*).
\]
Hence, $\xi$ is nilpotent, and therefore either $\xi = 0$ or $\xi \neq 0$ with $\xi^2 = 0$.  
Since $\omega_{\mathfrak{b}}$ is even, then $\dim(\mathfrak{b}_{\bar{1}})=1$.  If $\xi|_{\mathfrak{b}_{\bar{1}}} = 0$, choose a basis $\{e_1, e_2 \mid f_1\}$ such that 
\[
\omega_{\mathfrak{b}} = e_1^* \wedge e_2^*+f_1^* \wedge f_1^*.
\]
Since $b_0$ is even, we can write
$
b_0 = \alpha e_1 + \beta e_2$, $ \alpha, \beta \in \mathbb{K}.
$
Then, the matrix of $\xi$ with respect to $\{e_1, e_2 \mid f_1\}$ is
\[
\xi =
\begin{pmatrix}
0 & 0 & 0 \\
0 & 0 & 0 \\
a & 0 & 0
\end{pmatrix},
\qquad a \in \mathbb{K}.
\]
Since $\xi^*\circ \xi=-\xi\circ \xi^*$, it follows that  $a=0$ and hence $\xi=0$.

Now, if $\xi|_{\mathfrak{b}_{\bar{1}}} \neq 0$, choose a basis $\{e_1, e_2 \mid f_1\}$ such that 
\[
\omega_{\mathfrak{b}} = e_1^* \wedge e_2^* +f_1^* \wedge f_1^*.
\]
Since $b_0$ is even, we can write
$
b_0 = \alpha e_1 + \beta e_2$,     $ \alpha, \beta \in \mathbb{K}$.

Then, the matrix of $\xi$ with respect to $\{e_1, e_2 \mid f_1\}$ is
\[
\xi =
\begin{pmatrix}
0 & 0 & a \\
0 & 0 & 0 \\
\alpha & \beta & 0
\end{pmatrix},
\qquad a \neq 0, \quad \alpha, \beta\in \mathbb{K}.
\]
Since $\xi^2 = 0$, we deduce that $\alpha = \beta = 0$. Moreover, from $\xi^* \circ \xi =- \xi \circ \xi^* $, it follows that $a = 0$, and hence $\xi = 0$, which is a contradiction.
\end{proof}
\sssbegin{Proposition}
Any five-dimensional flat quasi-Frobenius non-abelian Lie superalgebra is isomorphic to the following:

\begin{itemize}
\item[$(i)$]  $(\g^1, \br, \om)$ \textup{(}in the basis $\{x_1, x_2, x_3, x_3 \, | \,y_1\}$\textup{)}:
$$[x_1, x_2]=x_3, \quad \om= x_1^*\wedge x_4^*+x_2^* \wedge x_3^*+y_1^* \wedge y_1^*,$$
\item[$(ii)$]  $(\g^2, \br, \om)$ \textup{(}in the basis $\{x_1, x_2 \, | \,y_1, y_2, y_3\}$\textup{)}:
$$
	[x_1,y_1]=y_2,\quad [y_1, y_1]= x_2\mbox{ and }\om=2x_1^*\wedge x_2^*-y_1^*\wedge y_2^*+y_3^* \wedge y_3^*.
		$$
\end{itemize}
\end{Proposition}

\begin{proof}
Let $(\mathfrak{g},[\; , \;], \om)$ be  a flat  orthosymplectic quasi-Frobenius non-abelian Lie superalgebra of total dimension $5$. According to Theorem~\ref{central}, such a superalgebra is an odd or even flat double extension of the three-dimensional abelian quasi-Frobenius Lie superalgebra $(\mathfrak{b},[\;,\;]_\fb,  \omega_{\mathfrak{b}})$ by means of  $(\xi, b_0)$. By Prop. \ref{solution admissible dim3}, we have:  

\underline{The case where $(\xi, b_0)$ is an even admissible pair}. We have the following cases: 
\begin{itemize}
\item[$(i)$]  If $\xi=0$ and $b_0=\alpha e_{1}+\beta e_2$, and $\om=e_1^*\wedge e^*_2+f_1^*\wedge f^*_1$, where $\alpha,\beta\in \mathbb{K}$,  then according to \eqref{Liebrackets1}, there exists a basis $\{e,d,e_{1},e_2\mid f_1\}$ of $\mathfrak{g}$ such that the Lie brackets are given by 
		\begin{equation*}
	\left[d,e_1\right]=-\beta e,\;\left[d,e_2\right]=\alpha e, \quad \om= e^*\wedge d^*+e_1^* \wedge e_2^*+f_1^*\wedge f^*_1.
		\end{equation*} 
		Since $(\g, \br)$ is non-abelian, then $(\alpha,\beta)\neq(0,0)$. If for instance  $\beta\neq0$, then we put  $$(x_1,x_2,x_3,x_4, y_1)=\left (-\frac{1}{\beta}e_1, d, -e, -\beta e_2-\alpha e_1, f_1 \right ).$$ 
		Thus $(\mathfrak{g},[\;,\;], \om)$ is symplectomorphic to $(\g^1, \br, \om)$.
    \item[$(ii)$]  If $\xi\neq 0$ and $b_0=\alpha e_{1}$ where $\alpha\in \mathbb{K}$,  then according to \eqref{Liebrackets1}, there exists a basis $\{e,d,e_{1},e_2\mid f_1\}$ of $\mathfrak{g}$ such that the Lie brackets are given by 
$$\left[d,e_2\right]=\alpha e-3ae_1, \,\, a\neq 0.$$
		Put  $$(x_1,x_2,x_3,x_4, y_5)=\left(3a d+\alpha e_2,d+\frac{\alpha+1}{3a}e_2, \alpha e-3ae_1,-\frac{\alpha+1}{3a}e+e_1, f_1\right).$$ 
		Thus, in this new basis, the Lie brackets and the  form are reduced to
		\[
		[x_1,x_2]=x_3\quad \text{ and } \quad \om=x_1^*\wedge x_4^*+x_2^*\wedge x_3^*+y_1^*\wedge y_1^*,
		\]
		which implies that $(\mathfrak{g},[\;, \;], \om)$ is symplectomorphic to $(\g^1, \br, \om)$.\\
\end{itemize} 
\underline{The case where $(\xi, b_0)$ is an odd admissible pair.}  It follows that $\xi=0$ and $b_0=\alpha e_{1}+\beta e_2$, and $\om=e_1^*\wedge e^*_2+ f_1^*\wedge f^*_1$,  where $\alpha,\beta\in \mathbb{K}$,  then according to \eqref{Liebrackets2}, there exists a basis $\{ e_{1},e_2\mid e,d,f_1\}$ of $\mathfrak{g}$ such that the Lie brackets are given by 
		\begin{equation*}
[d,d]= 2\al e_1+2\beta e_1,\;\;	\left[d,e_1\right]=-\beta e,\;\left[d,e_2\right]=\alpha e, \quad \om= -e^*\wedge d^*+e_1^* \wedge e_2^*+f_1^*\wedge f^*_1.
\end{equation*} 
Since $(\g, \br)$ is non-abelian, then $(\alpha,\beta)\neq(0,0)$. If for instance  $\beta\neq0$, then we put  $$(x_1, x_2\mid y_1,y_2,y_3)=\left (\frac{1}{\beta}e_1, 2\beta e_2+2\alpha e_1, d,  e, f_1 \right ).$$ 
		Thus, $(\mathfrak{g},[\; , \;], \om)$ is symplectomorphic to $(\g^2, \br, \om)$.
\end{proof}

\bigskip

\noindent {\bf Acknowledgments:} We sincerely thank the referee for the valuable comments, which have improved the clarity and readability of the manuscript.

\bigskip

\noindent \underline{{\bf Declarations:}}
\bigskip

\noindent {\bf FUNDING STATEMENT}:  SB's research was supported by GRANT: AD-065. \\
{\bf COMPETING INTEREST  STATEMENT}: the authors declare no competing interest.


\end{document}